\documentclass[leqno,10pt]{amsart}

\usepackage{amssymb}
\usepackage{amsmath}
\usepackage{amsthm}
\usepackage{verbatim}

\usepackage{color}
\usepackage{geometry}
\renewcommand{\div}{\diver}
\newcommand{\loc}{\mathrm{loc}}
\newcommand{\Lip}{\mathrm{Lip}}
\newcommand{\metric}{\langle,\rangle}
\newcommand{\RR}{\mathbb R}
\newcommand{\PP}{\mathbb P}

\newcommand{\Mbar}{\overline{M}}
\newcommand{\gbar}{\overline{g}}
\newcommand{\nablabar}{\overline{\nabla}}

\newcommand{\ebar}{\overline{e}}
\newcommand{\Nbar}{\overline{N}}

\newcommand{\DeltaP}{\mathrm{\Delta}^{\PP}}
\newcommand{\Deltahat}{\widehat{\mathrm{\Delta}}}

\newcommand{\dist}{\gamma}
\newcommand{\disthat}{\widehat{\gamma}}
\newcommand{\distP}{\gamma^{\PP}}

\newcommand{\BP}{B^{\PP}}
\newcommand{\Bhat}{\widehat{B}}

\newcommand{\Riembar}{\overline{\Riem}}
\newcommand{\Rbar}{\overline{R}}
\newcommand{\RiemP}{\Riem^{\PP}}
\newcommand{\RP}{R^{\PP}}

\newcommand{\Ricbar}{\overline{\Ric}}
\newcommand{\RicP}{\Ric^{\PP}}
\newcommand{\Richat}{\widehat{\Ric}}

\newcommand{\HessP}{\Hess^{\PP}}
\newcommand{\Hesshat}{\widehat{\Hess}}

\newcommand{\kulknom}{\mathbin{\bigcirc\mspace{-15mu}\wedge\mspace{3mu}}}
\newcommand{\Hmean}{\widehat{H}}
\newcommand{\II}{\mathrm{I\hspace{-1pt}I}}
\newcommand{\Hvec}{\mathbf{H}}

\newcommand{\fgrad}{z}

\newcommand{\diver}{\mathrm{div}}
\newcommand{\vol}{\mathrm{vol}}
\DeclareMathOperator{\supp}{supp}
\DeclareMathOperator{\essinf}{ess\,inf}
\DeclareMathOperator{\esssup}{ess\,sup}

\newcommand{\trace}{\mathrm{trace}}
\newcommand{\Ric}{\mathrm{Ric}}
\newcommand{\Hess}{\mathrm{Hess}}

\newcommand{\cut}{\mathrm{cut}}
\newcommand{\Riem}{\mathrm{Riem}}
\newcommand{\di}{\mathrm{d}}

\theoremstyle{plain}
\newtheorem{thm}{Theorem}
\newtheorem{crl}{Corollary}
\newtheorem{prp}{Proposition}
\newtheorem{lemma}{Lemma}

\numberwithin{equation}{section}

\theoremstyle{definition}

\newtheorem{rmk}{Remark}

\title{Spacelike hypersurfaces in standard static spacetimes}

\author[G. Colombo]{Giulio Colombo}
\address{Dipartimento di Matematica, Universit\`a degli Studi di Milano, 20133 Milano, Italy}
\curraddr{}
\email{giulio.colombo@unimi.it}
\thanks{}

\author[J. A. S. Pelegr\'in]{Jos\'e A. S. Pelegr\'in}
\address{Departamento de Geometr\'ia y Topolog\'ia, Universidad de Granada, 18071 Granada, Spain}
\curraddr{Departamento de Matem\'atica Aplicada y Estad\'istica, Universidad 
	CEU San Pablo, 28003 Madrid, Spain}
\email{jpelegrin@ugr.es, jose.sanchezpelegrin@ceu.es}
\thanks{}

\author[M. Rigoli]{Marco Rigoli}
\address{Dipartimento di Matematica, Universit\`a degli Studi di Milano, 20133 Milano, Italy}
\curraddr{}
\email{marco.rigoli@unimi.it}
\thanks{}

\begin{document}

\maketitle

\begin{abstract}
In this work we study spacelike hypersurfaces immersed in spatially open standard static spacetimes with complete spacelike slices. Under appropriate lower bounds on the Ricci curvature of the spacetime in directions tangent to the slices, we prove that every complete CMC hypersurface having either bounded hyperbolic angle or bounded height is maximal. Our conclusions follow from general mean curvature estimates for spacelike hypersurfaces. In case where the spacetime is a Lorentzian product with spatial factor of nonnegative Ricci curvature and sectional curvatures bounded below, we also show that a complete maximal hypersurface not intersecting a spacelike slice is itself a slice. This result is obtained from a gradient estimate for parametric maximal hypersurfaces.
\end{abstract}

\bigskip



\section{Introduction and main results}

When searching for general solutions of Einstein's field equations in a spacetime $\Mbar$, it is customary to assume the a priori existence of an infinitesimal symmetry (see \cite{Da}, \cite{Ea}). The symmetry often comes from a Killing or, more generally, a conformal vector field $X$ on $\Mbar$, see for instance \cite{Du}. In this case, the spacetime can be classified depending on the causal character of the symmetry. Thus, a Lorentzian manifold admitting a timelike Killing vector field $X$ is called a stationary spacetime due to the fact that observers along the vector field $X$ see a metric that does not change. Moreover, if the timelike Killing vector field is irrotational, that is, the distribution orthogonal to the field is involutive, then a local warped product structure appears and the spacetime is called static \cite{BEE}. When this structure is global the spacetime is called a standard static spacetime. 

More precisely, by a standard static spacetime $(\Mbar,\gbar)$ we mean a product $\Mbar = \PP\times \RR$, with $\PP$ a (connected) orientable manifold of dimension $m\geq 2$, endowed with the Lorentzian metric
\begin{equation} \label{me}
\gbar = \pi_{\PP}^{\ast}(\sigma) - (h\circ\pi_{\PP})^2 \pi_{\RR}^{\ast}(\di t^2),
\end{equation}
where $\pi_{\PP} : \Mbar \to \PP$, $\pi_{\RR} : \Mbar \to \RR$ are the projections onto the factors of the product, and where $\sigma$ and $h$ are respectively a Riemannian metric and a smooth positive function on $\PP$. Thus, the spacetime $\Mbar$ is a warped product in the sense of \cite{O'N}, with base $(\PP,\sigma)$, fiber $(\RR,-\di t^2)$ and warping function $h$. Each tangent vector $X\in T_p\Mbar$, $p=(x,t)\in\Mbar$, can be decomposed as $X=\mathcal{H}X+\mathcal{V}X$ with $\mathcal{H}X$ tangent to the leaf $\PP\times\{t\}$ and $\mathcal{V}X$ tangent to the fiber $\{x\}\times\RR$, that is, $\mathcal{H}X\in\ker(\pi_{\RR})_{\ast}$ and $\mathcal{V}X\in\ker(\pi_{\PP})_{\ast}$. Following \cite{O'N}, vectors $X=\mathcal{H}X$ tangent to leaves will be called horizontal and vectors $X=\mathcal{V}X$ tangent to fibers will be called vertical. Note that a horizontal vector $X$ is always spacelike so for such $X$ we can set $|X|=\sqrt{\gbar(X,X)}\geq0$.

When $h\equiv1$ the resulting standard static spacetime $(\Mbar,\gbar)$ is a semi-Riemannian product with factors $(\PP,\sigma)$ and $(\mathbb R,-\di t^2)$; in the following, a manifold of this type will be called a Lorentzian product.

On a standard static spacetime, the role of the vector field $X$ responsible of the infinitesimal symmetry is played by the global timelike Killing vector field $\partial_t$ on $\Mbar$. In \cite{O'N} it is proved that any static spacetime is locally isometric to a standard static one. Moreover, in \cite{ARR} and \cite{Sa} sufficient conditions are given for a static spacetime to be standard. The importance of standard static spacetimes is due to the fact that they include some classical spacetimes like Lorentz-Minkowski spacetime $\mathbb{L}^n$, Einstein static universe as well as models that describe a universe with one spherically symmetric non-rotating mass, such as a star or a black hole, as it happens in the exterior Schwarzschild spacetime \cite{SW}.

In any standard static spacetime $\Mbar = \PP\times_h \RR$ there exists a distinguished foliation whose leaves are given by the totally geodesic level hypersurfaces of the function $\pi_{\RR}$. They are known as the spacelike slices $\PP\times\{t_0\}$, $t_0\in\RR$. We recall that in the spatially closed case, several uniqueness results have been obtained on the splitting of these spacetimes in terms of their usual orthogonal decomposition (see \cite{ARR} and \cite{SS}). However, the problem of guaranteeing the uniqueness of the splitting for spatially open standard static spacetimes remains open, in fact there exist spacetimes with different splittings of type (\ref{me}), for instance $\mathbb{L}^n$. Other nontrivial cases are considered in \cite{Tod} and in \cite{GO}, where the authors describe the general structure of a standard static spacetime admitting more than one decomposition and give some uniqueness results under suitable curvature assumptions.

In this paper we focus our attention on spatially open standard static spacetimes. Their importance comes from the fact that, despite the historical relevance of spatially closed models, observations suggest that our physical universe is actually spatially open \cite{Chiu}. Moreover, spatially closed spacetimes lead to a  violation of the holographic principle \cite{BR}, making spatially open models more suitable for a possible quantum theory of gravity \cite{Bo}.

Given an $m$-dimensional manifold $M$, an immersion
\[
\psi : M \to \Mbar
\]
in a standard static spacetime $\Mbar = \PP\times_h\RR$ is said to be spacelike if $g=\psi^{\ast}\gbar$ is a Riemannian metric on $M$. In this case, $M$ is called a spacelike hypersurface and $\psi:(M,g)\to(\Mbar,\gbar)$ is an isometric immersion. Roughly speaking, each spacelike hypersurface represents the physical space of some observer in a given instant of their time and their study has been crucial in General Relativity \cite{MT}. Among other reasons, their interest rely on the key role they play in the proof of the positivity of the gravitational mass \cite{SY}, their importance in the study of the structure of singularities in the space of solutions of Einstein's equations \cite{AMM} and the fact that the initial value problem for the Einstein's field equation in General Relativity is formulated in terms of a spacelike hypersurface (see, for instance \cite{Ri} and references therein). Furthermore, in Causality Theory, the existence of a certain spacelike hypersurface can determine the causal properties of the spacetime. For instance, a spacetime is globally hyperbolic if and only if it admits a Cauchy hypersurface \cite{Ger}. Indeed, any globally hyperbolic spacetime is diffeomorphic to $\mathbb{R} \times S$, being $S$ a smooth spacelike Cauchy hypersurface, see Theorem 1 of \cite{BS}. 

\bigskip

\noindent Let $\psi:M\to\Mbar$ be a spacelike hypersurface immersed in the standard static spacetime $\Mbar=\PP\times_h\RR$. The map
\[
\pi = \pi_{\psi} = \pi_{\PP}\circ\psi : M\to\PP
\]
is an immersion and therefore a local diffeomorphism. The tensor field $\hat\sigma=\pi^{\ast}\sigma$ is a Riemannian metric on $M$ that satisfies $\hat\sigma\geq g$ in the sense of quadratic forms, where $g=\psi^{\ast}\gbar$ is the Riemannian metric induced by the immersion $\psi$.

If $\PP$ is noncompact then $M$ itself is noncompact, since $\pi$ is continuous and open. If $(\PP,\sigma)$ is complete then $(M,\hat\sigma)$ is complete if and only if $\pi:M\to\PP$ is a (topological) covering map. Vice versa, if $(M,\hat\sigma)$ is complete, then $\pi:(M,\hat\sigma)\to(\PP,\sigma)$ is a Riemannian covering map and $(\PP,\sigma)$ itself is complete. For a justification of these statements we refer the reader to Proposition 23 and Lemma 17 of \cite{Pe}. We remark that completeness of $(M,g)$ implies completeness of $(M,\hat\sigma)$, since $\hat\sigma\geq g$. We also remark that if $\pi:M\to\PP$ is a covering map of degree 1, that is, a global diffeomorphism, then $\psi(M)\subseteq\PP\times\RR$ is the graph of a smooth function $u:\PP\to\RR$. Viceversa, the graph of a smooth function $u:\PP\to\RR$ satisfying
\[
h^2 \sigma(Du,Du) < 1 \qquad \text{on } \PP
\]
is the image of the spacelike hypersurface $\psi_u:\PP\to\Mbar$ defined by $\psi_u(x)=(x,u(x))$ for each $x\in\PP$, and the corresponding map $\pi_u=\pi_{\PP}\circ\psi_u:\PP\to\PP$ is the identity map.

The manifold $M$ is orientable as a consequence of the time-orientability of $\Mbar$ ensured by the presence of the global timelike vector field $\partial_t$. In particular, there exists a unique unit timelike normal vector field $N\in\mathfrak X^{\bot}(M)$ with the same time-orientation of $\partial_t$, that is, satisfying $\gbar(N,\partial_t)<0$ everywhere. The wrong way Cauchy inequality implies that
\[
\gbar\left(N,\frac{1}{h\circ\pi}\,\partial_t\right)\leq-1.
\]
This enables us to define the hyperbolic cosine of the hyperbolic angle $\theta$ between $N$ and $\partial_t$ by setting
\[
\cosh\theta = -\gbar\left(N,\frac{1}{h\circ\pi}\,\partial_t\right).
\]
The mean curvature function of $\psi$ in the direction of $N$ will be denoted by $H$. Both $\cosh\theta$ and $H$ are smooth functions on $M$.

\bigskip

In Section 2 below we derive the geometric and analytical equations relevant for our pourposes and in Section 3 we obtain, as first result, some lower bounds on the hyperbolic angle of a spacelike hypersurface under different geometric assumptions. Our first result is the following

\begin{thm} \label{intro:thm:cosh1}
	Let $\Mbar=\PP\times_h\RR$ be a standard static spacetime with complete, noncompact base $(\PP,\sigma)$. Suppose that for some constant $G_0>0$ it holds
	\begin{equation} \label{QBEriccibound}
	\Ricbar(X,X) \geq - mG_0|X|^2
	\end{equation}
	for each horizontal vector $X\in T\Mbar$. Let $\psi:M\to\Mbar$ be a spacelike hypersurface such that $\pi:M\to\PP$ is a covering map. If $|H|\geq H_0$ on $M$ for some $H_0>0$, then
	\begin{equation} \label{limsupcosh}
	\limsup_{M\ni x\to\infty} \cosh\theta(x) \geq \sqrt{1+\frac{H_0^2}{G_0}}.
	\end{equation}
\end{thm}

\noindent As a consequence, but see also Corollary \ref{intro:crl:cosh2} below, we have

\begin{crl} \label{intro:crl:cosh1}
	Let $\Mbar=\PP\times_h\RR$ be a standard static spacetime with complete, noncompact base $(\PP,\sigma)$. Suppose that
	\begin{equation} \label{QBEriccinonneg}
	\Ricbar(X,X) \geq 0
	\end{equation}
	for each horizontal vector $X\in T\Mbar$. If $\psi:M\to\Mbar$ is a spacelike immersed hypersurface such that $\pi:M\to\PP$ is a covering map and the hyperbolic angle is bounded, then
	\[
	\inf_M |H| = 0.
	\]
	In particular, if $\psi$ has constant mean curvature then it is maximal.
\end{crl}

\begin{rmk} \label{rmk:ricci}
	As a consequence of Lemma \ref{curvaturesbar} below we have that for each horizontal vector $X\in T\Mbar$ it holds
	\[
	\Ricbar(X,X) = \RicP(X_0,X_0) - \frac{\HessP(h)(X_0,X_0)}{h}, \quad |X|^2=\sigma(X_0,X_0)
	\]
	with $X_0=(\pi_{\PP})_{\ast}X\in T\PP$. Therefore, assumptions (\ref{QBEriccibound}) and (\ref{QBEriccinonneg}) can be regarded as generalized curvature assumptions on the weighted manifold $(\PP,\sigma,h)$, since $\RicP - 	\frac{\HessP(h)}{h}$ is the modified Bakry-Emery Ricci tensor
	\begin{align*}
	\RicP_{m,m+1}(\DeltaP_{-\log h})=\RicP + \HessP(-\log h) + \di(-\log h)\otimes \di(-\log h) = \RicP - \frac{\HessP(h)}{h}
	\end{align*}
	introduced by Qian, \cite{Qi}, and ubiquitous in the study of generalized $k$-Einstein manifolds, Ricci solitons, etc. See for instance \cite{MRS}, \cite{CMMR} and references therein. A lower bound of the form (\ref{QBEriccibound}) implies an upper bound on the growth of the weighted volume of geodesic balls in $\PP$, with weight function $h$, namely, inequality \eqref{volumeGbound} below (cfr.~\cite{MRS}). Indeed, a bound of the form (\ref{limsupcosh}) can be obtained when (\ref{QBEriccibound}) is directly replaced by a volume growth assumption on $\PP$, in the case where the covering map $\pi:M\to\PP$ has finite degree. This is the content of the next result. We recall that in the non-parametric case where $M=\PP$, $\psi=\psi_u$ and $\psi(M)=\psi_u(\PP)$ is an entire graph over $\PP$, the map $\pi_u$ has degree $1$.
\end{rmk}

\begin{thm} \label{intro:thm:cosh2}
	Let $\Mbar^{m+1}=\PP^m\times_h\RR$ be a standard static spacetime with complete, noncompact base $(\PP,\sigma)$. Let $G_0>0$ be a given constant and suppose that for some (hence, any) $q\in\PP$
	\begin{equation} \label{volumeGbound}
	\liminf_{r\to+\infty}\frac{\log\left(\int_{\BP_r}h\right)}{r} \leq m\sqrt{G_0},
	\end{equation}
	where $\BP_r=B^{\sigma}_r(q)$ is the geodesic ball of $(\PP,\sigma)$ centered at $q$ with radius $r$. Let $\psi:M^m\to\Mbar$ be a spacelike immersed hypersurface such that $\pi:M\to\PP$ is a covering map of finite degree. If $|H|\geq H_0$ on $M$ for some $H_0>0$, then
	\[
	\limsup_{M\ni x\to\infty} \cosh\theta(x) \geq \sqrt{1+\frac{H_0^2}{G_0}}.
	\]
\end{thm}

\begin{crl} \label{intro:crl:cosh2}
	Let $\Mbar^{m+1}=\PP^m\times_h\RR$ be a standard static spacetime with complete, noncompact base $(\PP,\sigma)$ such that for some (hence any) $q\in\PP$
	\[
	\liminf_{r\to+\infty}\frac{\log\left(\int_{\BP_r}h\right)}{r} = 0
	\]
	where $\BP_r=B^{\sigma}_r(q)$ is the geodesic ball of $(\PP,\sigma)$ centered at $q$ with radius $r$. Let $\psi:M^m\to\Mbar$ be a spacelike immersed hypersurface such that $\pi:M\to\PP$ is a covering map of finite degree. If the hyperbolic angle is bounded then
	\[
	\inf_M |H|=0.
	\]
	In particular, if $\psi$ has constant mean curvature then it is maximal.
\end{crl}

In Sections 4 and 5 we give several ``half-space'' theorems for spacelike hypersurfaces. In particular, in Section 4 we first focus on the case where $\Mbar$ is a Lorentzian product and we prove the next

\begin{thm} \label{intro:thm:hs2}
	Let $\Mbar^{m+1}=\PP^m\times\RR$ be a Lorentzian product with complete, noncompact base $(\PP,\sigma)$. Suppose that for some constant $G_0>0$ it holds
	\begin{equation} \label{riccibound'}
	\Ricbar(X,X) \geq - (m-1)G_0|X|^2
	\end{equation}
	for each horizontal vector $X\in T\Mbar$. Let $\psi:M\to\Mbar$ be a spacelike hypersurface such that $\pi:M\to\PP$ is a covering map.
	\begin{itemize}
		\item [(a)] If $H\geq0$ on $M$ and
		\[
		\liminf_{M\ni x\to\infty} H(x) > 0
		\]
		then $\psi(M)$ is not contained in any ``lower half-space'' of the form $\PP\times(-\infty,t_0]$, $t_0\in\RR$.
		\item [(b)] If $H\leq0$ on $M$ and
		\[
		\limsup_{M\ni x\to\infty} H(x) < 0
		\]
		then $\psi(M)$ is not contained in any ``upper half-space'' of the form $\PP\times[t_0,+\infty)$, $t_0\in\RR$.
	\end{itemize}
\end{thm}

\noindent The following is an immediate consequence of Theorem \ref{intro:thm:hs2}.

\begin{crl} \label{intro:crl:hs1-2}
	Let $\Mbar^{m+1}=\PP^m\times\RR$ be a Lorentzian product with complete, noncompact base $(\PP,\sigma)$. Suppose that for some constant $G_0>0$
	\begin{equation}
	\Ricbar(X,X) \geq - (m-1) G_0|X|^2
	\end{equation}
	for each horizontal vector $X\in T\Mbar$. Let $\psi:M\to\Mbar$ be a spacelike hypersurface such that $\pi:M\to\PP$ is a covering map. If $\psi(M)$ is contained in a slab $\PP\times[t_0,t_1]$, $-\infty<t_0<t_1<+\infty$, and $H$ does not change sign on $M$ then
	\[
	\liminf_{M\ni x\to\infty}|H(x)|=0.
	\]
	In particular, if $\psi$ has constant mean curvature then it is maximal.
\end{crl}

The argument that proves Theorem \ref{intro:thm:hs2} can be easily adapted to the case where $\Mbar$ is a standard static spacetime with radially symmetric base $(\PP,\sigma)$ and with warping factor given by a radial function $h$ on $(\PP,\sigma)$. These assumptions on the structure of $\Mbar$, although restrictive, are satisfied by several classical solutions of Einstein equations. As an example, we consider the case where $\Mbar$ is the Schwarzschild spacetime $\Mbar^{m+1} = \PP_S^m \times_{h_S} \RR$ with
\[
	(\PP_S,\sigma_S) = \left( (\rho_S,+\infty) \times \mathbb{S}^{m-1}, \frac{\di\rho^2}{V(\rho)} + \rho^2 \metric_{ \mathbb{S}^{m-1} } \right),
\]
\[
	V(\rho) = 1 - 2\mu\rho^{2-m}, \qquad h_S = \sqrt{V(\rho)}, 
\]
where $\rho$ is the standard coordinate on the interval $(\rho_S,+\infty)$, $\mu>0$ is a mass parameter and $\rho_S = (2\mu)^{1/(m-2)}$. In this setting, we prove the next
\begin{thm} \label{intro:thm:schwarz}
	Let $\Mbar = \PP_S\times_{h_S}\RR$ be the Schwarzschild spacetime of dimension $m+1$, with $m\geq 3$. Let $\psi:M\to\Mbar$ be a spacelike hypersurface such that $\pi:M\to\PP$ is a covering map.
	\begin{itemize}
		\item [(a)] If
		\[
			\liminf_{\substack{\rho(\pi(x)) \to +\infty \\ x \in M}} H(x) > 0
		\]
		then $\psi(M)$ is not contained in any lower half-space of the form $\PP\times(-\infty,t_0]$, $t_0\in\RR$.
		\item [(b)] If
		\[
			\limsup_{\substack{\rho(\pi(x)) \to +\infty \\ x \in M}} H(x) < 0
		\]
		then $\psi(M)$ is not contained in any upper half-space of the form $\PP\times[t_0,+\infty)$, $t_0\in\RR$.
	\end{itemize}
\end{thm}

In Section 5 we prove further results of the above type in the case where $\Mbar$ is a standard static spacetime, under different geometric assumptions on $\Mbar$ and $\psi$. In particular, when the hyperbolic angle of the hypersurface is bounded and  suitable bounds on the growth of $h$ and of the volume of $\PP$ are satisfied, we obtain conclusions similar to those in Theorem \ref{intro:thm:hs2}.

\begin{thm} \label{intro:thm:hs3}
	Let $\Mbar^{m+1}=\PP^m\times_h\RR$ be a standard static spacetime with complete, noncompact base $(\PP,\sigma)$ and let $\psi:M\to\Mbar$ be a spacelike hypersurface with bounded hyperbolic angle and such that $\pi:M\to\PP$ is a covering map. Assume that one of the following conditions is satisfied:
	\begin{itemize}
		\item [(i)] for some constant $G_0>0$
		\[
		\Ricbar(X,X)\geq-mG_0|X|^2
		\]
		for each horizontal vector $X\in T\Mbar$, and for some (hence, any) point $q\in\PP$
		\[
		\limsup_{\PP\ni p\to\infty}\frac{h(p)}{d_{\sigma}(p,q)}<+\infty,
		\]
		where $d_{\sigma}$ is the distance on $\PP$ induced by the metric $\sigma$;
		\item [(ii)] $\pi$ has finite degree and for some $\mu\in[0,2)$ and for some (hence, any) point $q\in\PP$
		\begin{align*}
		\liminf_{r\to+\infty}\frac{\log\left(\int_{\BP_r}h\right)}{r^{2-\mu}}< & +\infty, \\
		\limsup_{\PP\ni p\to\infty}\frac{h(p)}{d_{\sigma}(p,q)^{\mu}}< & +\infty,
		\end{align*}
		where $\BP_r=B^{\sigma}_r(q)$ is the geodesic ball of $(\PP,\sigma)$ centered at $q$ with radius $r$.
	\end{itemize}
	Then:
	\begin{itemize}
		\item [(a)] if $H\geq0$ on $M$ and
		\[
		\liminf_{M\ni x\to\infty}H(x)>0
		\]
		then $\psi(M)$ is not contained in any lower half-space of the form $\PP\times(-\infty,t_0]$, $t_0\in\RR$;
		\item [(b)] if $H\leq0$ on $M$ and
		\[
		\limsup_{M\ni x\to\infty}H(x)<0
		\]
		then $\psi(M)$ is not contained in any upper half-space of the form $\PP\times[t_0,+\infty)$, $t_0\in\RR$.
	\end{itemize}
\end{thm}

\noindent As a direct consequence of Theorem \ref{intro:thm:hs3} we obtain

\begin{crl} \label{intro:crl:hs3}
	Let $\Mbar=\PP\times_h\RR$ be a standard static spacetime with complete, noncompact base $(\PP,\sigma)$ and let $\psi:M\to\Mbar$ be a spacelike hypersurface with bounded hyperbolic angle and such that $\pi:M\to\PP$ is a covering map. Assume that either condition (i) or condition (ii) of Theorem \ref{intro:thm:hs3} is satisfied. If $\psi(M)$ is contained in a slab $\PP\times[t_0,t_1]$, $-\infty<t_0<t_1<+\infty$, and $H$ does not change sign on $M$ then
	\[
	\liminf_{M\ni x\to\infty}|H(x)|=0.
	\]
	In particular, if $\psi$ has constant mean curvature, then it is maximal.
\end{crl}

Clearly, the assumptions contained in Theorem \ref{intro:thm:hs3} on the function $h$ are satisfied at once (with $\mu=0$ in setting $(ii)$) when $h$ is a bounded function on $\PP$. In this case we can also replace the hypothesis that $\psi:M\to\Mbar$ has bounded hyperbolic angle with the assumption that $M$ is complete in the induced metric. This allows us to reach the same conclusions of Theorem \ref{intro:thm:hs3} under different geometric conditions.

\begin{thm} \label{intro:thm:hs4}
	Let $\Mbar=\PP\times_h\RR$ be a standard static spacetime with complete, noncompact base $(\PP,\sigma)$ and let $\psi:M\to\Mbar$ be a spacelike complete hypersurface. Assume that $h$ is bounded and that for some (hence, any) $o\in M$
	\begin{equation} \label{linfi}
	\liminf_{r\to+\infty}\frac{\log(\vol(B_r))}{r^2} < +\infty,
	\end{equation}
	where $B_r=B^g_r(o)$ is the geodesic ball of $(M,g)$ centered at $o$ with radius $r$ and the volume is measured with respect to the induced volume element on $M$.
	\begin{itemize}
		\item [(a)] If $H\geq0$ on $M$ and
		\[
		\liminf_{M\ni x\to\infty}H(x)>0
		\]
		then $\psi(M)$ is not contained in any lower half-space of the form $\PP\times(-\infty,t_0]$, $t_0\in\RR$;
		\item [(b)] if $H\leq0$ on $M$ and
		\[
		\limsup_{M\ni x\to\infty}H(x)<0
		\]
		then $\psi(M)$ is not contained in any upper half-space of the form $\PP\times[t_0,+\infty)$, $t_0\in\RR$.
	\end{itemize}
\end{thm}

\noindent We then deduce

\begin{crl} \label{intro:crl:hs4}
	Let $\Mbar=\PP\times_h\RR$ be a standard static spacetime with complete, noncompact base $(\PP,\sigma)$ and let $\psi:M\to\Mbar$ be a spacelike complete hypersurface. Assume that $h$ is bounded and that for some (hence, any) $o\in M$ condition (\ref{linfi}) is satisfied. If $\psi(M)$ is contained in a slab $\PP\times[t_0,t_1]$, $-\infty<t_0<t_1<+\infty$, and $H$ does not change sign on $M$ then
	\[
	\liminf_{M\ni x\to\infty}|H(x)|=0.
	\]
	In particular, if $\psi$ has constant mean curvature, then it is maximal.
\end{crl}

A different, more restrictive bound on the growth of the volume of $(M,g)$ forces the image $\psi(M)$ of the hypersurface to lie in a spacelike slice. Indeed, the following uniqueness result holds.

\begin{thm} \label{intro:thm:u1}
	Let $\Mbar^{m+1}=\PP^m\times_h\RR$ be a standard static spacetime with complete, noncompact base $(\PP,\sigma)$ and let $\psi:M\to\Mbar$ be a spacelike complete hypersurface. Assume that $h$ is bounded and that, for some $o\in M$,
	\begin{equation} \label{notl1}
	\int_R^{+\infty}\frac{\di r}{\vol(\partial B_r)} = +\infty
	\end{equation}
	for some (hence, any) $R>0$, where $\vol(\partial B_r)$ is, for a.~e.~$r\in\RR^+$, the Hausdorff $(m-1)$-dimensional measure of the boundary of the geodesic ball $B_r=B^g_r(o)$ of $(M,g)$ centered at $o$ with radius $r$.
	\begin{itemize}
		\item [(a)] If $H\geq0$ on $M$ then either $\psi(M)$ is a totally geodesic slice or it is not contained in any lower half-space of the form $\PP\times(-\infty,t_0]$, $t_0\in\RR$;
		\item [(b)] if $H\leq0$ on $M$ then either $\psi(M)$ is a totally geodesic slice or it is not contained in any upper half-space of the form $\PP\times[t_0,+\infty)$, $t_0\in\RR$.
	\end{itemize}
\end{thm}

\begin{crl} \label{intro:crl:u1}
	Let $\Mbar=\PP\times_h\RR$ be a standard static spacetime with complete, noncompact base $(\PP,\sigma)$ and let $\psi:M\to\Mbar$ be a spacelike complete hypersurface. Assume that $h$ is bounded and that, for some $o\in M$, condition (\ref{notl1}) is satisfied, and also assume that one of the following conditions is satisfied:
	\begin{itemize}
		\item [(a)] $\psi$ has constant mean curvature and $\psi(M)$ is contained in a slab $\PP\times[t_0,t_1]$, $-\infty<t_0<t_1<+\infty$;
		\item [(b)] $\psi$ is maximal and there exists some $t_0\in\RR$ such that $\psi(M)$ does not intersect the slice $\PP\times\{t_0\}$.
	\end{itemize}
	Then, $\psi(M)$ is a totally geodesic slice.
\end{crl}

The companion result to Theorem \ref{intro:thm:u1} in the non-parametric case is the following Calabi-Bernstein type result for spacelike graphs in standard static spacetimes. Condition \eqref{hnotl1} below can be regarded as an upper bound on the weighted volume of the boundaries of geodesic balls of $\PP$, with weight $h^2$. Note that when $h$ is bounded on $\PP$, \eqref{hnotl1} is certainly satisfied if an analogous condition is imposed on the nonweighted volume of such boundaries, as in \eqref{notl1}.

\begin{thm} \label{intro:thm:u2}
	Let $(\PP,\sigma)$ be a complete, noncompact Riemannian manifold and consider a positive function $h \in C^\infty (\PP)$ such that, for some $q\in\PP$,
	\begin{equation} \label{hnotl1}
	\int_R^{+\infty}\frac{\di r}{\int_{\partial\BP_r}h^2}=+\infty,
	\end{equation}
	for some (hence, any) $R>0$, where $\BP_r=B^{\sigma}_r(q)$ is the geodesic ball of $(\PP,\sigma)$ centered at $q$ with radius $r$. Let $H \in C^\infty(\PP)$ be a nonnegative function. Then, constant functions are the only entire bounded above solutions of the equation
	\begin{equation} \label{Hequation-sigma}
	\div\left(\frac{hDu}{\sqrt{1-h^2|Du|^2}}\right)+\frac{\sigma(Dh,Du)}{\sqrt{1-h^2|Du|^2}}=mH \quad \text{on } \PP
	\end{equation}
	satisfying $h|Du|<1$ pointwise on $\PP$ and
	\begin{equation}
	\int_R^{+\infty}\frac{\di r}{\int_{\partial \BP_r}\frac{h^2}{\sqrt{1-h^2|Du|^2}}}=+\infty.
	\end{equation}
	In particular, there exists such a solution if and only if $H \equiv 0$.
\end{thm}

We observe that in order to obtain the above results we have used appropriate forms of the comparison principle and of the weak maximum principle, valid under the assumptions of completeness and volume growth bounds like the one in (\ref{linfi}), or we have guaranteed a parabolic setting for the appropriate operator as in Theorem \ref{intro:thm:u1}.

\bigskip

In the last part of the paper we prove an upper bound for the hyperbolic angle of a maximal hypersurface immersed in a Lorentzian product $\Mbar$ whose base satisfies a uniform negative curvature bound from below, provided the image of the hypersurface is contained in a half-space of $\Mbar$.

\begin{thm} \label{intro:thm:coshb}
	Let $\Mbar^{m+1}=\PP^m\times\RR$ be a Lorentzian product with complete, noncompact base $(\PP,\sigma)$. Suppose that there exist two constants $G>0$, $B>0$ such that the Ricci curvature of $(\PP,\sigma)$ satisfies
	\begin{equation} \label{riccibound}
	\Ric \geq -(m-1)G\,\sigma
	\end{equation}
	on $\PP$ (in the sense of quadratic forms) and that the sectional curvatures of $(\PP,\sigma)$ are bounded from below by $-B$ on $\PP$, that is, for each $p\in\PP$ and for each $2$-plane $\Pi\subseteq T_p\PP$,
	\begin{equation} \label{sectbound}
	K_p(\Pi) \geq -B.
	\end{equation}
	Let $\psi:M\to\Mbar$ be a maximal spacelike hypersurface such that $\pi:M\to\PP$ is a covering map. If there exists $t_0\in\RR$ such that $\psi(M)$ does not intersect the slice $\PP\times\{t_0\}$, then
	\begin{equation} \label{habound}
	\cosh\theta \leq e^{(m-1)\sqrt{2G}|\tau-t_0|}
	\end{equation}
	on $M$, where $\tau=\pi_{\RR}\circ\psi$.
\end{thm}

\noindent As an immediate consequence of Theorem \ref{intro:thm:coshb} we obtain the following Calabi-Bernstein result for Lorentzian products.

\begin{crl} \label{intro:crl:u2}
	Let $\Mbar=\PP\times\RR$ be a Lorentzian product with complete, noncompact base $(\PP,\sigma)$. Suppose that the Ricci curvature of $(\PP,\sigma)$ is nonnegative and that the sectional curvatures of $(\PP,\sigma)$ are uniformly bounded from below on $\PP$. Let $\psi:M\to\Mbar$ be a maximal spacelike hypersurface such that $\pi:M\to\PP$ is a covering map. If there exists $t_0\in\RR$ such that $\psi(M)$ does not intersect the slice $\PP\times\{t_0\}$, then $\psi(M)$ itself is a slice $\PP\times\{t_1\}$ for some $t_1\in\RR$, $t_1\neq t_0$.
\end{crl}

Another consequence of inequality (\ref{habound}) is an upper bound for $\cosh\theta$ for maximal hypersurfaces contained in slabs.

\begin{crl} \label{intro:crl:u3}
	Let $\Mbar=\PP\times\RR$ be a Lorentzian product with complete, noncompact base $(\PP,\sigma)$. Suppose that the Ricci curvature of $(\PP,\sigma)$ satisfies (\ref{riccibound}) for some $G>0$ and that the sectional curvatures of $(\PP,\sigma)$ are uniformly bounded from below on $\PP$. Let $\psi:M\to\Mbar$ be a maximal spacelike hypersurface such that $\pi:M\to\PP$ is a covering map. If $\psi(M)$ is contained in a slab $\PP\times[t_0,t_1]$ of height $\delta=t_1-t_0$, with $-\infty<t_0<t_1<+\infty$, then
	\begin{equation} \label{habound2}
	\cosh\theta \leq e^{(m-1)\sqrt{2G}\delta}
	\end{equation}
	on $M$.
\end{crl}

We remark that in our last results we do not require any a priori upper bound on the hyperbolic angle $\theta$. The inequality (\ref{habound}) in Theorem \ref{intro:thm:coshb} is obtained as a gradient estimate for the height function $\tau$ considering $M$ with the metric $\hat\sigma$ defined above; this allows us to rely on an essentially non-parametric argument, which in the first steps is similar to that used in \cite{RSS} to obtain a gradient bound for minimal graphs in Riemannian products, even if $\psi(M)$ is not assumed to be a graph in $\Mbar$. Examples of gradient estimates obtained with similar techniques can be traced back to \cite{Ko}, \cite{Se} and the references therein. Uniqueness results for maximal spacelike hypersurfaces in Lorentzian products have also been obtained in \cite{Al}, \cite{AAA} in the case where $(\PP,\sigma)$ is a surface of nonnegative Gaussian curvature.

\section{The geometric setting}

We start this section with a preliminary result about the curvature tensors of a standard static spacetime.

\begin{lemma} \label{curvaturesbar}
	Let $\Mbar^{m+1}=\PP^m\times_h\RR$ be a standard static spacetime with warped metric
	\[
	\gbar=-h(x)^2 \di t^2+\metric_{\PP}
	\]
	where $h\in C^{\infty}(\PP)$, $h>0$ and $t$ is the standard coordinate on $\RR$. Then
	\begin{itemize}
		\item [(a)] the Riemannian curvature $(0,4)$-tensor $\Riembar$ of $\Mbar$ is given by
		\begin{equation} \label{riembar}
		\Riembar = \RiemP + \left(h \, \HessP(h)\right)\kulknom(\di t\otimes \di t),
		\end{equation}
		where $\RiemP$, $\HessP$ are the Riemannian curvature $(0,4)$-tensor and the Hessian of $\PP$, respectively, and$~\kulknom$ denotes the Kulkarni-Nomizu product of symmetric $(0,2)$-tensors;
		\item [(b)] the Ricci tensor of $\Mbar$ is given by
		\begin{equation} \label{ricbar}
		\Ricbar = \RicP - \frac{\HessP(h)}{h} + \left(h\,\DeltaP h\right)\di t\otimes \di t,
		\end{equation}
		where $\RicP$, $\DeltaP$ are the Ricci tensor and the Laplace-Beltrami operator of $\PP$, respectively.
	\end{itemize}
\end{lemma}

\begin{rmk}
	With a little abuse of notation, in the RHS's of formulas (\ref{riembar}) and (\ref{ricbar}) we are omitting writing the pullback $\pi_{\PP}^{\ast}$ for tensors and functions defined on $\PP$ to avoid introducing an unnecessarily complicated notation, so that $\RiemP$, $h\cdot\DeltaP h$, etc.~stand for $\pi_{\PP}^{\ast}(\RiemP)$, $(h\cdot\DeltaP h)\circ\pi_{\PP}$, etc. We also point out that we are adopting the definitions
	\begin{align*}
		\Riembar(X_1,X_2,X_3,X_4) & = \gbar(\Rbar(X_3,X_4)X_2,X_1), \\
		\RiemP(Y_1,Y_2,Y_3,Y_4) & = \sigma(\RP(Y_3,Y_4)Y_2,Y_1),
	\end{align*}
	for $X_i\in T\Mbar$, $Y_i\in T\PP$, $1\leq i\leq 4$, with $\Rbar$ the curvature operator of $\Mbar$ given by
	\begin{equation} \label{defR}
	\Rbar(V,W)Z=\nablabar_V\nablabar_W Z - \nablabar_W\nablabar_V Z - \nablabar_{[V,W]}Z
	\end{equation}
	for each $V,W,Z\in\mathfrak X(\Mbar)$ and with $\RP$ the similarly defined curvature operator of $\PP$. Moreover, we recall that for a given real vector space $V$ the Kulkarni-Nomizu product of two bilinear symmetric forms $\alpha,\beta:V\times V\to\RR$ is the $4$-linear form $\alpha\kulknom\beta:V\times V\times V\times V\to\RR$ defined by
	\begin{align*}
	(\alpha\kulknom\beta)(X_1,X_2,X_3,X_4) & = \alpha(X_1,X_3)\beta(X_2,X_4) + \alpha(X_2,X_4)\beta(X_1,X_3) \\
	& \phantom{=\;} - \alpha(X_1,X_4)\beta(X_2,X_3) - \alpha(X_2,X_3)\beta(X_1,X_4)
	\end{align*}
	for each $X_1,X_2,X_3,X_4\in V$.
\end{rmk}

\begin{proof}[Proof of Lemma \ref{curvaturesbar}]
	We apply Proposition 7.42 of \cite{O'N}. Note that the definition of the Riemannian curvature operator $\Rbar$ given by \cite{O'N}, see Lemma 3.35 therein, differs from the above \eqref{defR} by a minus sign, that is, $\Rbar(V,W)Z = \Rbar_{WV}Z = - \Rbar_{VW}Z$, where $\Rbar_{(\;\cdot\;,\;\cdot\;)}(\;\cdot\;)$ is the notation adopted in \cite{O'N} to denote the curvature operator.
	
	Let $p=(t,x)\in\overline{M}$ be a given point and let $X,Y,Z,V,W \in T_p\overline{M}$ be given tangent vectors at $p$, with $X,Y,Z$ tangent to the leaf $\{t\} \times \PP$ and $V,W$ tangent to the timelike curve $\RR \times \{x\}$. From formulas (1) and (3) of Proposition 7.42 of \cite{O'N} we get
	\begin{align*}
		\Riembar(\;\cdot\;,Z,X,Y) & = \Riem(\;\cdot\;,Z,X,Y), \\
		\Riembar(\;\cdot\;,Z,V,X) & = - \frac{\HessP(h) (X,Z)}{h}\gbar(V,\;\cdot\;), \\ \Riembar(\;\cdot\;,W,X,Y) & = 0, \\
		\Riembar(\;\cdot\;,W,X,V) & = - \gbar(V,W) \frac{\HessP(h) (X,\;\cdot\;)}{h}.	
	\end{align*}
	By the symmetry properties of the curvature tensor, see Proposition 3.36 of \cite{O'N}, we can rewrite the formulas above as
	\begin{align*}
		\Riembar(X,Y,Z,\;\cdot\;) & = \Riem(X,Y,Z,\;\cdot\;), \\
		\Riembar(X,Y,W,\;\cdot\;) & = 0, \\
		\Riembar(X,V,Z,\;\cdot\;) & = - \frac{\HessP(h) (X,Z)}{h} \gbar(V,\;\cdot\;), \\ 
		\Riembar(X,V,W,\;\cdot\;) & = \frac{\HessP(h) (X,\;\cdot\;)}{h} \gbar(V,W).
	\end{align*}
	A direct computation shows that the RHS of \eqref{riembar} also satisfy the identities above, since $\gbar(V,\;\cdot\;) = - h^2 (\di t\otimes \di t)(V,\;\cdot\;)$. By the symmetry properties of $\overline{\mathrm{Riem}}$ again, these identities uniquely determine its action on $T_p\overline{M}$. The identity \eqref{ricbar} is a straightforward consequence of \eqref{riembar} by the definition of Ricci tensor.
\end{proof}

\noindent Now consider a spacelike hypersurface $\psi:M\to\Mbar$ immersed into a standard static spacetime $\Mbar=\PP\times_h\RR$. In the Introduction we have defined the global unit normal vector field $N\in\mathfrak X^{\bot}(M)$ with the same time-orientation of the Killing vector field $\partial_t\in\mathfrak X(\Mbar)$. The hyperbolic cosine of the hyperbolic angle $\theta$ between $N$ and $\partial_t$ is the smooth function on $M$ given by
\[
\cosh\theta = -\gbar\left(N,\frac{1}{h\circ\pi}\,\partial_t\right).
\]

The global shape operator $A:TM\to TM$ of $\psi$ in the direction of $N$ is defined by
\[
g(AX,Y) = \gbar(\II(X,Y),N)
\]
for every $X,Y\in TM$, where $\II:TM\times TM\to TM^{\bot}$ is the second fundamental tensor of $\psi$ and $g=\psi^{\ast}\gbar$ is the Riemannian metric on $M$ induced by the ambient manifold $(\Mbar,\gbar)$. The mean curvature function $H\in C^{\infty}(M)$ satisfying $\Hvec=HN$, where $\Hvec$ is the mean curvature vector of $\psi$, is given by
\[
H = -\frac{1}{m}\,\trace_g(A).
\]
Each $p\in M$ has an open neighbourhood $U\subseteq M$ such that $\psi|_U$ is an embedding. For any smooth vector field $\tilde N$ extending $N$ on a neighbourhood $V\subseteq\Mbar$ of $\psi(U)$ we have
\[
H = \frac{1}{m}\,\div_{\gbar}(\tilde N)\circ\psi
\]
on $U$, as a consequence of Weingarten's equation
\[
g(AX,Y) = -\gbar(\nablabar_{\psi_{\ast}X}\tilde N,\psi_{\ast}Y).
\]
Hereafter, $\nablabar$ and $\div_{\gbar}$ denote the Levi-Civita connection of $(\Mbar,\gbar)$ and the corresponding divergence operator, and $\di\psi = \psi_{\ast} : TM \to T\Mbar$ is the pushforward of tangent vectors induced by $\psi$, that is, the differential of $\psi$.

Let $t$ be the standard coordinate on the factor $\RR$ of $\Mbar$. The composition of $t$ with the map $\pi_{\RR}\circ\psi:M\to\RR$ defines the smooth vertical height function
\[
\tau = t\circ\pi_{\RR}\circ\psi
\]
on $M$. We also set
\[
\hat h = h\circ\pi = h\circ\pi_{\PP}\circ\psi.
\]

\bigskip

\noindent We first consider the manifold $M$ with the induced Riemannian metric $g=\psi^{\ast}\gbar$. Let $\nabla$, $\mathrm{\Delta}$ and $\div_g$ be the Levi-Civita connection, the Laplace-Beltrami operator and the divergence operator on $(M,g)$, respectively. The gradient of the height function $\tau$ satisfies
\[
g(\nabla\tau,\;\cdot\;) = \di \tau = \psi^{\ast} \di t = -\frac{1}{\hat h^2} \, \gbar(\partial_t,\psi_{\ast}\;\cdot\;)
\]
and therefore is given by
\[
\nabla\tau = -\frac{1}{\hat h^2}\,\partial_t^{\top},
\]
where $\partial_t^{\top}\in\mathfrak X(M)$ is the tangential part of $\partial_t$ along $\psi$. In other words, $\partial_t^{\top}$ is the unique vector field such that $\psi_{\ast}\partial_t^{\top}|_p$ is the orthogonal projection of $\partial_t|_{\psi(p)}$ onto $\psi_{\ast}T_pM$ for each $p\in M$. Hence, we have
\[
|\nabla\tau|^2_g = \frac{1}{\hat h^4}[\gbar(\partial_t,\partial_t)+\gbar(\partial_t,N)^2] = \frac{1}{\hat h^2}(\cosh^2\theta-1)
\]
or, equivalently,
\begin{equation} \label{cosh-g}
\cosh\theta = \sqrt{1+\hat h^2|\nabla\tau|^2_g}.
\end{equation}

The height function $\tau$ satisfies the differential equation
\begin{equation} \label{Hequation-g1}
mH\cosh\theta = \hat h\mathrm{\Delta}\tau + 2g(\nabla\hat h,\nabla\tau)
\end{equation}
on $M$. Indeed, considering a local orthonormal reference frame $\{E_1,\dots,E_n\}$ on $(M,g)$ we compute
\begin{align*}
\mathrm{\Delta}\tau & = -\sum_{i=1}^n g\left(\nabla_{E_i}\left(\frac{1}{\hat h^2}\,\partial_t^{\top}\right),E_i\right) \\ & = -\sum_{i=1}^n E_i\left(\frac{1}{\hat h^2}\right)g(\partial_t^{\top},E_i) - \frac{1}{\hat h^2}\sum_{i=1}^n g(\nabla_{E_i}\partial_t^{\top},E_i) \\ & = - \frac{2}{\hat h}\sum_{i=1}^n g(E_i,\nabla \hat h)g(\nabla\tau,E_i) + \frac{\cosh\theta}{\hat h}\sum_{i=1}^n \gbar(\nablabar_{E_i}N,E_i) \\ & = -\frac{2}{\hat h}g(\nabla\hat h,\nabla\tau) - \frac{\cosh\theta}{\hat h}\sum_{i=1}^n g(AE_i,E_i),
\end{align*}
where in the third identity we have used the equality
\[
\partial_t^{\top} = \partial_t + \gbar(N,\partial_t)N = \partial_t - \hat h\cosh\theta N
\]
and the fact that $\gbar(\nablabar_{E_i}\partial_t,E_i)=0$ since $\partial_t$ is Killing.

Equation (\ref{Hequation-g1}) can be rewritten as
\begin{equation} \label{Hequation-g2}
\frac{mH\cosh\theta}{\hat h} = \frac{1}{\hat h^2}\div_g(\hat h^2\nabla\tau) = \mathrm{\Delta}_{-2\log\hat h}\tau,
\end{equation}
where the operator appearing in the RHS is the drifted Laplacian
\[
\mathrm{\Delta}_{-2\log\hat h} = \mathrm{\Delta} - g(\nabla(-2\log\hat h),\nabla\;\cdot\;) = \mathrm{\Delta} + \frac{2}{\hat h}g\left(\nabla\hat h,\nabla\;\cdot\;\right),
\]
or even in the form
\[
mH\cosh\theta = \frac{1}{\hat h}\div_g(\hat h^2\nabla\tau) = \div_{-\log\hat h,g}(\hat h\nabla\tau)
\]
where $\div_{-\log\hat h,g}$ is the weighted divergence operator
\begin{equation} \label{wdiv-g}
\div_{-\log\hat h,g} = \div_g - g(\nabla(-\log\hat h),\;\cdot\;) = \div_g + \frac{1}{\hat h}g(\nabla\hat h,\;\cdot\;).
\end{equation}

\bigskip

\noindent Our aim is now to provide a further expression for $mH$ in terms of $\tau$ by considering a second metric on $M$ different from the induced metric $g=\psi^{\ast}\gbar$. Towards this end we observe that the smooth map $\pi:M\to\PP$ is an immersion because $\psi$ is spacelike. Since $\dim M = m = \dim\PP$, we have that $\pi$ is a local diffeomorphism. This allows us to define a second Riemannian metric
\[
\hat\sigma = \pi^{\ast}\sigma
\]
on $M$. Let $D$, $\div_{\hat\sigma}$ be the Levi-Civita connection and the corresponding divergence operator on $(M,\hat\sigma)$, respectively. From the definition of $\gbar$ and $g$ it follows that
\begin{equation} \label{hatsigma-g}
\hat\sigma = g + \hat h^2 \di \tau^2.
\end{equation}
Thus, the gradient $D\tau$ of $\tau$ satisfies
\begin{align*}
\hat\sigma(D\tau,\;\cdot\;) & = \di\tau = g(\nabla\tau,\;\cdot\;) = \hat\sigma(\nabla\tau,\;\cdot\;) - \hat h^2 \di\tau(\nabla\tau)\di\tau \\
& = \hat\sigma(\nabla\tau,\;\cdot\;) - \hat h^2|\nabla\tau|^2_g\hat\sigma(D\tau,\;\cdot\;).
\end{align*}
As a consequence using (\ref{cosh-g})
\[
D\tau = \frac{\nabla\tau}{1+\hat h^2|\nabla\tau|^2_g} = \frac{\nabla\tau}{\cosh^2\theta}.
\]
From (\ref{hatsigma-g}) we obtain
\[
|D\tau|^2_{\hat\sigma} = \frac{1}{\cosh^4\theta}\left(|\nabla\tau|^2_g+\hat h^2 \di \tau(\nabla\tau)^2\right) = \frac{|\nabla\tau|^2_g}{1+\hat h^2|\nabla\tau|^2_g}.
\]
It follows that
\[
h^2|D\tau|^2_ {\hat\sigma} = \frac{\hat h^2|\nabla\tau|^2_g}{1+\hat h^2|\nabla\tau|^2_g} < 1
\]
so that the above can be equivalently written as
\[
|\nabla\tau|^2_g = \frac{|D\tau|^2_{\hat\sigma}}{1-\hat h^2|D\tau|^2_{\hat\sigma}}.
\]
From (\ref{cosh-g}) we infer
\begin{equation} \label{cosh-hatsigma}
\cosh\theta = \frac{1}{\sqrt{1-\hat h^2|D\tau|^2_{\hat\sigma}}}.
\end{equation}
The following proposition will be crucial for the results in the next sections.

\begin{prp}
	In the above notations we have
	\begin{equation} \label{Hequation-hatsigma1}
	mH = \frac{1}{\hat h}\div_{\hat\sigma}\left(\frac{\hat h^2 D\tau}{\sqrt{1-\hat h^2|D\tau|^2_{\hat\sigma}}}\right) = \div_{-\log\hat h,\hat\sigma}\left(\frac{\hat h D\tau}{\sqrt{1-\hat h^2|D\tau|^2_{\hat\sigma}}}\right),
	\end{equation}
	where the weighted divergence operator $\div_{-\log\hat h,\hat\sigma}$ is defined in way similar to that in (\ref{wdiv-g}).
\end{prp}

\begin{proof}
	Let $p\in M$ be a given point. Since $\pi$ is a local isometry, there exists an open neighbourhood $U_p\subseteq M$ of $p$ such that the restriction $\psi|_{U_p} : U_p\to\Mbar$ is an embedding and $\pi|_{U_p}:U_p\to \PP$ is an isometric diffeomorphism onto the image $U:=\pi(U_p)$. We fix the index ranges
	\[
		1\leq i,j,k,\dots\leq m, \quad 1\leq a,b,c,\dots\leq m+1.
	\]
	Up to restricting $U_p$, we can assume that there exists a local orthonormal coframe $\{\theta^i\}_{i=1}^m$ for $(\PP,\sigma)$ defined on $U$, with corresponding Levi-Civita connection forms $\{\theta^i_j\}_{i,j=1}^m$ defined by the structural equations
	\begin{equation} \label{structP}
		\begin{split}
			\di \theta^i & = -\theta^i_j\wedge\theta^j, \\
			\theta^i_j+\theta^j_i & = 0 \qquad \qquad \text{for } 1\leq i,j\leq m.
		\end{split}
	\end{equation}
	We recall that orthonormality of the coframe $\{\theta^i\}$ means that the metric $\sigma$ is expressed as $\sigma=\delta_{ij}\theta^i\otimes\theta^j$, with $\delta$ the Kronecker symbol, and we also recall that the Levi-Civita connection forms $\{\theta^i_j\}$ are the unique $1$-forms such that the covariant derivatives of the elements of the local (orthonormal) frame $\{e_i\}_{i=1}^m$ dual to $\{\theta^i\}$ are given by
	\[
		De_j = \theta^i_j\otimes e_i.
	\]
	We can define a local Lorentz orthonormal coframe $\{\omega^a\}_{a=1}^{m+1}$ on $U\times_h\RR\subseteq\Mbar$ by setting
	\begin{equation} \label{defcoframebar}
		\omega^i=\theta^i \text{ for } 1\leq i\leq m, \quad \omega^{m+1} = h \di t.
	\end{equation}
	In this case, Lorentz orthonormality means that the metric $\gbar$ is given by
	$$
		\gbar=\delta_{ij}\omega^i\otimes\omega^j-\omega^{m+1}\otimes\omega^{m+1}.
	$$
	The corresponding Levi-Civita connection forms $\{\omega^a_b\}_{a,b=1}^{m+1}$ are defined by the structural equations
	\begin{align*}
		\di \omega^a & = -\omega^a_b\wedge\omega^b, \\
		\omega^i_j+\omega^j_i = \omega^i_{m+1}-\omega^{m+1}_i & = \omega^{m+1}_{m+1} = 0 \quad \text{for } 1\leq i,j\leq m.
	\end{align*}
	A straightforward computation using \eqref{structP} and \eqref{defcoframebar} shows that
	\begin{equation} \label{LCcoframebar}
		\omega^i_j = \theta^i_j, \quad \omega^{m+1}_i = \frac{h_i}{h}\omega^{m+1},
	\end{equation}
	where $\di h=h_i\theta^i$ on $\PP$.
	
	Since $\pi|_{U_p}$ is an isometry onto the image, we have that $\{\pi^{\ast}\theta^i\}$ is a local orthonormal coframe on $U_p$ and that $\{\pi^{\ast}\theta^i_j\}$ are the corresponding Levi-Civita connection forms. From now on, we will omit writing the pullback $\pi^{\ast}$. We set $\alpha=\psi^{\ast}\omega^{m+1}$ on $M$. Note that $\alpha=\hat h \, \di \tau$. Writing $\alpha=\alpha_i\theta^i$ on $U_p$ and letting $\{\ebar_a\}_{a=1}^{m+1}$ be the local orthonormal frame for $(\Mbar,\gbar)$ dual to $\{\omega^a\}$ we can easily verify that the local orthonormal frame $\{e_i\}_{i=1}^m$ for $(M,\hat\sigma)$ dual to $\{\theta^i\}$ must satisfy
	\begin{equation} \label{pushforwardei}
	(\di \psi)_q(e_i)_q = (\bar e_i)_{\psi(q)} + \alpha_i(q)(\bar e_{m+1})_{\psi(q)}
	\end{equation}
	for each $q\in U_p$. From (\ref{pushforwardei}) it follows that for any smooth extension $\Nbar=\Nbar^a\ebar_a$ of $N$ on a neighbourhood of $\psi(U_p)\subseteq\Mbar$ the functions $N^a:=\Nbar^a\circ\psi\in C^{\infty}(U_p)$, $1\leq a\leq m+1$, must satisfy
	\begin{equation} \label{Nm+1cosh}
	N^{m+1} = -\gbar(N,\ebar_{m+1}) = \omega^{m+1}(N) = -\gbar\left(N,\frac{1}{h}\,\partial_t\right) = \cosh\theta
	\end{equation}
	and
	\begin{equation} \label{Nicosh}
	N^i = \alpha_i N^{m+1} = \alpha_i\cosh\theta
	\end{equation}
	for $1\leq i\leq m$.
	
	Let $\II=\II_{ij}\theta^i\otimes\theta^j\otimes N$ be the second fundamental tensor of the immersion $\psi$ on $U_p$. By Weingarten's equation we have
	\[
	\II_{ij} = \gbar\left(\nablabar_{(\di \psi)e_i}\Nbar,(\di \psi)e_j\right)
	\]
	for $1\leq i,j\leq m$. Writing $\nablabar \, \Nbar=\Nbar^a_b\omega^b\otimes\ebar_a$ we have, by the properties of covariant differentiation,
	\begin{equation} \label{covbarN}
	\Nbar^a_b\omega^b = \di \Nbar^a + \Nbar^b\omega^a_b
	\end{equation}
	for $1\leq a\leq m+1$. From (\ref{pushforwardei}) and (\ref{covbarN}) we deduce
	\begin{align*}
	\II_{ij} & = \gbar\left(\Nbar^a_b\omega^b((\di\psi)e_i)\ebar_a,(\di\psi)e_j\right) \\
	& = \delta_{jk}(\Nbar^k_b\omega^b)((\di\psi)e_i) - \alpha_j(\Nbar^{m+1}_b\omega^b)((\di\psi)e_i) \\
	& = \left[\delta_{jk}(\di N^k + N^a\psi^{\ast}\omega^k_a)-\alpha_j(\di N^{m+1}+N^a\psi^{\ast}\omega^{m+1}_a)\right](e_i)
	\end{align*}
	and from (\ref{Nicosh}) we get
	\begin{align*}
	\II_{ij}\theta^i & = \di(\alpha_j\cosh\theta) - \alpha_k\cosh\theta\cdot\psi^{\ast}\omega^k_j + \cosh\theta\cdot\psi^{\ast}\omega^{m+1}_j \\
	& \phantom{=\;} - \alpha_j(\di\cosh\theta) - \alpha_j\alpha_k\cosh\theta\cdot\psi^{\ast}\omega^k_{m+1} \\
	& = \cosh\theta\left[\di\alpha_j - \alpha_k\psi^{\ast}\omega^k_j + \psi^{\ast}\omega^{m+1}_j - \alpha_j\alpha_k\psi^{\ast}\omega^k_{m+1}\right].
	\end{align*}
	Using (\ref{LCcoframebar}) we can further write
	\begin{equation} \label{secfundtens}
	\II_{ij}\theta^i = \cosh\theta\left[\di\alpha_j - \alpha_k\theta^k_j + \left(\hat h_j - \alpha_j\alpha_k\hat h^k\right)\frac{\alpha}{\hat h}\right],
	\end{equation}
	where $\di\hat h = \hat h_k\theta^k$ on $U_p$. Again, by the properties of covariant differentiation we have
	\begin{equation} \label{covalpha}
	(\di\alpha_j - \alpha_k\theta^k_j)\otimes\theta^j = D\alpha.
	\end{equation}
	
	Note that the metric $g=\psi^{\ast}\gbar$ is given by $g=g_{ij}\theta^i\otimes\theta^j$ with $g_{ij} = \delta_{ij} - \alpha_i\alpha_j$ and that, by (\ref{cosh-hatsigma}),
	\[
	\cosh^2\theta = \frac{1}{1-\delta^{ij}\alpha_i\alpha_j}.
	\]
	A straightforward computation shows that the elements of the inverse matrix $(g^{ij}) = (g_{ij})^{-1}$ are
	\[
	g^{ij} = \delta^{ij} + \cosh^2\theta\,\alpha_i\alpha_j.
	\]
	Therefore, the mean curvature function
	\[
	H = -\frac{1}{m}\trace_g(A) = \frac{1}{m}\trace_g(\II_{ij}\theta^i\otimes\theta^j)
	\]
	is given by
	\begin{align*}
	mH & = \cosh\theta\cdot\trace_{\hat\sigma}(D\alpha) + \cosh^3\theta\cdot D\alpha(\hat h D\tau,\hat h D\tau) \\
	& \phantom{=\;} + \cosh\theta\cdot \di\hat h(D\tau) - \cosh\theta\cdot \di\hat h(D\tau)\cdot\delta^{ij}\alpha_i\alpha_j \\
	& \phantom{=\;} + \cosh^3\theta\cdot \di\hat h(D\tau)\cdot\delta^{ij}\alpha_i\alpha_j - \cosh^3\theta\cdot \di\hat h(D\tau)(\delta^{ij}\alpha_i\alpha_j)^2 \\
	& = \trace_{\hat\sigma}(D(\cosh\theta\cdot\alpha)) + \cosh\theta\cdot \di\hat h(D\tau) \\
	& = \div_{\hat\sigma}\left(\frac{\hat h D\tau}{\sqrt{1-\hat h^2|D\tau|^2}}\right) + \frac{\hat\sigma(D\hat h,D\tau)}{\sqrt{1-\hat h^2|D\tau|^2}} \\
	\end{align*}
	that is, (\ref{Hequation-hatsigma1}).
\end{proof}

\begin{rmk}
	We observe that, if the map $\pi:M\to\PP$ is injective, $\psi(M)\subseteq\PP\times\RR$ is the graph of a smooth function $u\in C^{\infty}(\Omega)$, $\Omega\subseteq\PP$ an open domain. In this case $\Omega=\pi(M)$, $u=\tau\circ\pi^{-1}$ and $\pi:(M,\hat\sigma)\to(\Omega,\sigma|_{\Omega})$ is an isometry. Thus equation (\ref{Hequation-hatsigma1}) reduces to the prescribed mean curvature equation
	\[
	m\,(H\circ\pi^{-1}) = \frac{1}{h} \, \div_{\sigma}\left(\frac{h^2 Du}{\sqrt{1-h^2|Du|^2_{\sigma}}}\right),
	\]
	with $Du$ the gradient of $u$ in $(\Omega,\sigma|_{\Omega})$.
\end{rmk}

\section{Hyperbolic angle estimates}

Throughout this section we consider $M$ furnished with the metric $\hat\sigma$ and we denote by $H$ and $\tau$ the mean curvature function and the vertical height function, respectively, of the immersion $\psi:M\to\Mbar$, as defined in the previous section; in particular $H$ is the mean curvature function in the direction of the normal $N$ to the immersion. Integration over domains $\Omega\subseteq M$ will always be intended with respect to the volume element induced by $\hat\sigma$ and integration over boundaries of sufficiently regular domains is intended with respect to the corresponding $(m-1)$-dimensional Hausdorff measure. In particular, for each relatively compact domain $\Omega\subseteq M$ we define the weighted volume
\[
\vol_{-\log\hat h}(\Omega) = \int_{\Omega}\hat h
\]
and, if $\partial\Omega$ is sufficiently regular, we also define
\[
\vol_{-\log\hat h}(\partial\Omega) = \int_{\partial\Omega}\hat h.
\]
This apparently cumbersome notation for the weighted volume is due to the standard notation
\[
\vol_f(\Omega)=\int_{\Omega}e^{-f}
\]
for $f\in C^{\infty}(M)$. For each relatively compact domain $\Omega\subseteq M$ we shall also denote the weighted integral mean of $H$ over $\Omega$ by
\begin{equation} \label{defHmean}
\Hmean(\Omega)=\frac{\int_{\Omega} H\hat h}{\vol_{-\log\hat h}(\Omega)}.
\end{equation}

\bigskip

\noindent Our aim is to deduce two consequences of equation (\ref{Hequation-hatsigma1}). The first is based on a clever idea of Salavessa, \cite{Sal}. Suppose that $(M,\hat\sigma)$ is complete, noncompact. We introduce the weighted Cheeger constant $\mathfrak C_{-\log\hat h}$ of the weighted manifold $(M,\hat\sigma,\hat h)$ by setting
\begin{equation} \label{cheeger}
\begin{split}
\mathfrak C_{-\log\hat h} & = \inf\left\{\frac{\int_{\partial\Omega}\hat h}{\int_\Omega\hat h} : \text{$\Omega\subseteq M$ relatively compact domain, $\partial\Omega$ smooth} \right\} \\
& = \inf\left\{\frac{\vol_{-\log\hat h}(\partial\Omega)}{\vol_{-\log\hat h}(\Omega)} : \text{$\Omega\subseteq M$ relatively compact domain, $\partial\Omega$ smooth} \right\}.
\end{split}
\end{equation}
Note also that the definition (\ref{cheeger}) of the Cheeger constant has to be changed in case $M$ is compact.

Introducing the operator $\Deltahat_{-\log\hat h}=\Deltahat+\hat\sigma\left(\frac{D\hat h}{\hat h},D\;\cdot\;\right)$ where $\Deltahat$ is the Laplace-Beltrami operator of $(M,\hat\sigma)$ and indicating with
\[
\lambda_1^{\Deltahat_{-\log\hat h}}(M)=\inf_{\substack{\varphi\in C^{\infty}(M)\\ \varphi\not\equiv 0}}\frac{\int_{\Omega}|D\varphi|^2\hat h}{\int_{\Omega}\varphi^2\hat h}
\]
its spectral radius, following the original argument in Cheeger, \cite{Che}, it is not hard to show that
\begin{equation} \label{cheegerthm}
\lambda_1^{\Deltahat_{-\log\hat h}}(M) \geq \frac{1}{4}\left(\mathfrak C_{-\log\hat h}\right)^2.
\end{equation}
We are now ready to prove
\begin{prp}
	Let $\Mbar^{m+1}=\PP^m\times_h\RR$ be a standard static spacetime with complete, noncompact base $(\PP,\sigma)$ and let $\psi:M\to\Mbar$ be a spacelike hypersurface such that $\pi:M\to\PP$ is a covering map. Let $\mathfrak C_{-\log\hat h}$ be the weighted Cheeger constant of $(M,\hat\sigma,\hat h)$ and let $\cosh\theta$ be the hyperbolic cosine of the hyperbolic angle of the immersion $\psi$. Suppose that
	\begin{equation} \label{sal:boundcosh}
	\sup_M\cosh\theta = \cosh\theta^{\ast} < +\infty.
	\end{equation}
	Then, the mean curvature $H$ in the direction of $N$ satisfies
	\begin{equation} \label{sal:Hbound}
	\inf_{\substack{\Omega\subseteq M,\\\overline{\Omega}\text{ compact}}} |\Hmean(\Omega)| \leq \frac{\sqrt{\cosh^2\theta^{\ast}-1}}{m} \, \mathfrak C_{-\log\hat h}.
	\end{equation}
	In particular, if $H$ is a constant, $\hat h\notin L^1(M,\hat\sigma)$ and the function $r\mapsto\vol_{-\log\hat h}(\Bhat_r)$ has subexponential growth for some (hence, any) $o\in M$, where $\Bhat_r=B^{\hat\sigma}_r(o)$, then $\psi:M\to\Mbar$ is maximal.
\end{prp}

\begin{proof}
	Let $\hat\sigma=\pi^{\ast}\sigma$ and $\Omega\subseteq M$ a relatively compact domain with smooth boundary and outward unit normal vector $\nu$. We apply the divergence theorem to (\ref{Hequation-hatsigma1}) to obtain
	\begin{equation} \label{divthm}
	m|\Hmean(\Omega)|\vol_{-\log\hat h}(\Omega) = \left|\int_{\Omega}mH\hat h\right| = \left|\int_{\partial\Omega}\frac{\hat h^2\hat\sigma(D\tau,\nu)}{\sqrt{1-\hat h^2|D\tau|^2}}\right| \leq \int_{\partial\Omega}\frac{\hat h^2|D\tau|}{\sqrt{1-\hat h^2|D\tau|^2}}
	\end{equation}
	where $\tau=t\circ\pi_{\PP}\circ\psi$ and where the measure we are considering is with respect to $\hat\sigma$. Note that this latter is complete because $g=\psi^{\ast}\gbar$ is so by assumption. As we have already observed we have the validity of (\ref{cosh-hatsigma}). From (\ref{sal:boundcosh}) and (\ref{divthm}) it follows that
	\[
	m|\Hmean(\Omega)|\vol_{-\log\hat h}(\Omega) \leq \int_{\partial\Omega}\hat h\sqrt{\cosh^2\theta-1} \leq \sqrt{\cosh^2\theta^{\ast}-1}\int_{\partial\Omega}\hat h
	\]
	or, in other words,
	\[
	|\Hmean(\Omega)| \leq \frac{\sqrt{\cosh^2\theta^{\ast}-1}}{m}\frac{\vol_{-\log\hat h}(\partial\Omega)}{\vol_{-\log\hat h}(\Omega)}.
	\]
	Letting $\Omega$ run over all relatively compact domains in $M$ with smooth boundary and taking the infimum, from the previous inequality we deduce the validity of (\ref{sal:Hbound}).
	
	To complete the proof we observe, see \cite{BMR}, that if $\vol_{-\log\hat h}(M)=+\infty$ and $\vol_{-\log\hat h}(B^{\hat\sigma}_r)$ has subexponential growth, then
	\[
	\lambda_1^{\Deltahat_{-\log\hat h}}(M)=0
	\]
	and the final part of the Proposition follows from (\ref{sal:Hbound}) and (\ref{cheegerthm}).
\end{proof}

Consider now the case where $M=\PP$ and $\psi:M\to\Mbar=\PP\times_h\RR$ is a graph given by the function $u:\PP\to\RR$, that is, $\psi(x)=\psi_u(x)=(x,u(x))$. In this case $\pi(x)=\pi_{\PP}\circ\psi(x)=x$ and therefore
\[
\hat\sigma = \pi^{\ast}\sigma = (\pi_{\PP}\circ\psi)^{\ast}(\sigma) = \sigma, \quad \hat h = h\circ\pi = h, \quad \tau = t\circ\pi = t.
\]
It follows that (\ref{Hequation-hatsigma1}) becomes exactly the mean curvature equation of the graph $\psi_u$, that is,
\[
mH = \frac{1}{h} \, \div_{\sigma}\left(\frac{h^2 Du}{\sqrt{1-h^2|Du|^2_{\sigma}}}\right).
\]
As a consequence we have the following
\begin{crl}
	Let $\psi:\PP\to\PP^m\times_h\RR$ be a spacelike graph given by the function $u:\PP\to\RR$. Assume that $(\PP,\sigma)$ is complete and let $\mathfrak C_{-\log h}^{\PP}$ be the weighted Cheeger constant of $(\PP,\sigma,h)$. Assume the hyperbolic cosine of the hyperbolic angle is bounded above by $\cosh\theta^{\ast}<+\infty$. Then the mean curvature $H$ of the graph satisfies
	\[
	\inf_{\substack{\Omega\subseteq M,\\\overline{\Omega}\text{ compact}}} |\Hmean(\Omega)| \leq \frac{\sqrt{\cosh^2\theta^{\ast}-1}}{m} \, \mathfrak C_{-\log h}^{\PP}.
	\]
	In particular, if $H$ is a constant, $h\notin L^1(\PP,\sigma)$ and $\vol_{-\log h}(B^{\sigma}r)$ has subexponential growth, then $\psi_u$ is maximal.
\end{crl}

\begin{rmk}
	Note that completeness of $g=\psi^{\ast}_u\gbar$ implies that of $\hat\sigma=\sigma$ in this case. However completeness of $\hat\sigma$ is generally less stringent.
\end{rmk}

\bigskip

\noindent We shall now analyze more consequences of equation (\ref{Hequation-hatsigma1}) that will lead to a proof of Theorems \ref{intro:thm:cosh1} and \ref{intro:thm:cosh2}. In particular, we develop the observations above restricting ourselves to consider the case where $\Omega$ runs over all geodesic balls $\Bhat_r=B^{\hat\sigma}_r(o)$ of $(M,\hat\sigma)$ centered at a fixed point $o\in M$. Towards this end, set $\disthat = d_{\hat\sigma}(o,\;\cdot\;)$, with $d_{\hat\sigma}$ the distance on $M$ induced by the metric $\hat\sigma$. When $\Bhat_r$ is a relatively compact subset in $M$, we have the validity of the next two lemmas. The first, a weighted divergence theorem with low regularity assumptions, is well known in the non-weighted version. We report its proof here for the sake of completeness.

\begin{lemma} \label{lemma:weakdivthm}
	Let $r>0$ be such that $\Bhat_r$ is relatively compact in $M$. Then
	\begin{equation} \label{weakdivthm}
	\int_{\partial \Bhat_s}\hat h\,\hat\sigma\left(\frac{\hat h D\tau}{\sqrt{1-\hat h^2|D\tau|^2}},D\disthat\right) = m\int_{\Bhat_s}\hat h H
	\end{equation}
	for a.e.~$s\in(0,r)$.
\end{lemma}

\begin{proof}
	The function $\disthat$ is Lipschitz continuous and therefore the integral in the LHS of (\ref{weakdivthm}) is well defined for a.e.~$s\in(0,r)$. Let $f:(0,r)\to\RR$ be defined by
	\[
	f(s) = \int_{\Bhat_s}\hat h\,\hat\sigma\left(\frac{\hat h D\tau}{\sqrt{1-\hat h^2|D\tau|^2}},D\disthat\right).
	\]
	By the coarea formula (see \cite{SY2}, p.~89) we have
	\begin{equation} \label{coarea}
	\begin{split}
	f(s_1)-f(s_0) & = \int_{\Bhat_{s_1}\setminus \Bhat_{s_0}} \hat h\,\hat\sigma\left(\frac{\hat h D\tau}{\sqrt{1-\hat h^2|D\tau|^2}},D\disthat\right) \\
	& = \int_{s_0}^{s_1}\left[\int_{\partial \Bhat_s}\hat h\,\hat\sigma\left(\frac{\hat h D\tau}{\sqrt{1-\hat h^2|D\tau|^2}},D\disthat\right)\right] \di s
	\end{split}
	\end{equation}
	for each $0<s_0<s_1<r$. On the other hand, we also have $f\in C^1((0,r))$ and
	\begin{equation} \label{fderivative}
	f'(s) = m\int_{\Bhat_s}\hat h H
	\end{equation}
	for each $s\in(0,r)$. Indeed, for every $0<t_0<t_1<r$ let $\psi_{t_0,t_1}\in\Lip_c(M)$ be defined by
	\[
	\psi_{t_0,t_1}(x) = \begin{cases} 1 & \text{if }x\in \Bhat_{t_0}, \\ \dfrac{t_1-\disthat(x)}{t_1-t_0} & \text{if }x\in \Bhat_{t_1}\setminus \Bhat_{t_0}, \\ 0 & \text{if }x\in M\setminus \Bhat_{t_1}. \end{cases}
	\]
	We use (\ref{Hequation-hatsigma1}) and we integrate by parts to get
	\begin{equation} \label{ffinitediff}
	m\int_{\Bhat_{t_1}}\hat hH\psi_{t_1,t_0} = \int_{\Bhat_{t_1}\setminus \Bhat_{t_0}} \hat h\,\hat\sigma\left(\frac{\hat h D\tau}{\sqrt{1-\hat h^2|D\tau|^2}},D\psi_{t_1,t_0}\right) = \frac{f(t_1)-f(t_0)}{t_1-t_0}.
	\end{equation}
	For any given $s\in(0,r)$, identity (\ref{fderivative}) follows by fixing either $t_0=s$ or $t_1=s$, letting $t_1-t_0\to0^+$ in (\ref{ffinitediff}) and applying the dominated convergence theorem. By (\ref{coarea}) we get
	\[
	\int_{s_0}^{s_1}\left[f'(s)-\int_{\partial \Bhat_s}\hat h\,\hat\sigma\left(\frac{\hat h D\tau}{\sqrt{1-\hat h^2|D\tau|^2}},D\disthat\right)\right] \di s = 0
	\]
	for each $0<s_0<s_1<r$. Then (\ref{weakdivthm}) follows for a.e.~$s\in(0,r)$.
\end{proof}

\begin{lemma} \label{lemma:logdiff}
	Let $r>0$ be such that $\Bhat_r$ is relatively compact in $M$. Then
	\begin{equation} \label{logdiff}
	\log\int_{\Bhat_r}\hat h - \log\int_{\Bhat_R}\hat h = \int_R^r\left[\frac{\int_{\partial \Bhat_s}\hat h}{\int_{\Bhat_s}\hat h}\right] \di s
	\end{equation}
	for each $R\in(0,r)$.
\end{lemma}

\begin{proof}
	Let $f:(0,r)\to\RR$ be defined by
	\[
	f(s) = \int_{\Bhat_s}\hat h.
	\]
	By the coarea formula again, we have that
	\[
	f(s_1)-f(s_0) = \int_{\Bhat_{s_1}\setminus \Bhat_{s_0}}\hat h = \int_{s_0}^{s_1}\left[\int_{\partial \Bhat_s}\hat h\right] \di s
	\]
	for each $0<s_0<s_1<r$. The function $f$ is absolutely continuous on $(0,r)$, hence it is a.e.~differentiable on $(0,r)$ and there exists $\varphi\in L^1((0,r))$ such that
	\[
	f(s_1)-f(s_0) = \int_{s_0}^{s_1}\varphi(s) \di s
	\]
	for each $0<s_0<s_1<r$. Therefore,
	\[
	\varphi(s) = \int_{\partial \Bhat_s}\hat h
	\]
	for a.e.~$s\in(0,r)$. For each $R\in(0,r)$, the function $\log f$ is also absolutely continuous on $(R,r)$ and its a.e.~defined derivative is a.e.~equal to $\varphi/f$. Hence, (\ref{logdiff}) follows.
\end{proof}

We will also need the following ``long-range'' version of the mean value theorem.

\begin{lemma} \label{lemma:meanvalue}
	Let $f:\RR^+\to\RR$ be a measurable, locally integrable function. Let $0<R<R_1$, $\varepsilon>0$ be given. There exists $R_2>R_1$ such that
	\begin{equation} \label{meanvalue}
	\essinf_{(R_1,r)} f - \varepsilon < \frac{1}{r-R}\int_R^r f(t) \di t < \esssup_{(R_1,r)} f + \varepsilon
	\end{equation}
	for each $r>R_2$.
\end{lemma}

\begin{proof}
	We only prove the existence of $R_2>R_1$ such that the first inequality in (\ref{meanvalue}) is satisfied for each $r>R_2$. The proof of the second one is analogous. Set
	\[
	f_{\ast}(r) = \essinf_{(R_1,r)}f
	\]
	for each $r>R_1$. Suppose that $r>R_1+\varepsilon$ is such that
	\[
	\frac{1}{r-R}\int_R^r f(t) \di t \leq \essinf_{(R_1,r)}f - \varepsilon = f_{\ast}(r) - \varepsilon.
	\]
	Multiplying both sides by $r-R$ and noting that $f_{\ast}$ is nonincreasing, we get
	\begin{align*}
	-\varepsilon(r-R) & \geq\int_R^r[f(t)-f_{\ast}(t)] \di t \\
	& = \int_R^{R_1+\varepsilon}[f(t)-f_{\ast}(t)] \di t + \int_{R_1+\varepsilon}^r[f(t)-f_{\ast}(t)] \di t \\
	& \geq \int_R^{R_1+\varepsilon}[f(t)-f_{\ast}(R_1+\varepsilon)] \di t + \int_{R_1+\varepsilon}^r[f(t)-f_{\ast}(t)] \di t \\
	& \geq \int_R^{R_1+\varepsilon}[f(t)-f_{\ast}(R_1+\varepsilon)] \di t.
	\end{align*}
	Hence, it must be
	\begin{equation} \label{rupperbound}
	r \leq R - \frac{1}{\varepsilon}\int_R^{R_1+\varepsilon}[f(t)-f_{\ast}(R_1+\varepsilon)] \di t.
	\end{equation}
	Note that the RHS of (\ref{rupperbound}) is different from $+\infty$. Setting
	\[
	R_2 := \max\left\{R_1+\varepsilon,R-\frac{1}{\varepsilon}\int_R^{R_1+\varepsilon}[f(t)-f_{\ast}(R_1+\varepsilon)] \di t\right\}
	\]
	we have that the first inequality in (\ref{meanvalue}) holds for each $r>R_2$.
\end{proof}

We are now ready to prove the main result of this section.

\begin{thm} \label{thm:coshestimate1}
	Consider the manifold $(M,\hat\sigma)$. Then, with the notations introduced above, the following statements hold true.
	\begin{itemize}
		\item [(a)] If $r>0$ is such that $\Bhat_r$ is relatively compact in $M$, then
		\begin{equation} \label{boundarycoshest1}
		\max_{\partial \Bhat_s}\sqrt{\cosh^2\theta-1}\,\frac{\int_{\partial \Bhat_s}\hat h}{\int_{\Bhat_s}\hat h} \geq m|\Hmean(\Bhat_s)|
		\end{equation}
		for a.e.~$s\in(0,r)$ and
		\begin{equation} \label{disccoshest1}
		\max_{\overline{\Bhat_r}\setminus\Bhat_R}\sqrt{\cosh^2\theta-1}\,\frac{\log\int_{\Bhat_r}\hat h - \log\int_{\Bhat_R}\hat h}{r-R} \geq \min_{R\leq\rho\leq r}m|\Hmean(\Bhat_{\rho})|
		\end{equation}
		for each $R\in(0,r)$.
		\item [(b)] If $(M,\hat\sigma)$ is complete and $Q:\mathbb R^+\to\mathbb R^+$ is such that
		\begin{equation} \label{growthbound1}
		\liminf_{r\to+\infty}\frac{\log\int_{\Bhat_r}\hat h}{Q(r)} \leq 1,
		\end{equation}
		then
		\begin{equation} \label{limsupcoshest1}
		\limsup_{r\to+\infty}|\Hmean(\Bhat_r)| \leq \frac{1}{m}\left(\limsup_{M\ni x\to\infty}\sqrt{\cosh^2\theta(x)-1}\right)\left(\limsup_{r\to+\infty}\frac{Q(r)}{r}\right)
		\end{equation}
		as long as the RHS of \eqref{limsupcoshest1} does not present in the indeterminate form $(+\infty)\cdot0$ or $0\cdot(+\infty)$.
	\end{itemize}
\end{thm}

\begin{proof}
	\begin{itemize}
		\item [(a)] Suppose that $\Bhat_r$ is relatively compact in $M$. From (\ref{cosh-hatsigma}) as we have already seen we have
		\[
		\frac{\hat h|D\tau|}{\sqrt{1-\hat h^2|D\tau|^2}} = \sqrt{\cosh^2\theta-1},
		\]
		so, since $|D\disthat|=1$ a.e., (\ref{boundarycoshest1}) immediately follows from (\ref{weakdivthm}) for a.e.~$s\in(0,r)$. As a consequence,
		\[
		\max_{\overline{\Bhat_r}\setminus\Bhat_R}\sqrt{\cosh^2\theta-1}\,\frac{\int_{\partial \Bhat_s}\hat h}{\int_{\Bhat_s}\hat h} \geq \min_{R\leq\rho\leq r}m|\Hmean(\Bhat_{\rho})|
		\]
		for a.e.~$s\in(0,r)$. Integrating both sides on $[R,r]$ with respect to $s$, by (\ref{logdiff}) we get (\ref{disccoshest1}).
		\item [(b)] Suppose that $(M,\hat\sigma)$ is complete. By recursively applying Lemma \ref{lemma:meanvalue} we obtain the existence of a sequence $\{r_n\}_{n\in\mathbb N}$ such that $r_{n+1}>r_n>R$ and
		\[
		\frac{\log\int_{\Bhat_{r_{n+1}}}\hat h - \log\int_{\Bhat_R}\hat h}{r_{n+1}-R} > \essinf_{s\in(r_n,r_{n+1})}\frac{\int_{\partial \Bhat_s}\hat h}{\int_{\Bhat_s}\hat h} - \frac{1}{n}
		\]
		for each $n\in\mathbb N$, and satisfying $r_n\to+\infty$ as $n\to+\infty$ and
		\begin{equation} \label{attainedliminf}
		\lim_{n\to+\infty}\frac{\log\int_{\Bhat_{r_n}}\hat h}{Q(r_n)} = \liminf_{r\to+\infty}\frac{\log\int_{\Bhat_r}\hat h}{Q(r)}.
		\end{equation}
		This implies the existence of an increasing sequence $\{s_n\}_{n\in\mathbb N}$ such that $r_n<s_n<r_{n+1}$ and
		\[
		\frac{m|\Hmean(\Bhat_{s_n})|}{\max_{\partial \Bhat_{s_n}}\sqrt{\cosh^2\theta-1}} \leq \frac{\int_{\partial \Bhat_{s_n}}\hat h}{\int_{\Bhat_{s_n}}\hat h} \leq \frac{\log\int_{\Bhat_{r_{n+1}}}\hat h - \log\int_{\Bhat_R}\hat h}{r_{n+1}-R} + \frac{1}{n}
		\]
		for each $n\in\mathbb N$. By (\ref{growthbound1}) and (\ref{attainedliminf}) we have
		\begin{equation} \label{limsup:Q}
		\begin{split}
		\limsup_{n\to+\infty} \left[ \frac{\log\int_{\Bhat_{r_{n+1}}}\hat h - \log\int_{\Bhat_R}\hat h}{r_{n+1}-R} + \frac{1}{n} \right] & = \limsup_{n\to+\infty}\frac{\log\int_{\Bhat_{r_{n+1}}}\hat h}{r_{n+1}} \\
		& \leq \limsup_{n\to+\infty} \frac{Q(r_n)}{r_n} \\
		& \leq \limsup_{r\to+\infty} \frac{Q(r)}{r}.
		\end{split}
		\end{equation}
		Since $s_n\to+\infty$ as $n\to+\infty$ we have
		\begin{equation} \label{limsup:cosh}
		\limsup_{n\to+\infty}\max_{\partial \Bhat_{s_n}}\sqrt{\cosh^2\theta-1}\leq\limsup_{M\ni x\to\infty}\sqrt{\cosh^2\theta(x)-1}
		\end{equation}
		and
		\begin{align*}
		\limsup_{r\to+\infty}|\Hmean(\Bhat_r)| & \leq \frac{1}{m}\limsup_{n\to+\infty}\left[\left(\max_{\partial \Bhat_{s_n}}\sqrt{\cosh^2\theta-1}\right)\left(\frac{\log\int_{\Bhat_{r_{n+1}}}\hat h-\log_{\int_{\Bhat_R}}\hat h}{r_{n+1}-R}+\frac{1}{n}\right)\right].
		\end{align*}
		By (\ref{limsup:Q}) and (\ref{limsup:cosh}), the claim follows.
	\end{itemize}
\end{proof}

From \cite{MRS} we draw the following

\begin{lemma} \label{lemma:genBGthm}
	Let $r>0$ be such that $\Bhat_r$ is relatively compact in $M$. Let $G\in C^0([0,r])$ be such that
	\begin{equation} \label{modBEriccibound}
	\Richat - \frac{\Hesshat(\hat h)}{\hat h} \geq - m(G\circ\disthat)\,\hat\sigma
	\end{equation}
	on $\Bhat_r$, where $\Richat$ and $\Hesshat$ are the Ricci tensor and the Hessian operator of $(M,\hat\sigma)$ and inequality (\ref{modBEriccibound}) is intendend in the sense of quadratic forms. Let $k\in C^2([0,r])$ be a solution of the problem
	\begin{equation} \label{ksupersol}
	\begin{cases}
	k''-Gk \geq 0, \\
	k(0)=0, \quad k'(0)=1.
	\end{cases}
	\end{equation}
	Then
	\begin{itemize}
		\item [(a)] $k\circ\disthat>0$ on $\Bhat_r$,
		\item [(b)] the inequality
		\[
		\Deltahat_{-\log\hat h}\disthat \leq m\frac{k'\circ\disthat}{k\circ\disthat}
		\]
		holds pointwise on $\Bhat_r\setminus(\cut(o)\cup\{o\})$ and weakly on $\Bhat_r$,
		\item [(c)] the inequality
		\[
		\frac{\int_{\partial \Bhat_{s_0}}\hat h}{k(s_0)^m} \geq \frac{\int_{\partial \Bhat_{s_1}}\hat h}{k(s_1)^m}
		\]
		holds for a.~e.~$s_0,s_1\in(0,r)$ with $s_0<s_1$.
	\end{itemize}
\end{lemma}

Now, let $q=\pi(o)\in\PP$. We define the function
\[
\distP = d_{\sigma}(q,\;\cdot\;)
\]
on $\PP$, where $d_{\sigma}$ is the distance on $\PP$ induced by the Riemannian metric $\sigma$. For each $r>0$, we set
\[
\BP_r = B^{\sigma}_r(q) = \{y\in\PP:\distP(y)<r\}.
\]
The map $\pi:(M,\hat\sigma)\to(\PP,\sigma)$ is a local isometry, hence it is distance decreasing (see, for instance, Proposition 21 of \cite{Pe}). So, we have
\begin{equation} \label{gammaineq}
\distP\circ\pi \leq \disthat
\end{equation}
on $M$, and
\begin{equation} \label{ballinclusion}
\pi(\Bhat_r) \subseteq \BP_r
\end{equation}
for each $r>0$. Moreover, we have the equality sign in (\ref{gammaineq}) and (\ref{ballinclusion}) when $\pi$ is a diffeomorphism. Then, from Theorem \ref{thm:coshestimate1} we obtain the following

\begin{thm} \label{thm:coshestimate2}
	Let $G\in C^0(\RR^+_0)$ be such that
	\begin{equation} \label{modBEricciboundP}
	\RicP - \frac{\HessP(h)}{h} \geq - m(G\circ\distP)\,\sigma
	\end{equation}
	on $\PP$, where $\RicP$ and $\HessP$ are the Ricci tensor and the Hessian operator of $(\PP,\sigma)$. Let $k\in C^2(\RR^+_0)$ be a solution of problem (\ref{ksupersol}). If $G$ is nondecreasing, then the following statements hold true.
	\begin{itemize}
		\item [(a)] If $r>0$ is such that $\Bhat_r=B^{\hat\sigma}_r$ is relatively compact in $M$, then
		\begin{equation} \label{boundarycoshest2}
		\max_{\partial \Bhat_s}\sqrt{\cosh^2\theta-1}\,\frac{k(s)^m}{\int_0^s k(t)^m \di t} \geq m|\Hmean(\Bhat_s)|
		\end{equation}
		for a.e.~$s\in(0,r)$ and
		\begin{equation} \label{disccoshest2}
		\max_{\overline{\Bhat_r}\setminus\Bhat_R}\sqrt{\cosh^2\theta-1}\,\frac{\log\int_0^r k(t)^m \di t - \log\int_0^R k(t)^m \di t}{r-R} \geq \min_{R\leq\rho\leq r}m|\Hmean(\Bhat_{\rho})|
		\end{equation}
		for each $R\in(0,r)$.
		\item [(b)] If $(M,\hat\sigma)$ is complete, then
		\begin{equation} \label{limsupcoshest2}
		\limsup_{r\to+\infty}|\Hmean(\Bhat_r)| \leq \frac{1}{m}\left(\limsup_{M\ni x\to\infty}\sqrt{\cosh^2\theta(x)-1}\right)\left(\limsup_{r\to+\infty}\frac{\log\int_0^r k(t)^m \di t}{r}\right)
		\end{equation}
		with the same observation of Theorem \ref{thm:coshestimate1} about the RHS of (\ref{limsupcoshest2}).
	\end{itemize}
	If $\pi$ is a diffeomorphism, then the previous statements hold true without requiring the monotonicity of $G$.
\end{thm}

\begin{proof}
	We have
	\begin{align*}
	\Richat - \frac{\Hesshat(\hat h)}{\hat h} & = \pi^{\ast}\left(\RicP - \frac{\HessP(h)}{h}\right) \geq -m(G\circ\distP\circ\pi)\pi^{\ast}\sigma \\
	& = -m(G\circ\distP\circ\pi)\,\hat\sigma
	\end{align*}
	on $M$. If $G$ is nondecreasing, by (\ref{gammaineq}) we infer
	\begin{equation} \label{Gboundcomparison}
	- G\circ\distP\circ\pi \geq - G\circ\disthat
	\end{equation}
	on $M$ and therefore (\ref{modBEriccibound}) is satisfied. If $\pi$ is a diffeomorphism, then (\ref{gammaineq}) holds with equality sign even if $G$ is not monotonic and (\ref{modBEriccibound}) is again satisfied on $M$.
	
	Suppose that $\Bhat_r$ is relatively compact in $M$. From claim (c) of Lemma \ref{lemma:genBGthm} we have
	\[
	\int_{\Bhat_s}\hat h = \int_0^s\left[\int_{\partial \Bhat_t}\hat h\right] \di t = \int_0^s\frac{\int_{\partial \Bhat_t}\hat h}{k(t)^m} k(t)^m \di t \geq \frac{\int_{\partial \Bhat_s}\hat h}{k(s)^m}\int_0^s k(t)^m \di t
	\]
	and therefore
	\[
	\frac{\int_{\partial \Bhat_s}\hat h}{\int_{\Bhat_s}\hat h} \leq \frac{k(s)^m}{\int_0^s k(t)^m \di t}
	\]
	for a.e.~$s\in(0,r)$. Then (\ref{boundarycoshest2}) follows from (\ref{boundarycoshest1}).
	
	The other inequalities can be derived from (\ref{boundarycoshest2}) following the argument of the proof of Theorem \ref{thm:coshestimate1}.
\end{proof}

\noindent Theorem \ref{intro:thm:cosh1} of the Introduction follows from Theorem \ref{thm:coshestimate2}. Indeed we have

\begin{proof}[Proof of Theorem \ref{intro:thm:cosh1}]
	As pointed out in Remark \ref{rmk:ricci} in the Introduction as a consequence of (\ref{ricbar}), assumption (\ref{QBEriccibound}) of Theorem \ref{intro:thm:cosh1} implies (\ref{modBEricciboundP}) with $G\equiv G_0$ on $\RR^+_0$. For such $G$, a solution $k\in C^2(\RR^+_0)$ of problem (\ref{ksupersol}) is given by
	\[
	k(t)=\frac{1}{\sqrt{G_0}}\sinh(\sqrt{G_0}t).
	\]
	Since $M$ is connected and $H_0>0$, the function $H$ has constant sign, so $|\Hmean(\Bhat_r)|\geq\inf_{\Bhat_r}|H|\geq H_0$ for each $r>0$. By de l'H\^opital theorem we have
	\begin{align*}
	\lim_{r\to+\infty}\frac{\log\int_0^r k(t)^m \di t}{r} & = \lim_{r\to+\infty}\frac{k(r)^m}{\int_0^r k(t)^m \di t} = \lim_{r\to+\infty}\frac{mk'(r)}{k(r)} \\
	& = \lim_{r\to+\infty}m\sqrt{G_0}\coth(\sqrt{G_0}r)=m\sqrt{G_0}.
	\end{align*}
	Since $G$ is constant, therefore nondecreasing, we apply (\ref{limsupcoshest2}) to obtain
	\[
	H_0\leq\sqrt{G_0}\limsup_{M\ni x\to\infty}\sqrt{\cosh^2\theta(x)-1},
	\]
	or, equivalently,
	\[
	1+\frac{H_0^2}{G_0} \leq \limsup_{M\ni x\to+\infty}\cosh^2\theta(x).
	\]
	Then (\ref{limsupcosh}) follows.
\end{proof}

\begin{proof}[Proof of Corollary \ref{intro:crl:cosh1}]
	Set $H_0 = \inf_M |H|$ and suppose, by contradiction, $H_0 > 0$. For every real number $G_0>0$ and for each horizontal vector $X\in T\Mbar$ we have $\Ricbar(X,X) \geq 0 \geq - mG_0|X|^2$, so we can apply Theorem \ref{intro:thm:cosh1} to get
	\[
		\limsup_{M\ni x\to+\infty} \cosh\theta(x) \geq \sqrt{1+\frac{H_0^2}{G_0}} \qquad \text{for every } \, G_0 > 0.
	\]
	Since $H_0 > 0$, we have $\lim_{G_0 \to 0^+} \sqrt{1+H_0^2/G_0} = +\infty$, so we get
	\[
		\limsup_{M\ni x\to+\infty} \cosh\theta(x) = +\infty,
	\]
	contradicting the assumption that $\psi$ has bounded hyperbolic angle.
\end{proof}

\noindent Another consequence of Theorem \ref{thm:coshestimate1} is the following

\begin{thm} \label{thm:coshestimate3}
	Suppose that $(M,\hat\sigma)$ is complete and that $\pi:M\to\PP$ is a covering map of finite degree $d$. Then the following statements hold true.
	\begin{itemize}
		\item [(a)] Let $r>0$. Then
		\begin{equation} \label{disccoshest3}
		\max_{\overline{\Bhat_r}\setminus\Bhat_R} \sqrt{\cosh^2\theta-1}\,\frac{\log d + \log\int_{\BP_r} h - \log\int_{\Bhat_R}\hat h}{r-R} \geq \min_{R\leq\rho\leq r}m|\Hmean(\Bhat_{\rho})|
		\end{equation}
		for each $R\in(0,r)$.
		\item [(b)] If $Q:\mathbb R^+\to\mathbb R^+$ is such that
		\begin{equation} \label{growthbound3}
		\liminf_{r\to+\infty}\frac{\log d + \log\int_{\BP_r} h}{Q(r)} \leq 1,
		\end{equation}
		then
		\begin{equation} \label{limsupcoshest3}
		\limsup_{r\to+\infty}|\Hmean(\Bhat_r)| \leq \frac{1}{m}\left(\limsup_{M\ni x\to\infty}\sqrt{\cosh^2\theta(x)-1}\right)\left(\limsup_{r\to+\infty}\frac{Q(r)}{r}\right)
		\end{equation}
		with the same observation of Theorem \ref{thm:coshestimate1} about the RHS of (\ref{limsupcoshest3}).
	\end{itemize}
\end{thm}

\begin{proof}
	Let $r>0$. Since $(M,\hat\sigma)$ is complete, $\Bhat_r$ is relatively compact in $M$. By (\ref{ballinclusion}), we have $\pi(\Bhat_r)\subseteq \BP_r$ and therefore $\Bhat_r\subseteq \pi^{-1}(\pi(\Bhat_r)) \subseteq \pi^{-1}(\BP_r)$. Each $p\in\pi(\Bhat_r)$ has a connected open neighbourhood $U_p\subseteq\pi(\Bhat_r)$ such that, for each connected open neighbourhood $V$ of $p$ contained in $U_p$, the open subset $\pi^{-1}(V)\subseteq M$ has $d$ connected components and the restriction of $\pi$ to any of them is an isometry onto $U$. A Vitali covering argument shows that $\pi(\Bhat_r)$ can be covered, up to a subset of null measure, by a countable family $\{V_n\}_{n\in\mathbb N}$ of pairwise disjoint such open neighbourhoods of points $p_n$. For each $n\in\mathbb N$ we have
	\[
	\int_{\pi^{-1}(V_n)\cap \Bhat_r}\hat h \leq \int_{\pi^{-1}(V_n)}\hat h = d\int_{V_n}h.
	\]
	Since $\left(\cup_n\pi^{-1}(V_n)\right)\cap \Bhat_r = \pi^{-1}(\cup_n V_n)\cap \Bhat_r$ has zero measure in $\Bhat_r$, by the monotone convergence theorem we have
	\begin{equation} \label{findegcomparison}
	\int_{\Bhat_r}\hat h \leq \sum_{n=1}^{+\infty}\int_{\pi^{-1}(V_n)}\hat h \leq d\sum_{n=1}^{+\infty}\int_{V_n} h = d\int_{\pi(\Bhat_r)} h \leq d\int_{\BP_r} h.
	\end{equation}
	Hence, (\ref{disccoshest3}) follows from (\ref{disccoshest1}). Moreover, if (\ref{growthbound3}) is satisfied, then (\ref{growthbound1}) is also satisfied and therefore (\ref{limsupcoshest3}) follows from (\ref{limsupcoshest1}).
\end{proof}

\noindent Theorem \ref{intro:thm:cosh2} of the Introduction now easily follows.

\begin{proof}[Proof of Theorem \ref{intro:thm:cosh2}.]
	In Theorem \ref{thm:coshestimate3} choose $Q(r)=m\sqrt{G_0}r$. Then (\ref{volumeGbound}) of Theorem \ref{intro:thm:cosh2} implies (\ref{growthbound3}). Indeed,
	\[
		\liminf_{r\to+\infty}\frac{\log d + \log\int_{\BP_r} h}{Q(r)} = \lim_{r\to+\infty}\frac{\log d}{m\sqrt{G_0}r} + \liminf_{r\to+\infty}\frac{\log\int_{\BP_r} h}{m\sqrt{G_0}r} = 0 + \liminf_{r\to+\infty}\frac{\log\int_{\BP_r} h}{m\sqrt{G_0}r} \leq 1
	\]
	Arguing as in the proof of Theorem \ref{intro:thm:cosh1}, apply (\ref{limsupcoshest3}) to obtain the desired conclusion.
\end{proof}

\noindent Corollary \ref{intro:crl:cosh2} can be derived from Theorem \ref{intro:thm:cosh2} by reasoning as in the proof of Corollary \ref{intro:crl:cosh1}.

\section{A half-space theorem}

We reproduce some of the arguments of Section 3.3 of \cite{AMR}, with due differences, to prove a comparison principle at infinity for the Lorentzian mean curvature operator
\[
	Lu = \frac{1}{q}\,\div\left(\frac{q^2 Du}{\sqrt{1-q^2|Du|^2}}\right)
\]
on every end of a complete, noncompact Riemannian manifold $(\widetilde{M},\metric)$, where $q$ is a sufficiently regular positive function on $\widetilde{M}$. We recall that an end of a noncompact manifold $\widetilde{M}$ with respect to a compact subset $K \subseteq \widetilde{M}$ is any unbounded connected component of $\widetilde{M} \setminus K$. We start by stating and proving some preliminary lemmas.

\begin{lemma} \label{lemma:coerc}
	Let $(V,\langle,\rangle)$ be a real vector space with a given positive definite symmetric bilinear form. Let $X,Y\in V$ be such that $|X|,|Y|<1$, where $|X|=\sqrt{\langle X,X\rangle}$, $|Y|=\sqrt{\langle Y,Y\rangle}$. Then
	\[
	\left\langle\frac{X}{\sqrt{1-|X|^2}}-\frac{Y}{\sqrt{1-|Y|^2}},X-Y\right\rangle\geq0
	\]
	and the equality holds if and only if $X=Y$.
\end{lemma}

\begin{proof}
	If $|X|=|Y|$, then
	\[
	\left\langle\frac{X}{\sqrt{1-|X|^2}}-\frac{Y}{\sqrt{1-|Y|^2}},X-Y\right\rangle = \frac{|X-Y|^2}{\sqrt{1-|X|^2}}
	\]
	and the claim is proved. If $|Y|=t|X|$, $t\in[0,1)$, then
	\begin{align*}
	\left\langle\frac{X}{\sqrt{1-|X|^2}}-\frac{Y}{\sqrt{1-|Y|^2}},X-Y\right\rangle & = \frac{\langle X,X-Y\rangle}{\sqrt{1-|X|^2}} - \frac{\langle Y,X-Y\rangle}{\sqrt{1-t^2|X|^2}} \\
	& \geq \frac {(1-t)|X|^2}{\sqrt{1-|X|^2}} - \frac{(1-t)t|X|^2}{\sqrt{1-t^2|X|^2}} \\
	& = (1-t)|X|\left(\frac{|X|}{\sqrt{1-|X|^2}}-\frac{t|X|}{\sqrt{1-t^2|X|^2}}\right) \\
	& \geq 0,
	\end{align*}
	where the last inequality is strict unless $|X|=0$.
\end{proof}

\begin{lemma} \label{lemma:comparisonprinciple}
	Let $(\widetilde{M},\langle\;,\;\rangle)$ be a Riemannian manifold and let $\Omega\subseteq \widetilde{M}$ be a nonempty, connected, relatively compact open subset. Let $q$ be a measurable function on $\Omega$ such that $q > 0$ a.~e.~on $\Omega$ and let $u,v\in\Lip(\overline\Omega)$ be such that
	\begin{equation} \label{comparisonhp}
	\begin{cases}
	\esssup_{\overline\Omega} q|Du|<1, \\
	\esssup_{\overline\Omega} q|Dv|<1, \\[0.1cm]
	\div\left(\dfrac{q^2 Du}{\sqrt{1-q^2|Du|^2}}\right) \leq \div\left(\dfrac{q^2 Dv}{\sqrt{1-q^2|Dv|^2}}\right) & \text{weakly on } \Omega, \\
	u \geq v & \text{on } \partial\Omega.
	\end{cases}
	\end{equation}
	Then, $u\geq v$ on $\overline\Omega$.
\end{lemma}

\begin{proof}
	Let $\varepsilon>0$ be given. Let $\alpha\in C^1(\RR)$ be such that
	\[
	\begin{cases}
	\alpha(t) = 0 & \text{for each } t\geq-\varepsilon, \\
	\alpha'(t) < 0 & \text{for each } t<-\varepsilon.
	\end{cases}
	\]
	Consider a vector field $W$ on $\Omega$ such that
	\[
	\begin{cases}
	W = [\alpha\circ(u-v)]\left(\dfrac{q^2 Du}{\sqrt{1-q^2|Du|^2}}-\dfrac{q^2 Dv}{\sqrt{1-q^2|Dv|^2}}\right) & \text{a.~e.~on } \Omega, \\
	W \equiv 0 & \text{where } u\geq v-\varepsilon.
	\end{cases}
	\]
	Because of (\ref{comparisonhp}), the vector field $W$ is compactly supported in $\Omega$ and its weak divergence satisfies
	\begin{align*}
	\div W & \leq \left\langle\frac{q^2 Du}{\sqrt{1-q^2|Du|^2}}-\frac{q^2 Dv}{\sqrt{1-q^2|Dv|^2}},D[\alpha\circ(u-v)]\right\rangle \\
	& = [\alpha'\circ(u-v)]\left\langle\frac{q Du}{\sqrt{1-q^2|Du|^2}}-\frac{q Dv}{\sqrt{1-q^2|Dv|^2}},q Du - q Dv\right\rangle\leq0
	\end{align*}
	weakly on $\Omega$, where the last inequality follows from Lemma \ref{lemma:coerc}. By applying the divergence theorem to $W$ on an open subset $\tilde\Omega\subseteq\Omega$ with smooth boundary and such that $\supp W\subseteq\tilde\Omega$ we get
	\[
	0 = \int_{\tilde\Omega}\div W = \int_{\Omega}\div W.
	\]
	It follows that
	\[
	[\alpha'\circ(u-v)]\left\langle\frac{q Du}{\sqrt{1-q^2|Du|^2}}-\frac{q Dv}{\sqrt{1-q^2|Dv|^2}},q Du - q Dv\right\rangle = 0
	\]
	a.~e.~on $\Omega$. From Lemma \ref{lemma:coerc} again we deduce that $Du=Dv$ a.~e.~where $u<v-\varepsilon$. Since $\varepsilon>0$ is arbitrarily given, it follows that the Lipschitz function $(u-v)_-$ has almost everywhere vanishing gradient on $\Omega$, and therefore it is constant on $\overline\Omega$. From (\ref{comparisonhp}) it follows that it is identically zero on $\overline\Omega$ and the claim is proved.
\end{proof}

\begin{lemma} \label{lemma:strongmax}
	Let $(\widetilde{M},\langle,\rangle)$ be a Riemannian manifold, let $\Omega\subseteq \widetilde{M}$ be a nonempty, connected, relatively compact open subset. Let $q\in C^1(\overline\Omega)$, $q>0$ and let $u\in C^2(\overline\Omega)$ be such that
	\[
	\begin{cases}
	\max_{\overline\Omega}q|Du| < 1, \\[0.1cm]
	\div\left(\dfrac{q^2 Du}{\sqrt{1-q^2|Du|^2}}\right) \leq 0 & \text{on } \Omega, \\
	u \geq 0 & \text{on } \Omega.
	\end{cases}
	\]
	If $u(x)=0$ for some $x\in\Omega$, then $u\equiv 0$ on $\overline\Omega$.
\end{lemma}

\begin{proof}
	Set
	\[
	f = \frac{q^2}{\sqrt{1-q^2|Du|^2}} \in C^1(\overline\Omega).
	\]
	We have that $u$ satisfies
	\[
	\mathrm{\Delta} u + \frac{1}{f}\langle Df,Du\rangle \leq 0
	\]
	on $\Omega$. The strong maximum principle for the operator $\mathrm{\Delta}+\frac{1}{f}\langle Df,D\;\cdot\;\rangle$ (see \cite{AMR}, Theorem 3.10) yields the desired conclusion.
\end{proof}

We are now ready to state and prove the following technical analytical result, relating the unboundedness of a solution of the prescribed mean curvature equation to the existence of an appropriate function which acts as a potential function for the equation outside a bounded set.

\begin{thm} \label{thm:Kcomparison}
	Let $(\widetilde{M},\metric)$ be a complete, noncompact Riemannian manifold, $K\subseteq \widetilde{M}$ a nonempty compact subset and let $(M_K,\hat\sigma)$ be a connected component of $\widetilde{M}\setminus K$ endowed with the restriction $\hat\sigma$ of the metric $\metric$. For each $r>0$, let
	\[
		\Omega_r = \{ x \in M_K : \mathrm{dist}_{\metric}(x,K) < r \}
	\]
	and let $\partial \Omega_r$ denote the boundary of $\Omega_r$ with respect to the induced topology on $M_K$, that is,
	\[
		\partial \Omega_r = \{ x \in M_K : \mathrm{dist}_{\metric}(x,K) = r \}.
	\]
	Let $q$ be measurable function such that $q>0$ a.~e.~on $M_K$, let $\mathcal H \in L^1_{\loc}(M_K)$, and let $\tau \in \Lip_{\loc}(M_K)$ be a solution of inequality
	\begin{equation} \label{Hineq}
		\div_{\hat\sigma}\left(\frac{q^2 D\tau}{\sqrt{1-q^2|D\tau|^2}}\right) \geq \mathcal H \qquad \text{on } \, M_K
	\end{equation}
	satisfying
	\[
		\esssup_{\overline\Omega} q|D\tau| < 1
	\]
	for every relatively compact $\Omega \subseteq M_K$. Let $r > R > 0$ be given real numbers and suppose that there exists a function $u \in\Lip_{\loc}(M_K\setminus\Omega_R)$ satisfying
	\begin{equation} \label{Kpot}
		\begin{cases}
			u \geq 0 & \text{on } \partial \Omega_R, \\
			u < \displaystyle \max_{\partial\Omega_r} \tau - \max_{\partial\Omega_R} \tau & \text{on } \partial \Omega_r, \\
			u(x) \to +\infty & \text{as } \mathrm{dist}_{\metric}(x,K) \to+\infty, \, x \in M_K, \\
			\esssup_{\overline\Omega} q|Du| < 1 & \text{on every relatively compact } \Omega \subseteq M_K, \\[0.1cm]
			\div_{\hat\sigma}\left(\dfrac{q^2 Du}{\sqrt{1-q^2|Du|^2}}\right) \leq \mathcal H & \text{on } M_K\setminus\overline{\Omega_R}. \\
		\end{cases}
	\end{equation}
	Then
	\[
	\sup_{M_K} \tau = +\infty.
	\]
\end{thm}

\begin{proof}
	Suppose, by contradiction, that $\sup_{M_K} \tau < +\infty$. Let
	\[
		\varepsilon = \max_{\partial\Omega_r} \tau - \max_{\partial\Omega_R} \tau - \max_{\partial\Omega_r} u \qquad \text{and} \qquad v = \tau - \max_{\partial \Omega_R} \tau - \frac{\varepsilon}{2}.
	\]
	Note that $\varepsilon>0$ by the second assumption in \eqref{Kpot}. Since $v$ and $\tau$ only differ by an additive constant, we also have
	\[
		\varepsilon = \max_{\partial\Omega_r} v - \max_{\partial\Omega_R} v - \max_{\partial\Omega_r} u = \max_{\partial\Omega_r} v - \max_{\partial\Omega_R} \left(\tau - \max_{\partial \Omega_R} \tau - \frac{\varepsilon}{2}\right) - \max_{\partial\Omega_r} u = \max_{\partial\Omega_r} v - \max_{\partial\Omega_r} u + \frac{\varepsilon}{2}
	\]
	and therefore
	\[
		\max_{\partial\Omega_R} v = - \frac{\varepsilon}{2} < 0 \leq \min_{\partial\Omega_R} u, \quad \max_{\partial\Omega_r} v = \max_{\partial\Omega_r} u + \frac{\varepsilon}{2} > \max_{\partial\Omega_r} u, \quad \displaystyle \sup_{M_K} v = \sup_{M_K} \tau - \max_{\partial \Omega_R} \tau - \frac{\varepsilon}{2} < +\infty.
	\]\\
	This implies that
	\[
		\begin{cases}
			u - v > 0 & \text{on } \partial \Omega_R, \\
			u - v < 0 & \text{on a nonempty subset of } \partial \Omega_r, \\
			u(x) - v(x) \to + \infty & \text{as } \mathrm{dist}_{\metric}(x,K) \to+\infty, \, x \in M_K.
		\end{cases}
	\]
	So, the subset $\{ x \in M_K : u(x) < v(x) \}$ is nonempty and relatively compact in $M_K$. Let $\Omega$ be one of its connected components. We have $u = v$ on $\partial\Omega$ and
	\[
		\div\left(\frac{q^2 Du}{\sqrt{1-q^2|Du|^2}}\right) \leq \div\left(\frac{q^2 Dv}{\sqrt{1-q^2|Dv|^2}}\right) \qquad \text{on } \, \Omega.
	\]
	By Lemma \ref{lemma:comparisonprinciple}, we conclude that $u \geq v$ in $\Omega$, contradiction.
\end{proof}

As a direct application of Theorem \ref{thm:Kcomparison} we can prove Theorems \ref{intro:thm:hs2} and \ref{intro:thm:schwarz} from the Introduction. The following two lemmas are instrumental to guaranteeing, under the hypotheses of the theorems, the existence of suitable nonlinear potentials $u$ satisfying conditions \eqref{Kpot}.

\begin{lemma} \label{lemma:Kexistenceprod0}
	Let $(\widetilde{M},\metric)$ be a complete, noncompact Riemannian manifold of dimension $m$ with a fixed origin $o \in \widetilde{M}$ and let $\dist$ be the distance function from $o$ in the metric $\metric$. Suppose that the Ricci tensor $\Ric$ of $\widetilde{M}$ satisfies
	\begin{equation} \label{ricbound}
		\Ric \geq -(m-1)(G\circ\dist) \, \hat\sigma
	\end{equation}
	on $M$ for some nonnegative function $G\in C^0(\RR^+_0)$ and let $k\in C^2(\RR^+_0)$ be a solution of the problem
	\[
		\begin{cases}
			k'' - G k \geq 0, \\
			k(0) = 0, \quad k'(0) = 1.
		\end{cases}
	\]
	Let $R>0$, $A\in C^0([R,+\infty))$, $A>0$ be such that
	\begin{equation} \label{limA}
		\liminf_{s\to+\infty} \frac{1}{k(s)^{m-1}}\int_R^s A(t)k(t)^{m-1} \di t > 0.
	\end{equation}
	Then for each $r > R$, $\varepsilon > 0$ there exists a function
	\[
		u\in C^2(\widetilde{M}\setminus(\overline{B_R}\cup\cut(o)))\cap\Lip_{\loc}(\widetilde{M}\setminus B_R),
	\]
	with $B_s$ the geodesic ball of $(\widetilde{M},\metric)$ centered at $o$ with radius $s$, satisfying
	\begin{equation} \label{Kpotp}
		\begin{cases}
			u = 0 & \text{on } \partial B_R, \\
			u \leq \varepsilon & \text{on } \partial B_r, \\
			u(x) \to +\infty & \text{as } \gamma(x) \to+\infty, \\
			\esssup_{\overline\Omega} |Du| < 1 & \text{on every relatively compact } \Omega \subseteq \widetilde{M}, \\
			\div\left(\dfrac{Du}{\sqrt{1-|Du|^2}}\right) \leq A(\dist) & \text{on } \widetilde{M} \setminus\overline{B_R}. \\
		\end{cases}
	\end{equation}
\end{lemma}

\begin{proof}
	We prove the lemma by explicitely constructing, for any given $r > R$, $\varepsilon > 0$, a radial function $u = u_0(\dist)$, with $u_0 \in C^2([R,+\infty))$, satisfying \eqref{Kpotp}. From \eqref{ricbound} and by the comparison theorem for the Laplacian of the distance function (see \cite{AMR}) we have that
	\[
		\Delta \dist \leq (m-1) \frac{k'(\dist)}{k(\dist)}
	\]
	pointwise on $M\setminus\cut(o)$ and weakly on $M$, so if $u = u_0(\dist)$ for some
	\[
		u_0 \in C^2([R,+\infty)) \qquad \text{satisfying} \qquad 0 \leq u_0' < 1 \, \text{ on } \, [R,+\infty),
	\]
	it follows that
	\[
		\div\left(\frac{Du}{\sqrt{1-|Du|^2}}\right) = \frac{u_0'(\dist)}{\sqrt{1-u_0'(\dist)^2}}\Delta\dist + \left(\frac{u_0'}{\sqrt{1-(u_0')^2}}\right)'(\gamma) \leq (m-1)\frac{k'(\dist)}{k(\dist)}f(\dist) + f'(\dist)
	\]
	pointwise on $M\setminus(\overline{B_R}\cup\cut(o))$ and weakly on $M\setminus\overline{B_R}$, having set
	\[
		f = \frac{u_0'}{\sqrt{1-(u_0')^2}}, \qquad \text{or, equivalently,} \qquad u_0' = \frac{f}{\sqrt{1+f^2}}.
	\]
	So, we aim at obtaining $u_0 \in C^2([R,+\infty))$ satisfying
	\begin{equation} \label{Kpotu0}
		u_0(R) = 0, \qquad u_0(r) \leq \varepsilon, \qquad \lim_{s\to+\infty} u_0(s) = +\infty 
	\end{equation}
	and we look for $u_0$ in the form
	\begin{equation} \label{u0int}
		u_0(s) = \int_R^s \frac{f(t)}{\sqrt{1+f(t)^2}} \di t
	\end{equation}
	with $f \in C^1([R,+\infty))$ a nonnegative function such that
	\begin{equation} \label{fineq}
		f' + (m-1)\frac{k'}{k}f \leq A \qquad \text{on } \, (R,+\infty).
	\end{equation}
	
	Note that for every $C\in\RR$ the unique solution $f_C$ of the Cauchy problem
	\[
		\begin{cases}
			f_C' + (m-1)\dfrac{k'}{k}f_C = CA & \qquad \text{on } \, (R,+\infty), \\
			f_C(R) = 0 &
		\end{cases}
	\]
	is given by
	\[
		f_C(s) = \frac{C}{k(s)^{m-1}}\int_R^s A(t)k(t)^{m-1} \di t \qquad \text{for each } s \geq R
	\]
	and is of class $C^1([R,+\infty))$.	If $C\in(0,1]$, then $f = f_C$ satisfies \eqref{fineq} by positivity of $A$, and assumption \eqref{limA} implies that
	\[
		\liminf_{s\to+\infty} f_C(s) > 0, \qquad \text{so} \qquad \liminf_{s\to+\infty} \frac{f_C(s)}{\sqrt{1+f_C(s)^2}} > 0.
	\]
	The function $u_0$ defined as in \eqref{u0int} for $f = f_C$ satisfies
	\[
		\lim_{s\to+\infty} u_0'(s) > 0, \qquad \text{and therefore} \qquad \lim_{s\to+\infty} u_0(s) = +\infty,
	\]
	and also
	\begin{equation} \label{u0r}
		u_0(r) \leq \int_R^r f(t) \di t \leq C (r-R) \max_{s\in[R,r]} \frac{1}{k(s)^{m-1}}\int_R^s A(t)k(t)^{m-1} \di t.
	\end{equation}
	For any fixed $r>R$, $\varepsilon>0$, it is always possible to choose $C\in(0,1]$ small enough so that the RHS of \eqref{u0r} is not larger than $\varepsilon$. Hence, for a suitable choice of $C\in(0,1]$, the function $u_0$ given by \eqref{u0int} with $f = f_C$ satisfies conditions \eqref{Kpotu0} and therefore $u = u_0(\dist)$ satisfies all of the requirements in \eqref{Kpotp}.
\end{proof}

\begin{lemma} \label{lemma:Kexistencerad}
	Let $(\PP_0,\sigma_0)$ be a radially symmetric Riemannian manifold of dimension $m$, with
	\[
		\PP_0 = (\rho_0,+\infty) \times \mathbb{S}^{m-1}, \qquad \sigma_0 = \frac{\di \rho^2}{V(\rho)} + \rho^2 \metric_{\mathbb{S}^{m-1}},
	\]
	$\rho_0 \geq 0$, $\rho$ the standard coordinate on $(\rho_0,+\infty)$, $(\mathbb{S}^{m-1},\metric_{\mathbb{S}^{m-1}})$ the standard $(m-1)$-dimensional sphere and $V\in C^{\infty}((\rho_0,+\infty))$ a positive function such that for some (hence, any) $\varepsilon>0$
	\[
		\int_{\rho_0}^{\rho_0+\varepsilon} \frac{\di t}{\sqrt{V(t)}} < +\infty, \qquad \int_{\rho_0+\varepsilon}^{+\infty} \frac{\di t}{\sqrt{V(t)}} = +\infty,
	\]
	so that the bijection $\phi : (\rho_0,+\infty) \to (0,+\infty)$ given by
	\[
		\phi(\rho) = \int_{\rho_0}^{\rho} \frac{\di t}{\sqrt{V(t)}} \qquad \text{for each } \, \rho > \rho_0
	\]
	is well defined. Let $\rho_1>\rho_0$, $h_0 \in C^1([\rho_1,+\infty))$ be such that
	\begin{equation} \label{limhr}
		\int_{\rho_1}^{+\infty} \frac{\di t}{h_0(t)\sqrt{V(t)}} = +\infty
	\end{equation}
	and let $A_0\in C^0([\rho_1,+\infty))$, $A_0>0$, be such that
	\begin{equation} \label{limAr}
		\int_{\rho_1}^{+\infty} \frac{A_0(t)t^{m-1}}{\sqrt{V(t)}} \di t = +\infty, \qquad \liminf_{\rho\to+\infty} \frac{1}{h_0(\rho) \rho^{m-1}} \int_{\rho_1}^{\phi(\rho)} \frac{A_0(t)t^{m-1}}{\sqrt{V(t)}} \di t > 0.
	\end{equation}
	Then for each $\rho_2 > \rho_1$, $\beta \in \RR$ there exists a function
	\[
		u \in C^2( [\rho_1,+\infty) \times \mathbb{S}^{m-1} )
	\]
	satisfying
	\begin{equation} \label{Kpotr}
		\begin{cases}
			u = 0 & \text{on } \{\rho_1\} \times \mathbb{S}^{m-1}, \\
			u \leq \beta & \text{on } \{\rho_2\} \times \mathbb{S}^{m-1}, \\
			u(x) \to +\infty & \text{as } \rho(x) \to+\infty \\
			h_0(\rho)|Du| < 1 & \text{on } [\rho_1,+\infty) \times \mathbb{S}^{m-1}, \\
			\div\left(\dfrac{h_0(\rho)^2 Du}{\sqrt{1-h_0(\rho)^2|Du|^2}}\right) \leq A_0(\rho) & \text{on } [\rho_1,+\infty) \times \mathbb{S}^{m-1}. \\
		\end{cases}
	\end{equation}
\end{lemma}

	\begin{proof}
		The proof is analogous to that of Lemma \ref{lemma:Kexistenceprod0}. We first observe that $(\PP_0,\metric)$ is isometric to the product manifold
		\begin{equation} \label{Prad}
			\PP = (0,+\infty) \times \mathbb{S}^{m-1} \quad \text{with metric} \quad \metric_{\PP} = \di s^2 + g(s)^2 \metric_{\mathbb{S}^{m-1}},
		\end{equation}
		where $s$ is the standard coordinate on $(0,+\infty)$ and $g : (0,+\infty) \to (\rho_0,+\infty)$ is the inverse of the function $\phi : (\rho_0,+\infty) \to (0,+\infty)$. Indeed, an isometry between $\PP_0$ and $\PP$ is given by the map $\Phi : \PP_0 \to \PP$ defined by
		\[
			\Phi((\rho,p)) = (\phi(\rho),p) \qquad \text{for each } \, \rho > \rho_0, \, p \in \mathbb{S}^{m-1}.
		\]
		Let $\rho_2 > \rho_1 > \rho_0$ and $h_0$, $A_0$ be given as in the statement of the lemma. Setting $R = g(\rho_1)$, $r = g(\rho_2)$, we let $A = A_0\circ g \in C^0([R,+\infty))$, $h = h_0\circ g \in C^1([R,+\infty))$. Clearly, the thesis is equivalent to claiming the existence of a function $v \in C^2( [R,+\infty) \times \mathbb{S}^{m-1} )$ satisfying
		\begin{equation} \label{Kpotrv}
			\begin{cases}
				v = 0 & \text{on } \{R\} \times \mathbb{S}^{m-1}, \\
				v \leq \beta & \text{on } \{r\} \times \mathbb{S}^{m-1}, \\
				v(x) \to +\infty & \text{as } s(x) \to+\infty \\
				h(s)|Dv| < 1 & \text{on } [R,+\infty) \times \mathbb{S}^{m-1}, \\
				\div\left(\dfrac{h(s)^2 Dv}{\sqrt{1-h(s)^2|Dv|^2}}\right) \leq A(s) & \text{on } [R,+\infty) \times \mathbb{S}^{m-1}. \\
			\end{cases}
		\end{equation}
		We construct $v$ as a radial function $v=v_0(s)$, with $v_0 \in C^2([R,+\infty))$. For a radially symmetric manifold as in \eqref{Prad}, the Laplacian of the coordinate function $s$ is given by
		\[
			\Delta s = (m-1)\frac{g'(s)}{g(s)}.
		\]
		Following the lines of the proof of Lemma \ref{lemma:Kexistenceprod0}, we look for $v_0$ of the form
		\begin{equation} \label{v0int}
			v_0(s) = \int_R^s \frac{1}{h(t)}\frac{f(t)/h(t)}{\sqrt{1+f(t)^2/h(t)^2}} \di t
		\end{equation}
		for some $f\in C^1([R,+\infty))$, since this is equivalent to saying that
		\[
			v_0(R) = 0 \qquad \text{and} \qquad f = \frac{h^2 v_0'}{\sqrt{1-h^2 \cdot (v_0')^2}} \qquad \text{on } \, (R,+\infty),
		\]
		and the last inequality in \eqref{Kpotrv} can be restated as
		\[
			f' + (m-1)\frac{g'}{g}f \leq A \qquad \text{on } \, (R,+\infty). \Rightarrow (g^{m-1}f)' = C g^{m-1}A
		\]
		For every $C \in (0,1]$, $\beta_1 \in \RR$, the function $f_{C,\beta_1} \in C^1([R,+\infty))$ defined by
		\[
			f_{C,\beta_1} (s) = \frac{C}{g(s)^{m-1}} \int_R^s A(t)g(t)^{m-1} \di t + \frac{\beta_1}{g(s)^{m-1}} \qquad \text{for each } \, s \geq R
		\]
		satisfies
		\[
			f_{C,\beta_1}' + (m-1)\frac{g'}{g}f_{C,\beta_1} = CA \leq A \qquad \text{on } \, (R,+\infty)
		\]
		because $A>0$, since $A_0$ is positive. To ensure that $v_0$ defined as in \eqref{v0int} with $f = f_{C,\beta_1}$ also satisfies
		\[
			\lim_{s\to+\infty} v_0(s) = +\infty
		\]
		it is sufficient to have
		\[
			0 < \liminf_{s\to+\infty} \frac{f(s)}{h(s)} = \liminf_{s\to+\infty}  \frac{1}{h(s)g(s)^{m-1}} \left( C \int_R^s A(t)g(t)^{m-1} \di t + \beta_1 \right)
		\]
		and
		\[
			\int_R^{+\infty} \frac{\di s}{h(s)} = +\infty.
		\]
		By changing variables, these conditions can be restated as
		\begin{align*}
			\liminf_{\rho\to+\infty} \frac{1}{h_0(\rho) \rho^{m-1}}\left(C \int_{\rho_1}^{\phi(\rho)} \frac{A_0(t)t^{m-1}}{\sqrt{V(t)}} \di t + \beta_1 \right) > 0, \qquad \int_{\rho_1}^{+\infty} \frac{\di t}{h_0(t)\sqrt{V(t)}} = +\infty.
		\end{align*}
		and they are clearly satisfied for every $C\in(0,1]$, $\beta_1\in\RR$ under assumptions \eqref{limhr} and \eqref{limAr}. It only remains to guarantee that
		\[
			\beta \geq v_0(r) = \frac{1}{h(r)g(r)^{m-1}} \left( C \int_R^r A(t)g(t)^{m-1} \di t + \beta_1 \right).
		\]
		For any fixed value of $C\in(0,1]$, the last member of the previous expression diverges to $-\infty$ as $\beta_1 \to -\infty$, so it is always possible to choose $\beta_1 \in \RR$ so that the inequality is satisfied.
	\end{proof}

Theorem \ref{intro:thm:hs2} of the Introduction now follows at once.

\begin{proof}[Proof of Theorem \ref{intro:thm:hs2}]
	Let $o\in M$ be a fixed point, $\hat\sigma = \pi^{\ast}\sigma$ and let $\disthat$ be the distance function from $o$ in $(M,\hat\sigma)$. We prove the validity of statement $(a)$ of the theorem. The proof of statement $(b)$ is analogous, up to replacing the function $\tau$ in the argument below with the function $-\tau$. So, let us suppose that $\liminf_{M\ni x \to \infty} H(x) > 0$. Since $(\PP,\sigma)$ is complete, $(M,\hat\sigma)$ is also complete. So, the geodesic balls $\Bhat_r$ of $(M,\hat\sigma)$ centered at $o$ are relatively compact and therefore there exists $R_0 > 0$ such that $\inf_{M\setminus\overline{\Bhat_{R_0}}} H > 0$. Let $H_0 = \inf_{M\setminus\overline{\Bhat_{R_0}}} H$.
	
	We first apply Lemma \ref{lemma:Kexistenceprod0} with $(\widetilde{M},\metric) = (M,\hat\sigma)$. Assumption \eqref{riccibound'} implies \eqref{ricbound} with $G\equiv G_0$ on $\RR^+_0$. A solution $k\in C^2(\RR^+_0)$ of the problem
	\[
		\begin{cases}
			k''-Gk\geq0, \\
			k(0)=0, \quad k'(0)=1
		\end{cases}
	\]
	is given by
	\begin{equation} \label{solk}
		k(t)=\frac{1}{\sqrt{G_0}}\sinh(\sqrt{G_0}t).
	\end{equation}
	By de l'H\^opital theorem, for $k$ as in \eqref{solk} it follows that
		\[
			\lim_{r\to+\infty}\frac{mH_0}{k(r)^{m-1}}\int_R^r k(t)^{m-1} \di t = mH_0\cdot\lim_{r\to+\infty}\frac{k(r)}{(m-1)k'(r)}= \frac{mH_0}{(m-1)\sqrt{G_0}}
		\]
	Therefore, setting $A\equiv mH_0$ on $[R,+\infty)$, we have that \eqref{limA} is satisfied. By Lemma \ref{lemma:Kexistenceprod0}, for every $R>0$, $r>R$, $\varepsilon>0$ there exists a function $u\in C^2(M\setminus(\overline{\Bhat_{R_0+R}}\cup\cut(o)))\cap\Lip_{\loc}(M\setminus \Bhat_{R_0+R})$ such that
	\begin{equation} \label{Kpotu:hs}
		\begin{cases}
			u = 0 & \text{on } \partial \Bhat_{R_0+R}, \\
			u \leq \varepsilon & \text{on } \partial \Bhat_{R_0+r}, \\
			u(x) \to +\infty & \text{as } \disthat(x) \to+\infty, \\
			\esssup_{\overline\Omega} |Du| < 1 & \text{on every relatively compact } \Omega \subseteq M, \\
			\div\left(\dfrac{Du}{\sqrt{1-|Du|^2}}\right) \leq mH_0 & \text{on } M \setminus\overline{\Bhat_{R_0+R}}. \\
		\end{cases}
	\end{equation}
	
	We now argue by contradiction. Recall that the vertical height function $\tau$ of the immersion $\psi$ satisfies equation
	\[
		\div_{\hat\sigma} \left( \frac{D\tau}{\sqrt{1-|D\tau|^2}} \right) = mH \qquad \text{on } \, M.
	\]
	Suppose that $\tau^{\ast} = \sup_M \tau < +\infty$. Since $H \not\equiv 0$ on $M$, the function $\tau$ cannot be constant on $M$. Since $H\geq0$ on $M$, Lemma \ref{lemma:strongmax} implies that the nonnegative function $u = \tau^{\ast} - \tau$ cannot attain the zero value, that is, the value $\tau^{\ast}$ cannot be attained by $\tau$ at any point of $M$. More generally, $\tau$ cannot attain a local maximum at any point of $M$. We now apply Theorem \ref{thm:Kcomparison} with $(\widetilde{M},\metric) = (M,\hat\sigma)$, $K = \overline{\Bhat_{R_0}}$, $M_K$ a connected component of $M\setminus\overline{\Bhat_{R_0}}$ such that
	\[
		\tau^{\ast} = \sup_{M_K} \tau,
	\]
	$q \equiv 1$, $\mathcal H \equiv mH_0$. Note that for every $s>0$ the subset $\Omega_s$ defined as in the statement of Theorem \ref{thm:Kcomparison} is the connected component of $\Bhat_{R_0+s} \setminus \overline{\Bhat_{R_0}}$ contained in $M_K$ and therefore $\partial \Omega_s \subseteq \partial \Bhat_{R_0+s}$. Choose $R>0$ and $r>R$ such that, setting
	\[
		\tau_1 = \max_{\partial \Omega_R} \tau, \qquad \tau_2 = \max_{\partial \Omega_r} \tau,
	\]
	we have
	\[
		\tau_1 < \tau_2.
	\]
	Note that such an $r>R$ exists because $\sup_{M_K} \tau = \tau^{\ast} > \tau_1$. Let $u\in C^2(M\setminus(\overline{\Bhat_{R_0+R}}\cup\cut(o)))\cap\Lip_{\loc}(M\setminus \Bhat_{R_0+R})$ be a function satisfying \eqref{Kpotu:hs} with $0 < \varepsilon < \tau_2 - \tau_1$. Then, conditions \eqref{Kpot} in Theorem \ref{thm:Kcomparison} are satisfied and we conclude that $\tau^{\ast} = \sup_{M_K} \tau = +\infty$, contradiction.
	
	So, we have proved that $\sup_M\tau=+\infty$. Therefore, $\pi_{\RR}(\psi(M))=\tau(M)\subseteq\RR$ is not contained in any upper bounded interval of the form $(-\infty,t_0]$, $t_0\in\RR$ and so $\psi(M)$ is not contained in any lower half-space of the form $\PP\times(-\infty,t_0]=\pi_{\RR}^{-1}((-\infty,t_0])$, $t_0\in\RR$.
\end{proof}

\begin{proof}[Proof of Corollary \ref{intro:crl:hs1-2}]
	By assumption, the function $H$ does not change sign on $M$ and $\psi(M)$ is contained in a slab $\PP\times[t_0,t_1]$, so it is also true that $\psi(M) \subseteq \PP \times (-\infty,t_1]$ and $\psi(M) \subseteq \PP \times [t_0,+\infty)$. If $H\geq 0$ on $M$, then by statement $(a)$ of Theorem \ref{intro:thm:hs2} we conclude that $\liminf_{M\ni x \to \infty} H(x) \leq 0$, so it must be $\liminf_{M\ni x \to \infty} H(x) = 0$. Similarly, if $H\leq 0$ on $M$ then by statement $(b)$ we conclude that $\limsup_{M\ni x \to \infty} H(x) = 0$.
\end{proof}

To conclude this section, we prove Theorem \ref{intro:thm:schwarz} along the lines of the proof of Theorem \ref{intro:thm:hs2}.

\begin{proof}[Proof of Theorem \ref{intro:thm:schwarz}]
	As recalled in the Introduction, the Schwarzschild spacetime $\Mbar$ of dimension $m+1$ has the structure of a standard static spacetime $\Mbar = \PP_S \times_{h_S} \RR$ with
	\[
		(\PP_S,\sigma_S) = \left( (\rho_S,+\infty) \times \mathbb{S}^{m-1}, \frac{\di \rho^2}{V(\rho)} + \rho^2 \metric_{ \mathbb{S}^{m-1} } \right), \qquad V(\rho) = 1 - 2\mu\rho^{2-m}, \qquad h_S = \sqrt{V(\rho)}, 
	\]
	where $\rho$ is the standard coordinate on the interval $(\rho_S,+\infty)$, $\mu>0$ is a mass parameter and $\rho_S = (2\mu)^{1/(m-2)}$.
	
	We prove the validity of statement $(a)$ of the theorem, the proof of statement $(b)$ being analogous. By assumption, there exists $H_0>0$ and $\rho_1 > \rho_S$ such that the mean curvature function $H$ of the immersion satisfies $H \geq H_0$ on $\{x \in M : \rho(\pi(x)) \geq \rho_1 \}$. Note that $\PP_S$ is simply connected because $m\geq 3$. Since $\pi : M \to \PP_S$ is a covering map, we deduce that $\pi$ is in fact a diffeomorphism. We set $(\PP_0,\sigma_0) = (M,\pi^{\ast}\sigma_S)$.
	
	We want to apply Lemma \ref{lemma:Kexistencerad} with $\rho_0 = \rho_S$, $\rho_1$ and $V$ as above, $h_0 = h_S$ and $A_0 = mH_0 h_0(\rho)$, so we verify that the hypotheses of the lemma are satisfied. Note that
	\[
		\sqrt{V(\rho)} = \sqrt{1 - \left(\frac{\rho}{\rho_S}\right)^{2-m}} = \sqrt{1 - \left(1 + \frac{\rho - \rho_S}{\rho_S}\right)^{2-m}} = \sqrt{(m-2)\frac{\rho-\rho_S}{\rho_S} + o(\rho-\rho_S)}
	\]
	as $\rho \to \rho_S^+$, and
	\[
		h_0(\rho) = \sqrt{V(\rho)} \to 1 \qquad \text{as } \, \rho \to +\infty,
	\]
	so we have
	\[
		\int_{\rho_0}^{\rho_0+1} \frac{\di t}{\sqrt{V(t)}} < +\infty, \quad \int_{\rho_0+1}^{+\infty} \frac{\di t}{\sqrt{V(t)}} = +\infty, \quad \int_{\rho_1}^{+\infty} \frac{\di t}{h_0(t)\sqrt{V(t)}} = \int_{\rho_1}^{+\infty} \frac{\di t}{V(t)} = +\infty.
	\]
	Moreover, since $\phi'(\rho) = 1/\sqrt{V(\rho)} \to 1$ as $\rho\to+\infty$, we have $\phi(\rho) \sim \rho$ as $\rho \to +\infty$ and therefore
	\[
		\int_{\rho_1}^{+\infty} \frac{A_0(t)t^{m-1}}{\sqrt{V(t)}} \di t = +\infty, \qquad \liminf_{\rho\to+\infty} \frac{1}{h_0(\rho) \rho^{m-1}} \int_{\rho_1}^{\phi(\rho)} \frac{A_0(t)t^{m-1}}{\sqrt{V(t)}} \di t = \lim_{\rho\to+\infty} H_0 \rho = +\infty.
	\]
	Then, we can apply Lemma \ref{lemma:Kexistencerad} to obtain that for each $\rho_2 > \rho_1$, $\beta \in \RR$ there exists a function
	\[
	u \in C^2( [\rho_1,+\infty) \times \mathbb{S}^{m-1} )
	\]
	satisfying
	\begin{equation} \label{Kpotr:sch}
		\begin{cases}
			u = 0 & \text{on } \{\rho_1\} \times \mathbb{S}^{m-1}, \\
			u \leq \beta & \text{on } \{\rho_2\} \times \mathbb{S}^{m-1}, \\
			u(x) \to +\infty & \text{as } \rho(x) \to+\infty \\
			h_0(\rho)|Du| < 1 & \text{on } [\rho_1,+\infty) \times \mathbb{S}^{m-1}, \\
			\div\left(\dfrac{h_0(\rho)^2 Du}{\sqrt{1-h_0(\rho)^2|Du|^2}}\right) \leq m H_0 h_0(\rho) & \text{on } [\rho_1,+\infty) \times \mathbb{S}^{m-1}.
		\end{cases}
	\end{equation}
	
	We conclude by applying Theorem \ref{thm:Kcomparison}. Fix $\rho_{\ast} \in (\rho_S,\rho_1)$ and let $R_0 = \phi(\rho_{\ast})$. By changing variables as in the proof of Lemma \ref{lemma:Kexistencerad}, we see that $((\rho_{\ast},+\infty)\times\mathbb {S}^{m-1},\sigma_S)$ is isometric to
	\[
		(R_0,+\infty)\times\mathbb{S}^{m-1} \quad \text{with metric} \quad \di s^2 + g(s)^2\metric_{\mathbb{S}^{m-1}}
	\]
	for $g=\phi^{-1}$. We let $g_0 : [0,+\infty) \to (0,+\infty)$ be such that
	\[
		\begin{cases}
			g_0(s) = s & \text{on } \, \left[ 0, \frac{1}{2}R_0\right), \\
			g_0(s) = g(s) & \text{on } \, [ R_0, +\infty ).
		\end{cases}
	\]
	Then, $((0,+\infty)\times\mathbb{S}^{m-1},ds^2+g_0(s)^2\metric_{\mathbb{S}^{m-1}})$ has the structure of a radially symmetric complete Riemannian manifold $(\widetilde{M},\metric)$ with its origin $o$ removed, and $(\phi(\rho_{\ast}),+\infty)\times\mathbb{S}^{m-1}$ coincides with $M_K = \widetilde{M} \setminus K$, where $K = \overline{B_{R_0}(o)}$. Recall that (the composition with $\pi$ of) the vertical height function $\tau$ of the hypersurface $\psi$ satisfies
	\[
		\frac{1}{h_0} \, \div\left(\dfrac{h_0(\rho)^2 D\tau}{\sqrt{1-h_0(\rho)^2|D\tau|^2}}\right) = mH \geq mH_0 \qquad \text{on } \, M_K.
	\]
	We choose at will $\rho_2 > \rho_1$ and $\beta$ a real number such that
	\[
		\beta < \max_{\{\rho_2\}\times\mathbb{S}^{m-1}} \tau - \max_{\{\rho_1\}\times\mathbb{S}^{m-1}} \tau.
	\]
	As already observed, we have the existence of a function $u \in C^2([\rho_1,+\infty)\times\mathbb{S}^{m-1})$ satisfying \eqref{Kpotr:sch}. By applying Theorem \ref{thm:Kcomparison}, we obtain that $\sup_M \tau = \sup_{M_K} \tau = +\infty$ and we conclude as in the proof of Theorem \ref{intro:thm:hs2}.
\end{proof}

\section{Further half-space results for mean convex hypersurfaces}

In this section we will prove Theorems \ref{intro:thm:hs3}, \ref{intro:thm:hs4}, \ref{intro:thm:u1} and \ref{intro:thm:u2} of the Introduction. To this aim, we will apply two forms of the maximum principle for the drifted Laplace-Beltrami operator on a complete Riemannian manifold. The first result is a particular case of Theorem 4.1 of \cite{AMR} and we state it as follows.

\begin{prp} \label{AMRwmp}
	Let $(M_0,\metric)$ be a complete, noncompact Riemannian manifold, $f\in C^{\infty}(M_0)$ and let $Q:\RR^+\to\RR^+$ be a nondecreasing function such that
	\begin{equation} \label{wmp:Q}
	\lim_{r\to+\infty}\frac{Q(r)}{r^2}=0.
	\end{equation}
	Suppose that for some (hence, any) point $o\in M_0$,
	\begin{equation} \label{wmp:liminf}
	\liminf_{r\to+\infty}\frac{Q(r)\log\left(\int_{B^0_r}e^{-f}\right)}{r^2}<+\infty,
	\end{equation}
	where $B^0_r=B^{\metric}_r(o)$ is the geodesic ball of $(M_0,\metric)$ centered at $o$ with radius $r$. Given $\kappa\in C^0(\RR)$ and $u\in C^1(M_0)$ such that $u^{\ast}=\sup_{M_0}u<+\infty$, suppose that
	\begin{equation} \label{wmp:ineq}
	e^f\div(e^{-f}Du) \geq \frac{\kappa(u)}{Q\circ\gamma_0}
	\end{equation}
	on $\Omega_c=\{x\in M_0:u(x)>c\}$ for some $c<u^{\ast}$, where $\gamma_0(x)=d_{\metric}(o,x)$ for each $x\in M_0$. Then $\kappa(u^{\ast})\leq0$.
\end{prp}

The proof of Theorem \ref{intro:thm:hs3} is now straightforward.

\begin{proof}[Proof of Theorem \ref{intro:thm:hs3}]
	Consider the manifold $M$ with the metric $\hat\sigma=\pi^{\ast}\sigma$. By assumption, the hyperbolic cosine of the hyperbolic angle is bounded above by $\cosh\theta^{\ast}<+\infty$.
	
	If the conditions in (i) are met, then choose $o\in\pi^{-1}(q)\in M$. There exists a constant $C_1>0$ such that $\hat h=h\circ\pi\leq C_1(1+\distP\circ\pi)\leq C_1(1+\hat\disthat)$ on $M$, where $\distP=d_{\sigma}(q,\;\cdot\;)$, $\hat\disthat=d_{\hat\sigma}(o,\;\cdot\;)$ and the last inequality follows from (\ref{gammaineq}). Arguing as in the proof of Theorem \ref{thm:coshestimate2} and recalling Lemma \ref{lemma:logdiff} we have that
	\[
	\liminf_{r\to+\infty}\frac{\log\left(\int_{\Bhat_r}\hat h\right)}{r}\leq\liminf_{r\to+\infty}\frac{\log\left(\int_0^r\frac{1}{G_0^{m/2}}\sinh(\sqrt{G_0}t)^m \di t\right)}{r}<+\infty,
	\]
	where $\Bhat_r=B^{\hat\sigma}_r(o)$, and therefore
	\begin{align*}
	\liminf_{r\to+\infty}\frac{\log\left(\int_{\Bhat_r}\cosh\theta\cdot\hat h^2\right)}{r} & \leq \liminf_{r\to+\infty}\frac{\log\left(\cosh\theta^{\ast}\cdot\sup_{\Bhat_r}\hat h\right)+\log\left(\int_{\Bhat_r}\hat h\right)}{r} \\
	& \leq \liminf_{r\to+\infty}\frac{\log\left(C_1\cdot\cosh\theta^{\ast}\right)+\log(1+r)+\log\left(\int_{\Bhat_r}\hat h\right)}{r} < +\infty.
	\end{align*}
	
	If the conditions in (ii) are met, then choose again $o\in\pi^{-1}(q)\in M$. As above, we can find $C_2>0$ such that $\hat h\leq C_2(1+\disthat^{\mu})$ on $M$, with $\disthat=d_{\hat\sigma}(o,\;\cdot\;)$. Following the argument used in the proof of Theorem \ref{thm:coshestimate3} to deduce the inequality (\ref{findegcomparison}) we have that, for some $C_3>0$,
	\begin{align*}
	\liminf_{r\to+\infty}\frac{\log\left(\int_{\Bhat_r}\cosh\theta\cdot\hat h^2\right)}{r^{2-\mu}} & \leq \liminf_{r\to+\infty}\frac{\log\left(\cosh\theta^{\ast}\cdot\sup_{\Bhat_r}\hat h\right)+\log\left(\int_{\Bhat_r}\hat h\right)}{r^{2-\mu}} \\
	& \leq \liminf_{r\to+\infty}\frac{C_3+\log(1+r^{\mu})+\log\left(\int_{B^{\sigma}_r(q)}h\right)}{r^{2-\mu}} < +\infty.
	\end{align*}
	
	In both cases, we conclude that conditions (\ref{wmp:Q}) and (\ref{wmp:liminf}) are satisfied on $(M_0,\metric)=(M,\hat\sigma)$ for $f=-\log(\cosh\theta\cdot\hat h^2)$, with $Q$ given by either $Q(r)=C_1\cdot(1+r)$ or $Q(r)=C_2\cdot(1+r^{\mu})$ and such that $\hat h\leq Q\circ\disthat$ on $M$.
	
	Using (\ref{cosh-hatsigma}) and (\ref{Hequation-hatsigma1}) we have that $\tau\in C^{\infty}(M)$ satisfies
	\[
	\frac{1}{\cosh\theta\cdot\hat h^2} \, \div_{\hat\sigma}\left(\cosh\theta\cdot\hat h^2D\tau\right) = \frac{mH}{\cosh\theta\cdot\hat h}.
	\]
	\begin{itemize}
		\item [(a)] Suppose, by contradiction, that $\psi(M)\subseteq\PP\times(-\infty,t_0]$ for some $t_0\in\RR$ and that the mean curvature function $H$ satisfies $H\geq0$ on $M$ and
		\[
		\liminf_{M\ni x\to\infty}H(x)>0.
		\]
		The function $\tau$ then satisfies $\tau^{\ast}=\sup_M\tau < +\infty$. Fix a compact subset $R>0$ such that
		\begin{equation}
		H_0 = \inf_{M\setminus\overline{\Bhat_R}}H>0.
		\end{equation}
		Arguing as in the proof of Theorem \ref{thm:Kcomparison} we see that $\tau_1=\max_{\overline{\Bhat_R}}\tau<\tau^{\ast}$. We have that $\tau$ satisfies
		\[
		\frac{1}{\cosh\theta\cdot\hat h^2} \, \div_{\hat\sigma}\left(\cosh\theta\cdot\hat h^2D\tau\right) = \frac{mH}{\cosh\theta\cdot\hat h}\geq\frac{mH_0}{\cosh\theta^{\ast}\cdot\hat h}\geq\frac{mH_0}{\cosh\theta^{\ast}}\cdot\frac{1}{Q\circ\disthat}
		\]
		on $\Omega_c=\{x\in M:\tau(x)>c\}\subseteq M\setminus\overline{\Bhat_R}$ for any $c\in(\tau_1,\tau^{\ast})$. So, taking $\kappa\equiv mH_0/\cosh\theta^{\ast}$ on $\RR$, we see that $\kappa(\tau^{\ast})>0$ and by applying Proposition \ref{AMRwmp} we obtain the desired contradiction.
		\item [(b)] The argument is the same, up to considering $-\tau$ instead of $\tau$.
	\end{itemize}
\end{proof}

The proof of Theorem \ref{intro:thm:hs4} goes along the same lines of that of Theorem \ref{intro:thm:hs3}.

\begin{proof}[Proof of Theorem \ref{intro:thm:hs4}]
	Consider the manifold $M$ with the metric $g=\psi^{\ast}\gbar$. Since $h$ is bounded on $\PP$, we have that $\hat h=h\circ\pi$ satisfies $\hat h^{\ast}=\sup_M\hat h<+\infty$. Condition (\ref{linfi}) implies that
	\[
	\liminf_{r\to+\infty}\frac{\log\left(\int_{B_r}\hat h^2\right)}{r^2} \leq \liminf_{r\to+\infty}\frac{2\log\hat h^{\ast}+\log(\vol(B_r))}{r^2} < +\infty,
	\]
	so that conditions (\ref{wmp:Q}) and (\ref{wmp:liminf}) are satisfied on $(M_0,\metric)=(M,\hat\sigma)$ for $f=-\log\hat h^2$ and $Q\equiv1$. Using (\ref{Hequation-g2}) we have that $\tau\in C^{\infty}(M)$ satisfies
	\[
	\frac{1}{\hat h^2} \, \div_g\left(\hat h^2\nabla\tau\right)=\frac{mH\cosh\theta}{\hat h}.
	\]
	\begin{itemize}
		\item [(a)] Suppose, by contradiction, that $\psi(M)\subseteq\PP\times(-\infty,t_0]$ for some $t_0\in\RR$ and that the mean curvature function $H$ satisfies $H\geq0$ on $M$ and
		\[
		\liminf_{M\ni x\to\infty}H(x)>0.
		\]
		The function $\tau$ then satisfies $\tau^{\ast}=\sup_M\tau < +\infty$. Fix $R>0$ such that
		\begin{equation} \label{wmp:H0}
		H_0 = \inf_{M\setminus\overline{B_R}}H>0.
		\end{equation}
		As in the proof of Theorem \ref{intro:thm:hs3}, $\tau_1=\max_{\overline{B_R}}\tau<\tau^{\ast}$. We have that $\tau$ satisfies
		\[
		\frac{1}{\hat h^2}\div_g\left(\hat h^2\nabla\tau\right) = \frac{mH\cosh\theta}{\hat h}\geq\frac{mH_0}{\hat h^{\ast}}\frac{1}{Q\circ\dist}
		\]
		on $\Omega_c=\{x\in M:\tau(x)>c\}\subseteq M\setminus\overline{B_R}$ for any $c\in(\tau_1,\tau^{\ast})$, where $\dist=d_g(o,\;\cdot\;)$ on $M$. So, taking $\kappa\equiv mH_0/\hat h^{\ast}$ on $\RR$, we see that $\kappa(\tau^{\ast})>0$ and by applying Proposition \ref{AMRwmp} we obtain the desired contradiction.
		\item [(b)] The argument is the same, up to considering $-\tau$ instead of $\tau$.
	\end{itemize}
\end{proof}

Corollaries \ref{intro:crl:hs3} and \ref{intro:crl:hs4} follow from Theorems \ref{intro:thm:hs3} and \ref{intro:thm:hs4} by reasoning as in the proof of Corollary \ref{intro:crl:hs1-2} given at the end of the previous section.

\bigskip

The second analytical result we rely on is a particular case of Theorem 4.14 of \cite{AMR} which gives a sufficient condition to ensure parabolicity of a drifted Laplace-Beltrami operator on a complete Riemannian manifold.

\begin{prp} \label{AMRparab}
	Let $(M_0,\metric)$ be a complete, noncompact Riemannian manifold, $f\in C^{\infty}(M_0)$ and assume that for some point $o\in M_0$
	\[
	\int_R^{+\infty}\frac{\di r}{\int_{\partial B^0_r}e^{-f}} = +\infty
	\]
	for some (hence, any) $R>0$, where for a.~e.~$r\in\RR^+$ the integral over the boundary of the geodesic ball $B^0_r=B^{\metric}_r(o)$ of $(M_0,\metric)$ centered at $o$ with radius $r$ is intended with respect to the induced $(m-1)$-dimensional Hausdorff measure. Suppose that $u\in C^1(M_0)$ satisfies $u^{\ast}=\sup_{M_0}<+\infty$ and
	\[
	e^f\div(e^{-f}Du)\geq0
	\]
	in the weak sense on $M_0$. Then $u$ is constant.
\end{prp}

\begin{proof}[Proof of Theorem \ref{intro:thm:u1}]
	Consider the manifold $M$ with the metric $g=\psi^{\ast}\gbar$. Since $h$ is bounded on $\PP$,  we have that $\hat h=h\circ\pi$ satisfies $\hat h^{\ast}=\sup_M\hat h<+\infty$. Condition (\ref{notl1}) implies that
	\[
	\int_R^{+\infty}\frac{\di r}{\int_{\partial B_r}\hat h^2} \geq \frac{1}{(\hat h^{\ast})^2}\int_R^{+\infty}\frac{\di r}{\vol(\partial B_r)} = +\infty.
	\]
	\begin{itemize}
		\item [(a)] Suppose that $\psi(M)\subseteq\PP\times(-\infty,t_0]$ for some $t_0\in\RR$ and that $H\geq0$ on $M$. Then by (\ref{Hequation-g2}) we have that $\tau$ satisfies $\tau^{\ast}=\sup_M\tau<+\infty$ and
		\[
		\frac{1}{\hat h^2} \, \div_g(\hat h^2\nabla\tau)\geq0
		\]
		on $M$. Then we apply Proposition \ref{AMRparab} with $(M_0,\metric)=(M,g)$ and $f=-\log\hat h^2$ to conclude that $\tau$ is constant on $M$. Hence, $\psi(M)$ is contained in a slice $\PP\times\{t_1\}$ for some $t_1\in\RR$. Since $(M,g)$ is complete it follows that $\pi:M\to\PP$ is a covering map, so it must be $\psi(M)=\PP\times\{t_1\}$.
		\item [(b)] The argument is the same, considering $-\tau$ instead of $\tau$.
	\end{itemize}
\end{proof}

\begin{proof}[Proof of Theorem \ref{intro:thm:u2}]
	Let $u\in C^{\infty}(\PP)$ be a bounded above solution of equation (\ref{Hequation-sigma}) satisfying the conditions expressed in the statement of the Theorem. In particular, $u$ satisfies
	\[
	e^f\div_{\sigma}(e^{-f}Du)\geq0
	\]
	on $\PP$, with $f=-\log(h^2\cdot(1-h^2|Du|^2)^{-1/2})\in C^{\infty}(\PP)$. Then the desired conclusion follows by applying Proposition \ref{AMRparab}.
\end{proof}

\begin{rmk}
	A Riemannian manifold admitting only constant functions as upper bounded subharmonic functions is said to be parabolic. In \cite{AMR} the definition of parabolicity is extended to a wide family of elliptic differential operators $L$, including the drifted Laplace-Beltrami operator. Such an operator $L$ is said to be parabolic on $(M_0,\metric)$ if each function $u\in C^1(M)$ satisfying $u^{\ast}<+\infty$ and $Lu\geq 0$ on $M$ (in the weak sense) is constant. Therefore, Proposition \ref{AMRparab} gives a sufficient condition for the operator $\mathrm{\Delta}_f=\mathrm{\Delta}-\langle Df,D\;\cdot\;\rangle$ to be parabolic on $(M_0,\metric)$. For more elaborated results on parabolicity of elliptic operators, we refer the interested reader to Chapter 4 of \cite{AMR}. 
	
	A different way of ensuring that an upper bounded function $u\in C^1(M_0)$ satisfying $\mathrm{\Delta}_f u\geq0$ is constant is obtained by replacing assumption (\ref{notl1}) in Theorem \ref{intro:thm:u1} by
	\begin{equation} \label{ulp}
	\int_R^{+\infty}\frac{\di r}{\int_{\partial B^0_r} u^p e^{-f}} = +\infty,
	\end{equation}
	\noindent provided that $u$ is non-negative and $p>1$. The case $p = 1$ requires more care and an extra assumption on the behaviour of $u$ (see for instance Theorem C of \cite{RS}). Note that $u\in L^p(M_0, e^{-f})$ implies (\ref{ulp}). We can also modify assumption (\ref{hnotl1}) in Theorem \ref{intro:thm:u2} in an analogous way to obtain our non-parametric uniqueness result.
	
\end{rmk}

\section{The weak half-space theorem for maximal hypersurfaces}

In this last section of the paper we prove Theorem \ref{intro:thm:coshb} of the Introduction and the subsequent Corollary \ref{intro:crl:u2}. We will need the validity of the Lemma \ref{pseudojacobi} below, whose proof will be postponed in order not to interrupt the logic thread of the main argument.

\begin{proof}[Proof of Theorem \ref{intro:thm:coshb}.]
	We set $u=\tau-t_0$ on $M$. Without loss of generality, we assume that $u>0$ on $M$. Let $o\in M$ be given. We will prove that
	\[
	\cosh\theta(o) = \frac{1}{\sqrt{1-|Du(o)|^2}} \leq e^{(m-1)\sqrt{2G}u(o)}.
	\]
	The proof is divided into five steps.
	
	\textbf{Step 1}. We first introduce some auxiliary functions defined on a geodesic ball centered at $o$ which will be used in the following steps, then we obtain the existence of a point $x_0$ contained in the ball and satisfying a suitable inequality, namely \eqref{Thetaineq} below. Denote by $(\;\cdot\;)_+$ the positive part of a real valued function (that is, set $f_+(x) = \max\{f(x),0\}$ for each $x$ in the domain of $f$). Let $R > 2u(o)$ be given. On the geodesic ball $\Bhat_R=B^{\hat\sigma}_R(o)$ consider the functions
	\begin{equation} \label{defphi}
	\varphi := \left(1-\frac{\disthat^2}{R^2}-Cu\right)_+
	\end{equation}
	and
	\begin{equation} \label{defeta}
	\eta := e^{K\varphi} - 1,
	\end{equation}
	where $\disthat=d_{\hat\sigma}(o,\;\cdot\;)$ and where the constants $C\in(\frac{2}{R},\frac{1}{u(o)})$, $K>0$ are to be chosen later. Note that $\varphi$ has compact support in $\Bhat_R$ since $u>0$, and that $\eta$ has the same support of $\varphi$. Moreover, both functions do not vanish in a neighbourhood of $o$.
	
	For the sake of brevity, we set $\fgrad=\cosh\theta$ on $M$. Note that
	\begin{equation} \label{Theta}
	\fgrad := \frac{1}{\sqrt{1-|Du|^2}}
	\end{equation}
	on $M$. By the previous remarks we have that the function
	\begin{equation} \label{defzeta}
	\zeta := \eta\cdot\fgrad^{\frac{1}{m-1}}
	\end{equation}
	is continuous on $\Bhat_R$, has compact support in $\Bhat_R$ and attains a global maximum at a point $x_0\in \Bhat_R$. Hence we have
	\[
	\fgrad(o) \leq \left(\frac{e^{K\varphi(x_0)}-1}{e^{K\varphi(o)}-1}\right)^{m-1}\fgrad(x_0)
	\]
	and from (\ref{defphi}) we get
	\begin{equation} \label{Thetaineq}
	\begin{split}
	\fgrad(o) & \leq \left(\frac{e^K-1}{e^{-KCu(o)}e^K-1}\right)^{m-1}\fgrad(x_0) \\
	& = \left(\frac{e^K-1}{e^K-e^{KCu(o)}}\right)^{m-1}e^{(m-1)KCu(o)}\fgrad(x_0).
	\end{split}
	\end{equation}
	
	\textbf{Step 2}. In this step we elaborate on \eqref{Thetaineq} to deduce inequality \eqref{Kineq}, that holds true for any $\delta \in (0,1)$ with $f$ defined as in \eqref{deff0}-\eqref{deff}. Inequality \eqref{Kineq} will be used in Step 3 to relate the magnitude of $K$ with that of $|Du|$ at the point $x_0$.
	
	First assume that $\disthat^2$ is smooth in a neighbourhood $\Omega\subseteq \Bhat_R$ of $x_0$. Since $\varphi(x_0) > 0$, we have that $\eta$ is also smooth on $\Omega$. Define the linear second order elliptic differential operator $L$ acting on functions $v\in C^2(M)$ by
	\begin{equation} \label{defL}
	Lv := \Deltahat v + \fgrad^2\Hesshat(v)(Du,Du),
	\end{equation}
	where $\Deltahat$, $\Hesshat$ are the Laplace-Beltrami and the hessian operators of $(M,\hat\sigma)$. Note that the differential of the function $|Du|^2$ is given by
	\[
	\di |Du|^2 = 2\Hesshat(u)(Du,\;\cdot\;),
	\]
	so that
	\[
	\di \fgrad = \fgrad^3\Hesshat(u)(Du,\;\cdot\;).
	\]
	It follows that
	\begin{equation} \label{Luzero}
	Lu = \Deltahat u + \hat\sigma\left(Du,\frac{D\fgrad}{\fgrad}\right) = \frac{1}{\fgrad}\div_{\hat\sigma}(\fgrad Du) = 0
	\end{equation}
	on $M$, since $\div_{\hat\sigma}(\fgrad Du)=0$ as $\psi:M\to\Mbar$ is maximal.
	
	We have the validity of the following Lemma.
	\begin{lemma} \label{pseudojacobi}
		Let $\fgrad$ be as in (\ref{Theta}). Let $\Omega\subseteq M$ be an open subset, $\eta\in C^2(\Omega)$ be nonnegative and $\alpha$ be a constant such that
		\[
		0 < \alpha \leq \frac{1}{m-1}.
		\]
		If the function $\zeta:=\eta\cdot\fgrad^{\alpha}$ satisfies
		\[
		D\zeta(x_0) = 0
		\]
		for some $x_0\in\Omega$ and (\ref{Luzero}) holds on $\Omega$, then
		\begin{equation} \label{Lineq}
		L\zeta \geq \fgrad^{\alpha}\cdot(L\eta + \alpha\eta\fgrad^2\Ric(Du,Du))
		\end{equation}
		at $x_0$.
	\end{lemma}
	
	The hypotheses of the Lemma are met by $\zeta$ defined in (\ref{defzeta}) with $\alpha=\frac{1}{m-1}$, because $x_0$ is an extremal point of $\zeta$. Since $x_0$ is also a maximum point for $\zeta$ and $L$ is elliptic, we have $L\zeta\leq 0$ at $x_0$. Therefore, from (\ref{Lineq}) and (\ref{riccibound}) we deduce
	\begin{equation} \label{Lineq2}
	0 \geq L\eta - G\eta\fgrad^2|Du|^2 = \Deltahat\eta + \fgrad^2\Hesshat(\eta)(Du,Du) - G\eta\fgrad^2|Du|^2
	\end{equation}
	at $x_0$. By (\ref{defeta}) we have
	\[
	\Hesshat(\eta) = Ke^{K\varphi}\Hesshat(\varphi) + K^2e^{K\varphi}\di\varphi\otimes \di\varphi.
	\]
	By (\ref{defphi}) it follows
	\[
	\di\varphi = - \frac{1}{R^2}\di\disthat^2 - C \di u
	\]
	and
	\[
	\Hesshat(\varphi) = - \frac{1}{R^2}\Hesshat(\disthat^2) - C\Hesshat(u),
	\]
	hence we get
	\begin{equation} \label{hesseta}
	\begin{split}
	\Hesshat(\eta) & = Ke^{K\varphi}\cdot\left(- \frac{1}{R^2}\Hesshat(\disthat^2) - C\Hesshat(u)\right. \\
	& \phantom{= Ke^K\cdot\left(\right)} + \frac{K}{R^4}\di\disthat^2\otimes \di\disthat^2 + KC^2 \di u\otimes \di u \\
	& \phantom{= Ke^K\cdot\left(\right)} \left. + \frac{KC}{R^2}\di\disthat^2\otimes \di u + \frac{KC}{R^2} \di u\otimes \di\disthat^2\right).
	\end{split}
	\end{equation}
	Substituting (\ref{hesseta}) into (\ref{Lineq2}) and recalling that $e^{K\varphi} = \eta+1$, it follows
	\begin{align*}
	0 & \geq Ke^{K\varphi}\cdot\left(\frac{K}{R^4}|D\disthat^2|^2 + \frac{K}{R^4}\fgrad^2\hat\sigma(D\disthat^2,Du)^2 + KC^2|Du|^2(1+\fgrad^2|Du|^2) \right. \\
	& \phantom{\geq Ke^{K\varphi}\cdot\left(\right)} + \frac{2KC}{R^2}(1+\fgrad^2|Du|^2)\hat\sigma(D\disthat^2,Du) - \frac{1}{R^2}\Deltahat\disthat^2 \\
	& \phantom{\geq Ke^{K\varphi}\cdot\left(\right)} \left. - \frac{1}{R^2}\fgrad^2\Hesshat(\disthat)(Du,Du) - CLu - \frac{G}{K}\frac{\eta}{\eta+1}\fgrad^2|Du|^2 \right).
	\end{align*}
	Noting that $1+\fgrad^2|Du|^2 = \fgrad^2$ and recalling (\ref{Luzero}), we can further write
	\begin{equation} \label{Lineq3}
	\begin{split}
	0 & \geq \frac{K}{R^4}|D\disthat^2|^2 + \frac{K}{R^4}\fgrad^2\hat\sigma(D\disthat^2,Du)^2 + KC^2\fgrad^2|Du|^2 \\
	& + \frac{2KC}{R^2}\fgrad^2\hat\sigma(D\disthat^2,Du) - \frac{1}{R^2}\Deltahat\disthat^2 - \frac{1}{R^2}\fgrad^2\Hesshat(\disthat)(Du,Du) \\
	& - \frac{G}{K}\frac{\eta}{\eta+1}\fgrad^2|Du|^2.
	\end{split}
	\end{equation}
	
	Using Young inequality
	\begin{align*}
		\delta a^2 + \frac{1}{\delta} b^2 & = \left(\sqrt{\delta}|a| - \frac{|b|}{\sqrt{\delta}}\right)^2 + 2(\sqrt{\delta}|a|)\left(\frac{|b|}{\sqrt{\delta}}\right) \\
		& \geq 2(\sqrt{\delta}|a|)\left(\frac{|b|}{\sqrt{\delta}}\right) = 2 |ab| \geq - 2 ab \qquad \text{for } \, a, b \in \RR, \delta > 0,
	\end{align*}
	we estimate
	\begin{align*}
	\frac{K}{R^4}\fgrad^2\hat\sigma(D\disthat^2,Du)^2 & + KC^2\fgrad^2|Du|^2 + \frac{2KC}{R^2}\fgrad^2\hat\sigma(D\disthat^2,Du) \geq \\
	& \geq (|Du|^2-\delta)KC^2\fgrad^2 + \frac{K}{R^4}\left(1-\frac{1}{\delta}\right)\fgrad^2\hat\sigma(D\disthat^2,Du)^2
	\end{align*}
	for any $\delta\in\RR^+$. Supposing $\delta\in(0,1)$, we can also apply Cauchy-Schwarz inequality
	\[
	- \hat\sigma(D\disthat^2,Du)^2 \geq -|D\disthat^2|^2|Du|^2
	\]
	to obtain
	\begin{equation} \label{Lineq4}
	\begin{split}
	\frac{K}{R^4}\fgrad^2&\hat\sigma(D\disthat^2,Du)^2 + KC^2\fgrad^2|Du|^2 + \frac{2KC}{R^2}\fgrad^2\hat\sigma(D\disthat^2,Du) + \frac{K}{R^4}|D\disthat^2|^2 \geq \\
	& \geq (|Du|^2-\delta)KC^2\fgrad^2 + \frac{K}{R^4}|D\disthat^2|^2\fgrad^2|Du|^2\left(1-\frac{1}{\delta}\right) = \\
	& = KC^2\fgrad^2\cdot\left[\left(1+\frac{|D\disthat^2|^2}{C^2 R^4}\left(1-\frac{1}{\delta}\right)\right)|Du|^2-\delta\right].
	\end{split}
	\end{equation}
	
	By the comparison theorem for the Hessian of the distance function (see Theorem 1.4 of \cite{AMR}), we have that (\ref{sectbound}) implies
	\begin{equation} \label{hessrcomp}
	-\Hesshat(\disthat) \geq -\sqrt{B}\coth(\sqrt{B}\disthat)\cdot\left\{\hat\sigma-\di\disthat\otimes \di\disthat\right\}
	\end{equation}
	on the largest open subset $M_0\subseteq M$ where $\disthat$ is smooth. From (\ref{hessrcomp}) it follows that
	\begin{equation} \label{hessr2comp}
	\begin{split}
	-\Hesshat(\disthat^2) & \geq -2f_0(\disthat)\hat\sigma + 2[f_0(\disthat)-1]\di\disthat \otimes \di\disthat \\
	& = -2f_0(\disthat)\hat\sigma + \frac{f_0(\disthat)-1}{2\disthat^2}\di\disthat^2\otimes \di\disthat^2.
	\end{split}
	\end{equation}
	on $M_0$, where $f:\RR^+\to\RR^+$ is given by
	\begin{equation} \label{deff0}
	f_0(t) = \sqrt{B}t\coth(\sqrt{B}t)
	\end{equation}
	for $t\in\RR^+$.
	
	The function $\disthat^2$ is smooth on $M_0\cup\{o\}$, so the smooth tensor fields $\Hesshat(\disthat^2)$ and $D\disthat^2$ are defined on $M_0\cup\{o\}$. The function $f:\RR\to\RR$ given by
	\begin{equation} \label{deff}
	f(t) := \begin{cases} f_0(|t|) & \text{if } t\neq0, \\ 1 & \text{if } t=0 \end{cases}
	\end{equation}
	is smooth and even. The function $g:\RR\to\RR$ defined by
	\begin{equation} \label{defg}
	g(t) := \begin{cases} \frac{f(t)-1}{2t^2} & \text{if } t\neq 0, \\ \frac{1}{6} & \text{if } t=0 \end{cases}
	\end{equation}
	is also smooth and even. Therefore, $f(\disthat)$ and $g(\disthat)$ are smooth functions on $M_0\cup\{o\}$ and we can extend inequality (\ref{hessr2comp}) obtaining
	\begin{equation} \label{hessr2comp1}
	-\Hesshat(\disthat^2) \geq -2f(\disthat)\hat\sigma + g(\disthat)\di\disthat^2\otimes \di\disthat^2
	\end{equation}
	on $M_0\cup\{o\}$.
	
	From (\ref{hessr2comp1}) and $|D\disthat^2|=2\disthat$ we easily get
	\begin{equation} \label{lapr2comp}
	-\frac{1}{R^2}\Deltahat\disthat^2 \geq -\frac{2mf(\disthat)}{R^2} + \frac{4g(\disthat)\disthat^2}{R^2} = -2\frac{1+(m-1)f(\disthat)}{R^2}
	\end{equation}
	and
	\begin{equation} \label{hessr2comp2}
	\begin{split}
	-\frac{\fgrad^2}{R^2}\Hesshat(\disthat^2)(Du,Du) & \geq -\frac{2f(\disthat)}{R^2}\fgrad^2|Du|^2 + \frac{g(\disthat)}{R^2}\fgrad^2\hat\sigma(D\disthat^2,Du)^2 \\
	& \geq -\frac{2f(\disthat)}{R^2}\fgrad^2|Du|^2 - \frac{\fgrad^2}{2R^2\disthat^2}\hat\sigma(D\disthat^2,Du)^2 \\
	& \geq -2\frac{1+f(\disthat)}{R^2}\fgrad^2|Du|^2.
	\end{split}
	\end{equation}
	
	Using $|D\disthat^2|^2 = 4\disthat^2$, recalling that $f$ is nondecreasing on $\RR^+_0$ and that $\disthat<R$ at $x_0$, we can put (\ref{Lineq4}), (\ref{lapr2comp}) and (\ref{hessr2comp2}) into (\ref{Lineq3}) and it follows that $K$ must satisfy the inequality
	\begin{equation} \label{Kineq}
	\begin{split}
	0 & \geq C^2\fgrad^2\cdot\left[\left(1+\frac{4}{C^2R^2}\left(1-\frac{1}{\delta}\right)\right)|Du|^2-\delta\right]K^2 \\
	& -\frac{2}{R^2}\left[(1+f(R))\fgrad^2|Du|^2 + 1 + (m-1)f(R)\right]K \\
	& -G\fgrad^2|Du|^2.
	\end{split}
	\end{equation}
	
	\textbf{Step 3}. We recall that $R$ was assumed to satisfy $R>2u(o)$, which allowed us to suppose that $C$ satisfies
	\begin{equation} \label{CRineq}
	\frac{2}{R} < C < \frac{1}{u(o)}.
	\end{equation}
	Then, we can set
	\begin{equation} \label{defgamma}
	\gamma := \frac{CR}{2} \in (1,+\infty)
	\end{equation}
	and
	\[
	\delta := \frac{2}{1+\gamma^2} \in (0,1).
	\]
	In this Step we obtain an upper bound on $K$ under the assumption that
	\begin{equation} \label{largegradu}
	|Du|^2 > \frac{4\gamma^2}{(1+\gamma^2)^2} = \left(\frac{2\gamma}{1+\gamma^2}\right)^2 \in (0,1)
	\end{equation}
	at $x_0$. This bound will prove to be crucial to conclude the proof of the theorem in Step 5, arguing by contradiction. Suppose that (\ref{largegradu}) holds and let $\varepsilon\in\RR^+$ be the real number such that
	\begin{equation} \label{gamma1}
	|Du|^2 = \frac{4(1+\varepsilon)\gamma^2}{(1+\gamma^2)^2}
	\end{equation}
	at $x_0$. Then we have
	\begin{equation} \label{gamma2}
	\begin{split}
	\left(1+\frac{4}{C^2R^2}\left(1-\frac{1}{\delta}\right)\right)|Du|^2-\delta & = \left(1+\frac{1-\gamma^2}{2\gamma^2}\right)\cdot\frac{4(1+\varepsilon)\gamma^2}{(1+\gamma^2)^2} - \frac{2}{1+\gamma^2} \\
	& = \frac{2\varepsilon}{1+\gamma^2}
	\end{split}
	\end{equation}
	at $x_0$. Multiplying (\ref{Kineq}) by $z^{-2}$, which satisfies
	\[
	\fgrad^{-2} < \left(\frac{1-\gamma^2}{1+\gamma^2}\right)^2
	\]
	at $x_0$, and using (\ref{gamma1}) and (\ref{gamma2}) we obtain the inequality
	\begin{align*}
	0 & \geq \frac{2\varepsilon\gamma}{1+\gamma^2}C^2K^2 \\
	& - \left[\frac{8(1+\varepsilon)\gamma^2}{(1+\gamma^2)^2}\frac{1+f(R)}{R^2}+2\left(\frac{1-\gamma^2}{1+\gamma^2}\right)^2\frac{1+(m-1)f(R)}{R^2}\right]K \\
	& -\frac{4(1+\varepsilon)\gamma^2}{(1+\gamma^2)^2}G,
	\end{align*}
	which can be rephrased as
	\begin{equation} \label{Kineq2}
	K^2 + pK + q \leq 0
	\end{equation}
	with
	\begin{equation} \label{defpq}
	\begin{split}
	p & := -\frac{1}{C^2}\cdot\left[\frac{(1-\gamma^2)^2}{\varepsilon(1+\gamma^2)}\frac{1+(m-1)f(R)}{R^2}+\frac{4(1+\varepsilon)\gamma^2}{\varepsilon(1+\gamma^2)}\frac{1+f(R)}{R^2}\right], \\
	q & := -\frac{2(1+\varepsilon)\gamma^2}{\varepsilon(1+\gamma^2)}\frac{G}{C^2}.
	\end{split}
	\end{equation}
	As a consequence of (\ref{Kineq2}), we must have
	\begin{equation} \label{Kineq3}
	K \leq \frac{1}{2}(-\,p+\sqrt{p^2-4q}).
	\end{equation}
	
	\textbf{Step 4}. Let $K$ and $C$ as in Step 1 be fixed. Suppose that the fundamental assumption of Steps 2-3 does not hold, that is, the function $\disthat^2$ is not smooth at $x_0$. Then, clearly $x_0\neq 0$. Following the argument of Theorem 4.1 and Lemma 4.3 of \cite{RSS} it can be shown that there exists a unique unit speed minimizing geodesic $\mu:[0,\disthat(x_0)]\to M$ starting from $o$ and ending in $x_0$ and that, for each $\varepsilon'>0$ sufficiently small, the distance function $\disthat_{\varepsilon'}=d_{\hat\sigma}(o_{\varepsilon'},\;\cdot\;)$ from $o_{\varepsilon'} := \mu(\varepsilon')$ is smooth in a neighbourhood of $x_0$, that is, we are applying Calabi's trick. Moreover, the point $x_0$ is a local maximum for the function
	\[
	\zeta_{\varepsilon'} := (e^{K\varphi_{\varepsilon'}}-1)\cdot\fgrad^{\frac{1}{m-1}},
	\]
	where $\varphi_{\varepsilon'}$ is defined by
	\[
	\varphi_{\varepsilon'} := \left(1-\frac{(\disthat_{\varepsilon'}+\varepsilon')^2}{R^2}-Cu\right)_+.
	\]
	Therefore, the arguments of Steps 2 and 3 can be repeated for each $\varepsilon'>0$ small enough, and by letting $\varepsilon'\to0$ we still obtain that (\ref{largegradu}) can be satisfied at $x_0$ with $\gamma$ defined as in (\ref{defgamma}) only if $K$ satisfies (\ref{Kineq3}) for $p$ and $q$ as in (\ref{defpq}).
	
	\textbf{Step 5}. We are now ready to prove the validity of inequality
	\[
	z(o) \leq e^{(m-1)\sqrt{2G}u(o)}.
	\]
	Let $\{R_n\}_n$ be a nondecreasing positive sequence such that
	\begin{equation} \label{limRn}
	\lim_{n\to+\infty} R_n = +\infty.
	\end{equation}
	Let $\{C_n\}_n$ be a nonincreasing positive sequence such that
	\begin{equation} \label{asCn}
	\frac{1}{u(o)} > C_n \sim \frac{2\log R_n}{R_n}
	\end{equation}
	as $n\to+\infty$. Then
	\begin{equation} \label{limCn}
	\lim_{n\to+\infty} C_n = 0.
	\end{equation}
	For each $n\geq 1$, set
	\[
	\gamma_n := \frac{C_nR_n}{2}.
	\]
	Note that, for each $n$ large enough, conditions (\ref{CRineq}) are satisfied by $C=C_n$, $R=R_n$ and $\gamma$ defined as in (\ref{defgamma}) coincides with $\gamma_n$. Moreover,
	\begin{equation} \label{asgamma}
	\gamma_n \sim \log R_n \to +\infty
	\end{equation}
	as $n\to+\infty$.
	
	Let $\beta>1$ be given and let $\{K_n\}_n$ be a positive sequence such that
	\begin{equation} \label{condKnCn1}
	K_n > \beta\frac{\sqrt{2G}}{C_n}
	\end{equation}
	for each $n\geq 1$ and
	\begin{equation} \label{condKnCn2}
	\lim_{n\to+\infty} (K_nC_n) = \beta\sqrt{2G}.
	\end{equation}
	For $n\geq 1$ let $\varphi_n$, $\eta_n$ and $\zeta_n$ be the functions defined on the geodesic ball
	\[
	\Bhat_{R_n} = \{x\in M : \disthat(x) < R_n\}
	\]
	respectively by
	\begin{align*}
	\varphi_n & := \left(1-\frac{\disthat^2}{R_n^2} - C_n u\right)_+, \\
	\eta_n & := e^{K_n\varphi_n} - 1, \\
	\zeta_n & := \eta_n\cdot\fgrad^{\frac{1}{m-1}},
	\end{align*}
	and let $x_n\in \Bhat_{R_n}$ be a global maximum point for $\zeta_n$. As observed in Step 1, we have that $\fgrad(o)$ is bounded from above by each term of the sequence
	\begin{equation} \label{asbound}
	\left(\frac{e^{K_n}-1}{e^{K_n}-e^{K_nC_nu(o)}}\right)e^{(m-1)K_nC_nu(o)}\fgrad(x_n) \sim e^{\beta(m-1)\sqrt{2G}u(o)}\fgrad(x_n)
	\end{equation}
	as $n\to+\infty$, where the asymptotic relation follows from (\ref{limCn}), (\ref{condKnCn1}) and (\ref{condKnCn2}).
	
	We show that
	\begin{equation} \label{Thetato1}
	\lim_{n\to+\infty} \fgrad(x_n) = 1.
	\end{equation}
	Suppose, by contradiction, that (\ref{Thetato1}) is false. Then, there exists a strictly increasing sequence $\{n_k\}_k$ of positive integers such that
	\begin{equation} \label{absurdhp}
	\lim_{k\to+\infty} |Du(x_{n_k})|^2 = \ell \in (0,1].
	\end{equation}
	Up to choosing $n_1$ large enough, we have
	\[
	|Du(x_{n_k})|^2 > \frac{4\gamma_{n_k}^2}{(1+\gamma_{n_k}^2)^2}
	\]
	for each $k$, since, by (\ref{asgamma}),
	\[
	\frac{4\gamma_n^2}{(1+\gamma_n^2)^2} \sim \frac{4}{\log^2 R_n} \to 0
	\]
	as $n\to+\infty$. Defining $\{\varepsilon_k\}_k$ as the sequence of real numbers such that
	\[
	|Du(x_{n_k})|^2 = \frac{4(1+\varepsilon_k)\gamma_{n_k}^2}{(1+\gamma_{n_k}^2)^2}
	\]
	for each $k$, from (\ref{absurdhp}) we have that
	\begin{equation} \label{asepsilon}
	\varepsilon_k \sim \frac{\ell}{4}\gamma_{n_k}^2 \sim \frac{\ell}{4}\log^2 R_{n_k}
	\end{equation}
	as $k\to+\infty$. By Step 3 and Step 4, we deduce that the subsequence $\{K_{n_k}\}$ satisfies
	\begin{equation} \label{Knineq}
	K_{n_k} \leq \frac{1}{2}\left(-\,p_k + \sqrt{p_k^2 - 4q_k}\right)
	\end{equation}
	for each $k$, where $\{p_k\}_k$ and $\{q_k\}_k$ are the sequences of real numbers defined by
	\begin{equation} \label{defpkqk}
	\begin{split}
	p_k & := -\frac{1}{C_{n_k}^2}\cdot\left[\frac{(1-\gamma_{n_k}^2)^2}{\varepsilon_k(1+\gamma_{n_k}^2)}\frac{1+(m-1)f(R_{n_k})}{R_{n_k}^2}+\frac{4(1+\varepsilon_k)\gamma_{n_k}^2}{\varepsilon_{n_k}(1+\gamma_{n_k}^2)}\frac{1+f(R_{n_k})}{R_{n_k}^2}\right], \\
	q_k & := -\frac{2(1+\varepsilon_k)\gamma_{n_k}^2}{\varepsilon_k(1+\gamma_{n_k}^2)}\frac{G}{C_{n_k}^2}.
	\end{split}
	\end{equation}
	By (\ref{deff0}), (\ref{deff}) and (\ref{limRn})
	\[
	\lim_{k\to+\infty}\frac{f(R_{n_k})}{R_{n_k}^2} = \sqrt{B},
	\]
	therefore, using (\ref{asepsilon}) and (\ref{asCn}), we deduce
	\[
	p_k = O\left(\frac{1}{C_{n_k}^2 R_{n_k}}\right) = O\left(\frac{R_{n_k}}{\log^2 R_{n_k}}\right)
	\]
	and
	\[
	q_k \sim -G\frac{R_{n_k}^2}{2\log^2 R_{n_k}}
	\]
	as $k\to+\infty$. Hence,
	\begin{equation} \label{absurdconc1}
	\frac{1}{2}\left(-\,p_k+\sqrt{p_k^2-4q_k}\right) \sim \sqrt{-q_k} \sim \frac{\sqrt{2G}}{C_{n_k}} \qquad \text{as } k\to+\infty.
	\end{equation}
	Putting together (\ref{condKnCn2}), (\ref{Knineq}) and (\ref{absurdconc1}) it follows
	\[
	1 > \frac{1}{\beta} = \lim_{k\to+\infty} \frac{\sqrt{2G}}{K_{n_k}C_{n_k}} = \lim_{k\to+\infty} \frac{1}{2K_{n_k}}\left(-\,p_k+\sqrt{p_k^2-4q_k}\right) \geq 1
	\]
	which gives the desired contradiction. Therefore, (\ref{Thetato1}) is proved and by (\ref{Thetaineq}) and (\ref{asbound}) we get
	\begin{equation}
	\fgrad(o) \leq \lim_{n\to+\infty} e^{\beta(m-1)\sqrt{2G}u(o)}\fgrad(x_n) = e^{\beta(m-1)\sqrt{2G}u(o)}.
	\end{equation}
	Since $\beta>1$ is arbitrarily given, we conclude
	\begin{equation} \label{habound3}
	\fgrad(o) \leq e^{(m-1)\sqrt{2G}u(o)}.
	\end{equation}
\end{proof}

\begin{proof}[Proof of Corollary \ref{intro:crl:u2}]
	By Theorem \ref{intro:thm:coshb}, for each $x\in M$ we have the validity of the inequality
	\[
	\cosh\theta(x)\leq e^{(m-1)\sqrt{2G}|\tau(x)-t_0|}
	\]
	for any $G>0$. Lettin $G\to 0^+$, we obtain $\cosh\theta(x)=1$. Since $x\in M$ is arbitrary, it follows that $\cosh\theta\equiv1$ on $M$. By (\ref{cosh-hatsigma}), this is equivalent to $d\tau=0$ on $M$. Since $M$ is connected, this implies that $\tau$ is constant and therefore $\psi(M)$ is contained in a slice $\PP\times\{t_1\}$ for some $t_1\in\RR$. Clearly, $t_1\neq t_0$. Since $\pi:M\to\PP$ is a covering map, we conclude that $\psi(M)=\PP\times\{t_1\}$.
\end{proof}

\begin{proof}[Proof of Lemma \ref{pseudojacobi}]
	Up to restricting ourselves to a smaller neighbourhood of $x_0$, we can assume that a local orthonormal frame $\{e_i\}_{1\leq i\leq m}$ for $TM$ is defined on $\Omega\subseteq M$. Let $\{\theta^i\}_{1\leq i\leq m}$ be the coframe dual to $\{e_i\}_i$.
	
	We denote by $u_i$, $1\leq i\leq m$, the components of $\di u$, that is, $\di u=u_i\theta^i$, and by $u^i$ the components of the metrically equivalent vector field $Du=u^ie_i$. Note that orthogonality of the frame yields $u^i=u_i$ for $1\leq i\leq m$. Since $\{e_i\}_i$ is orthonormal, we have
	\[
	|Du|^2 = u_iu_i.
	\]
	Taking covariant derivative and recalling that the metric is parallel, we get
	\[
	\di |Du|^2 = 2u_iu_{ik}\theta^k,
	\]
	where, hereafter, we denote the components of the covariant derivative of a given tensor field by adding a lower index to the components of the field. Since $\Hesshat(u) = u_{ij}\theta^i\otimes\theta^j$, the relation above reads as
	\[
	\di |Du|^2 = 2\Hesshat(u)(Du,\;\cdot\;),
	\]
	as claimed in the proof of Theorem \ref{intro:thm:coshb}. Then, we also have
	\begin{equation} \label{covtheta}
	\begin{split}
	\fgrad_i & = \fgrad^3\delta^{kt}u_tu_{ki} = \fgrad^3u^ku_{ki} \\
	\fgrad_{ij} & = \fgrad^3\delta^{kt}u_{tj}u_{ki} + \fgrad^3\delta^{kt}u_tu_{kij} + 3\fgrad^2\delta^{kt}u_tu_{ki}\fgrad_j \\
	& = \fgrad^3\delta^{kt}u_{tj}u_{ki} + \fgrad^3u^ku_{kij} + 3\fgrad^5u^ku^tu_{ki}u_{tj},
	\end{split}
	\end{equation}
	where $\di\fgrad = \fgrad_i\theta^i$ and $\Hesshat(\fgrad) = \fgrad_{ij}\theta^i\otimes\theta^j$. Denoting by $R_{ijkt}$ the components of the Riemann curvature tensor of $(M,\hat\sigma)$ along the frame $\{\theta^i\}_i$ and using the identities
	\[
	u_{ij} = u_{ji}, \qquad u_{kij} = u_{ikj} = u_{ijk} + \delta^{ls}u_lR_{sikj}
	\]
	(see (1.110) and (1.116) of \cite{AMR}) we can also write
	\begin{equation} \label{hesstheta}
	\begin{split}
	\fgrad_{ij} & = \fgrad^3\delta^{kt}u_{ik}u_{jt} + \fgrad^3u^ku_{ijk} + \fgrad^3u^ku^tR_{kitj} + 3\fgrad^5u^ku^tu_{ik}u_{jt}.
	\end{split}
	\end{equation}
	Noting that the action of the operator $L$ on a generic function $v\in C^2(\Omega)$ is given by
	\begin{equation} \label{esprL}
	Lv = a^{ij}v_{ij}
	\end{equation}
	with $\Hesshat(v)=v_{ij}\theta^i\otimes\theta^j$ and $a^{ij}e_i\otimes e_j$ the tensor field with components
	\begin{equation} \label{defa}
	a^{ij} = \delta^{ij} + \fgrad^2u^iu^j,
	\end{equation}
	we obtain from (\ref{hesstheta})
	\begin{equation} \label{LTheta1}
	\begin{split}
	L\fgrad & = a^{ij}\fgrad_{ij} = \fgrad^3\delta^{ij}\delta^{kt}u_{ik}u_{jt} + 4\fgrad^5u^iu^j\delta^{kt}u_{ik}u_{jt} \\
	& + 3\fgrad^7u^iu^ju^ku^tu_{ik}u_{jt} + \fgrad^3u^k\delta^{ij}u_{ijk} + \fgrad^5u^iu^ju^ku_{ijk} \\
	& + \fgrad^3u^iu^jR_{ij},
	\end{split}
	\end{equation}
	with $\Richat = R_{ij}\theta^i\otimes\theta^j$ the Ricci curvature tensor of $(M,\hat\sigma)$.
	
	Observe that
	\begin{equation} \label{leib1}
	\delta^{ij}u_{ijk} = (\delta^{ij}u_{ij})_k - \delta^{ij}_ku_{ij} = (\delta^{ij}u_{ij})_k,
	\end{equation}
	since $\delta^{ij}e_i\otimes e_j$ is parallel, and that
	\begin{equation} \label{leib2}
	\fgrad^2u^iu^ju_{ijk} = (\fgrad^2u^iu^ju_{ij})_k - 2\fgrad^4u^iu^ju^tu_{ij}u_{kt} - 2\fgrad^2\delta^{it}u_{tk}u^ju_{ij}.
	\end{equation}
	Putting together (\ref{leib1}) and (\ref{leib2}) and recalling (\ref{esprL}) and (\ref{defa}) we get
	\begin{align*}
	\fgrad^3u^k\delta^{ij}u_{ijk} & + \fgrad^5u^iu^ju^ku_{ijk} = \fgrad^3(Lu)_k \\
	& - 2\fgrad^3u^iu^ju^ku^tu_{ij}u_{kt} - 2\fgrad^5u^ku^j\delta^{it}u_{tk}u_{ij}.
	\end{align*}
	Since $Lu=0$ on $\Omega$ by assumption, we can put this identity into (\ref{LTheta1}) to get
	\begin{equation} \label{LTheta2}
	L\fgrad = \fgrad^3\delta^{ij}\delta^{kt}u_{ik}u_{jt} + 2\fgrad^5u^iu^j\delta^{kt}u_{ik}u_{jt} + \fgrad^7u^iu^ju^ku^tu_{ik}u_{jt} + \fgrad^3u^iu^jR_{ij}.
	\end{equation}
	Using (\ref{defa}), expression (\ref{LTheta2}) can be rewritten as
	\begin{equation} \label{LTheta3}
	L\fgrad = \fgrad^3\cdot(a^{ij}a^{kt}u_{ik}u_{jt}+u^iu^jR_{ij}).
	\end{equation}
	
	Next,
	\begin{equation} \label{covzeta}
	\di\zeta = \fgrad^{\alpha}\left(\di\eta + \alpha\frac{\eta}{\fgrad}\di\fgrad\right)
	\end{equation}
	and having assumed $d\zeta=0$ at $x_0$, from (\ref{covzeta}) we deduce
	\begin{equation} \label{deta}
	\di\eta = - \alpha\frac{\eta}{\fgrad}\di\fgrad
	\end{equation}
	at $x_0$. Since
	\[
	\Hesshat(\zeta) = \fgrad^{\alpha}\Hesshat(\eta) + \di\fgrad^{\alpha}\otimes \di\eta + \di\eta\otimes \di\fgrad^{\alpha} + \eta\Hesshat(\fgrad^{\alpha}),
	\]
	from (\ref{deta}) and $d\fgrad^{\alpha} = \alpha\fgrad^{\alpha-1}d\fgrad$ we obtain
	\begin{equation} \label{hesszeta1}
	\Hesshat(\zeta) = \fgrad^{\alpha}\Hesshat(\eta) - 2\alpha^2\eta\fgrad^{\alpha-2}\di\fgrad\otimes \di\fgrad + \eta\Hesshat(\fgrad^{\alpha}).
	\end{equation}
	On the other hand,
	\[
	\Hesshat(\fgrad^{\alpha}) = \alpha\fgrad^{\alpha-1}\Hesshat(\fgrad) + \alpha(\alpha-1)\fgrad^{\alpha-2}\di\fgrad\otimes \di\fgrad
	\]
	and inserting into (\ref{hesszeta1}) we obtain
	\begin{equation} \label{hesszeta2}
	\Hesshat(\zeta) = \fgrad^{\alpha}\cdot\left(\Hesshat(\eta) + \frac{\alpha\eta}{\fgrad}\Hesshat(\fgrad) - \frac{\alpha(\alpha+1)\eta}{\fgrad^2}\di\fgrad\otimes \di\fgrad\right).
	\end{equation}
	By the definition of $L$ we get
	\begin{equation} \label{Lzeta1}
	L\zeta = \fgrad^{\alpha}\cdot\left[L\eta + \alpha\eta\left(\frac{L\fgrad}{\fgrad} - \frac{\alpha+1}{\fgrad^2}a^{ij}\fgrad_i\fgrad_j\right)\right].
	\end{equation}
	
	From (\ref{covtheta}) it follows that
	\begin{equation} \label{adtheta}
	a^{ij}\fgrad_i\fgrad_j = \fgrad^6a^{ij}u^ku^tu_{ik}u_{jt}.
	\end{equation}
	Putting (\ref{adtheta}) and (\ref{LTheta3}) into (\ref{Lzeta1}) we get
	\begin{equation}
	L\zeta = \fgrad^{\alpha}\cdot[L\eta + \alpha\eta(\fgrad^2 a^{ij}a^{kt}u_{ik}u_{jt} - (\alpha+1)\fgrad^4 a^{ij}u^ku^tu_{ik}u_{jt} + \fgrad^2u^iu^jR_{ij})]
	\end{equation}
	Since $\alpha$, $\eta$ and $\fgrad$ are nonnegative, the claim of the Lemma will follow by showing that
	\[
	a^{ij}a^{kt}u_{ik}u_{jt} - (\alpha+1)\fgrad^2a^{ij}u^ku^tu_{ik}u_{jt}
	\]
	is nonnegative.
		
	We define the tensor field $B=B^i_j\theta^j\otimes e_i$ by setting
	\begin{equation} \label{defM}
	B^i_j := a^{ik}u_{kj}
	\end{equation}
	for $1\leq i,j\leq m$. We also define the tensor field $a_{ij}\theta^i\otimes\theta^j$ by setting
	\begin{equation} \label{definva}
	a_{ij} := \delta_{ij}-u_iu_j
	\end{equation}
	for $1\leq i,j\leq m$. Observe that $a_{ij}\theta^i\otimes\theta^j=\psi^{\ast}\gbar$ on $\Omega$, hence it is a Riemannian metric on $\Omega$ and $(a_{ij})$ is the inverse of the matrix $(a^{ij})$, that is,
	\begin{equation} \label{inva}
	a^{ij}a_{jk} = \delta^i_k.
	\end{equation}
	Using (\ref{defM}) and (\ref{inva}), we have
	\begin{equation} \label{esprM1}
	\begin{split}
	a^{ij}a^{kt}u_{ik}u_{jt} - (\alpha+1)&\fgrad^2a^{ij}u^ku^tu_{ik}u_{jt} = \\
	& = B^j_kB^k_j - (\alpha+1)\fgrad^2B^j_ku^k\delta^l_ju_{lt}u^t \\
	& = B^i_jB^j_i - (\alpha+1)\fgrad^2B^j_ku^ka_{ij}a^{il}u_{lt}u^t \\
	& = B^i_jB^j_i - (\alpha+1)\fgrad^2a_{ij}B^j_ku^kB^i_tu^t.
	\end{split}
	\end{equation}
	The tensor field $B$ induces, at each point $x\in\Omega$, a linear operator
	\begin{align*}
	B(x) : T_xM & \to T_xM \\
	X & \mapsto B^i_jX^je_i =: B(x)X \qquad \forall X=X^je_j \in T_xM
	\end{align*}
	which is self-adjoint with respect to the scalar product $a_{ij}(x)\theta^i|_x\otimes\theta^j|_x$. Indeed, recalling (\ref{defM}), (\ref{inva}) we have
	\begin{align*}
	(a_{ij}\theta^i\otimes\theta^j)(BX,Y) & = a_{ij}B^i_tX^tY^j = a_{ij}a^{ik}u_{kt}X^tY^j \\
	& = a^{ki}a_{ij}u_{kt}X^tY^j = \delta^k_ju_{kt}X^tY^j = \Hesshat(u)(X,Y)
	\end{align*}
	for each $X=X^ie_i,Y=Y^ie_i\in T_xM$. We let $\{\lambda_i\}_{1\leq i\leq m}$ be the eigenvalues of $B$.
	
	We are now ready to show that
	\[
	B^i_jB^j_i - (\alpha+1)\fgrad^2a_{ij}B^j_ku^kB^i_tu^t
	\]
	is nonnegative and by (\ref{esprM1}) this will conclude the proof. Without loss of generality, we suppose that $(\lambda_1)^2 \geq (\lambda_i)^2$ for $2\leq i\leq m$. Let $B^2=B\circ B$. We have
	\begin{equation} \label{esprM3}
	B^i_jB^j_i = \trace(B^2) = \sum_{i=1}^m (\lambda_i)^2.
	\end{equation}
	On the other hand, we have
	\[
	a_{ij}B^j_ku^kB^i_tu^t = (a_{ij}\theta^i\otimes\theta^j)(BDu,BDu) = (a_{ij}\theta^i\otimes\theta^j)(B^2Du,Du).
	\]
	Since $\{(\lambda_i)^2\}_{1\leq i\leq m}$ are the eigenvalues of $B^2$, we have
	\begin{align*}
	(a_{ij}\theta^i\otimes\theta^j)(B^2Du,Du) & \leq (\lambda_1)^2(a_{ij}\theta^i\otimes\theta^j)(Du,Du) \\
	& = (\lambda_1)^2(u^iu_i-u_iu_ju^iu^j) = (\lambda_1)^2|Du|^2(1-|Du|^2) \\
	& = (\lambda_1)^2\frac{|Du|^2}{\fgrad^2}
	\end{align*}
	and therefore, using the fact that $|Du|<1$,
	\begin{equation} \label{esprM4}
	-(\alpha+1)\fgrad^2a_{ij}B^j_ku^kB^i_tu^t \geq -(\alpha+1)(\lambda_1)^2|Du|^2 \geq -(\alpha+1)(\lambda_1)^2.
	\end{equation}
	Now, recalling that $Lu=0$ on $\Omega$ we have
	\[
	0 = a^{ij}u_{ij} = \delta^k_ia^{ij}u_{jk} = \delta^k_iB^i_k = \trace(B) = \sum_{i=1}^m \lambda_i.
	\]
	The triangular inequality then implies
	\[
	|\lambda_1| = \left|0-\sum_{i=2}^m\lambda_i\right| \leq \sum_{i=1}^m |\lambda_i|
	\]
	and from Newton's inequality we get
	\begin{equation} \label{newton}
	(\lambda_1)^2 \leq \left(\sum_{i=2}^m |\lambda_i|\right)^2 \leq (m-1)\sum_{i=2}^m(\lambda_i)^2.
	\end{equation}
	Recalling that by assumption $\alpha\leq\frac{1}{m-1}$, using (\ref{esprM3}), (\ref{esprM4}) and (\ref{newton}) we finally obtain
	\begin{align*}
	B^i_jB^j_i - (\alpha+1)\fgrad^2a_{ij}B^j_ku^kB^i_tu^t & \geq \sum_{i=1}^m(\lambda_i)^2 - (\alpha+1)(\lambda_1)^2 \\
	& = \sum_{i=2}^m(\lambda_i)^2 - \alpha(\lambda_1)^2 \\
	& \geq \sum_{i=2}^m(\lambda_i)^2 - \frac{1}{m-1}(\lambda_1)^2 \geq 0
	\end{align*}
	and the Lemma is proved.
\end{proof}

\section*{Acknowledgements}

The second author is partially supported by Spanish MINECO and ERDF project MTM2016-78807-C2-1-P.

\end{document}